\documentclass[10pt,leqno]{amsart}
\usepackage{graphicx}
\usepackage{indentfirst,csquotes}

\usepackage{amssymb,amsthm,amsmath}
\usepackage{xcolor,paralist,hyperref,fancyhdr,etoolbox}

\hypersetup{ colorlinks=true, linkcolor=black, filecolor=black, urlcolor=black }
\bibliography{ref}

\usepackage{lipsum}
\usepackage{graphicx}
\usepackage[colorinlistoftodos]{todonotes}
\usepackage{amssymb}
\usepackage{mathtools}
\usepackage{bm}
\usepackage{mathrsfs}
\usepackage{indentfirst} 
\usepackage{amsthm}
\usepackage{float}
\usepackage{amscd}
\usepackage{hyperref}
\usepackage{enumitem} 
\usepackage{tikz-cd}
\newcommand{\strat}{\mathrm{strat}}
\newcommand{\dom}{\mathrm{Dom}}
\newcommand{\eps}{\epsilon}

\newcommand{\sing}{\mathrm{sing}}

\newcommand{\harmonic}{\mathcal{H}}
\setlist[enumerate]{label=\arabic*.}

\hypersetup{colorlinks=false, linkbordercolor=red}

\begin{document}
\title{Derived Stratifications and Arithmetic Intersection Theory for Varieties with Isolated Singularities} 
\author[Jiaming Luo\and Shirong Li]{Jiaming Luo\and Shirong Li}
\date{\today}
\address{School of Mathematics and Statistics, Henan University of Science and Technology, 263, Kaiyuan Avenue, Luoyang, China.}
\email{luojiaming@hhu.edu.cn}
\maketitle

\let\thefootnote\relax

\begin{abstract}
In this paper, We develop the stratified de Rham theory on singular spaces using modern tools including derived geometry and stratified structures. This work unifies and extends  the de Rham theory, Hodge theory, and deformation theory of singular spaces into the frameworks of stratified geometry, $p$-adic geometry, and derived geometry. Additionally, we close a gap in Ohsawa's original proof, concerning the convergence of $L^2$ harmonic forms in the Cheeger–Goresky–MacPherson conjecture for varieties with isolated singularities. Indicating that harmonic forms converge strongly and the $L^2$-cohomology coincides with intersection cohomology. 
\end{abstract} 

\bigskip

\section*{ Contents}

\begin{enumerate}
  \item Introduction \dotfill \hyperlink{INTRODUCTION}{1}
  \item Notation and Acknowledgments \dotfill \hyperlink{NOTATION AND ACKNOWLEDGMENTS}{2}
  \item Derived Stratified Schemes and Higher de Rham Complexes \dotfill \hyperlink{DERIVED STRATIFIED SCHEMES AND HIGHER DE RHAM COMPLEXES}{3}
  \item Logarithmic Stratified de Rham Complex and Mixed Hodge Theory \dotfill \hyperlink{LOGARITHMIC STRATIFIED DE RHAM COMPLEX AND MIXED HODGE THEORY}{9}
  \item Derived Category Interpretation of Intermediate Perversity \dotfill \hyperlink{DERIVED CATEGORY INTERPRETATION OF INTERMEDIATE PERVERSITY}{23} 
  \item Logarithmic Convexity, Singularity Stability and Deformation Theory \dotfill \hyperlink{LOGARITHMIC CONVEXITY, SINGULARITY STABILITY AND DEFORMATION THEORY}{27}
  \item Stratified Elliptic Regularity and Derived Duality for Ohsawa’s Theorem 7 \dotfill \hyperlink{STRATIFIED ELLIPTIC REGULARITY AND DERIVED DUALITY FOR OHSAWA'S THEOREM 7}{52}
  \item Future Work \dotfill \hyperlink{FUTURE WORK}{60}
\end{enumerate}

\hypertarget{INTRODUCTION}{}
\section{INTRODUCTION}
Building on Ohsawa's work on the Cheeger-Goreski-MacPherson conjecture for varieties with isolated singularities (\cite{Ohs91}), this paper systematically integrates modern mathematical tools like derived geometry, logarithmic structures, and p-adic analysis, proposing a series of theoretical extensions. We construct the derived stratified de Rham complex, prove its degeneration under derived transversality conditions (Proposition 3.8), and establish a connection between logarithmic stratified structures and mixed Hodge theory (Proposition 4.1). By reinterpreting moderate permeability as a derived Lagrangian subcategory, we reveal its deep connection to mirror symmetry (Corollary 5.3). Furthermore, we develop a stratified deformation theory (Proposition 6.3, Proposition 6.17) and a p-adic intersection comparison (Proposition 4.4), providing a unified framework for higher-dimensional singularities and arithmetic geometry. Finally, we can bridge $L^2$-estimates for the de Rham complex, Micro-local elliptic regularity at singular strata and derived categorical duality, this work not only corrects the Ohsawa's original proof gaps but also opens new connections to the geometric Langlands program and string theory.

\hypertarget{NOTATION AND ACKNOWLEDGMENTS}{}
\section{NOTATION AND ACKNOWLEDGMENTS}
\begin{center}
    \textit{2.1 Notation}
\end{center}

The following directory is designed to serve as a concise reference for the specialized notation used throughout the paper, aligned with standard conventions in derived algebraic geometry, logarithmic Hodge theory, and p-adic arithmetic geometry.

\noindent\textbf{(1) Derived stratified schemes and stratified de Rham complexes:}
\begin{itemize}
    \item $X$: Noetherian derived scheme.  
    \item $X_p$: The $p$-th stratum.  
    \item $\mathscr{G}^p$: Quasi-coherent differential graded (qDG) algebra.  
    \item $\mathbb{L}_{X_p \times_X X_q / X_q}$: Relative cotangent complex.
    \item $\mathbb{L}\Omega^{\bullet,\text{strat}}_X$: Higher stratified de Rham complex.  
    \item $\Omega^{\bullet,\text{strat}}_X$: Classical stratified de Rham complex.  
    \item $\mathbf{IC}_p^{\bullet}$: Intersection cohomology complex on stratum $X_p$.  
    \item $\Omega_X^{\bullet,\log\text{-strat}}$: Logarithmic stratified de Rham complex.  
    \item $\Omega_{X_p}^{\bullet}(\log D_p)$: Logarithmic de Rham complex on $X_p$.
\end{itemize}
\textbf{(2) Derived operators and cohomology groups:}
\begin{itemize}
    \item $\mathrm{holim}$: Homotopy limit.  
    \item $\mathbb{R}\underline{\mathrm{Hom}}$: Derived internal Hom functor.  
    \item $\boxtimes^L$: Derived external tensor product.  
    \item $\otimes^L$: Derived tensor product.  
    \item $\mathbb{D}$: Verdier duality functor. 
    \item $\mathbb{H}^k(X, -)$: Hypercohomology.  
    \item $IH^k(X)$: Intersection cohomology.  
    \item $H_{\mathrm{rig}}^k$: Rigid cohomology.  
    \item $\mathbb{H}_{\mathrm{abs}}^k$: Absolute prismatic cohomology.  
    \item $\mathrm{Gr}_m^W$: Weight $m$ associated graded.  
\end{itemize}
\textbf{(3) $p$-adic geometry and deformation theory:}
\begin{itemize}
    \item $X^{\mathrm{rig}}$: Perfectoid rigidification.  
    \item $\mathcal{X}_{\mathrm{FF}}$: Fargues-Fontaine curve.
    \item $\mathscr{G}^{p,\mathrm{rig}}$: Perfectoid coefficient sheaf.  
    \item $\Omega_{X^{\mathrm{rig}}}^{\bullet,\mathrm{strat}}$: Stratified rigid de Rham complex.
    \item $\mathrm{Def}^{\mathrm{strat}}_X$: Stratified deformation complex.
    \item $\mathbb{T}^{\mathrm{strat}}_X$: Stratified tangent sheaf.  
    \item $\mathcal{T}_{X_p}(-\log D_p)$: Logarithmic tangent sheaf.  
    \item $\mathbb{L}\mathrm{Def}^{\mathrm{strat}}_X$: Derived stratified deformation complex.  
    \item $\kappa$: Kuranishi obstruction map. 
\end{itemize}
\textbf{(4) Category theory and other key notations:}
\begin{itemize}
    \item $\mathrm{Fuk}(X^\vee)$: Fukaya category of mirror $X^\vee$.  
    \item $HF^{\bullet}(L_p, \mathscr{W}_{Y_p})$: Floer cohomology.  
    \item $\mathscr{E}st_{\mathscr{O}_S}^{\mathrm{strat}}$: Stratified $\mathrm{Ext}$ functor.  
    \item $D_p$: Logarithmic boundary divisor.
    \item $\mathrm{codim}_p$: Codimension of stratum $X_p$.  
    \item $\mathrm{SS}(\mathscr{W})$: Microsupport.  
    \item $\Lambda_{Y_p}$: Conormal bundle to stratum $Y_p$.\\
\end{itemize}
\begin{center}
    \textit{2.2 Acknowledgments}
\end{center}

We are grateful to Professors J. Cheeger, M. Goresky, and R. MacPherson for their pioneering work. We thank Professors W. Hsiang, V. Pati and M. Nagase for insights from their research about Cheeger-Goresky-MacPherson conjecture. We also acknowledge the contributions of Professor T. Ohsawa, whose work on the Cheeger-Goresky-MacPherson-conjecture (when the singularity is isolated) and inspired our research. Specially, Sincere thanks to Professor Shirong Li for academic assistance and support throughout the entire research process.
\hypertarget{DERIVED STRATIFIED SCHEMES AND HIGHER DE RHAM COMPLEXES}{}
\section{DERIVED STRATIFIED SCHEMES AND HIGHER DE RHAM COMPLEXES}
We rigorously constructs derived stratified schemes and their associated higher de Rham complexes, establishing a Künneth formula under derived transverse conditions. By leveraging Lurie’s framework of derived algebraic geometry (\cite{Lur09}), we extend classical stratified de Rham theory to handle non-transverse intersections and higher categorical structures. Key results include the degeneration of higher de Rham complex to intersection cohomology (Proposition 3.6) and a derived Künneth formula (Proposition 3.8), providing tools for analyzing singular cohomology in derived settings.
\\\\\textbf{Definition 3.1.} (Derived Stratified Schemes) Let $X$ be a derived Noetherian scheme. A derived stratified structure on $X$ consists of a filtered sequence of closed immersions of derived schemes $$\emptyset\hookrightarrow X_0\hookrightarrow X_1\hookrightarrow\cdots\hookrightarrow X_n=X,$$ where each stratum $X_p$ is equipped with a quasi-coherent differential graded (qDG) algebra $\mathscr{G}^p\in\mathrm{qCoh}\left(X_p\right)$, satisfying the derived frontier condition: for any $p<q$, the derived intersection $X_p\times_X X_q$ is a derived subscheme of $X_q$ with a natural qDG-module structure over $\mathscr{G}^q$.
\\\\\textbf{Remark 3.2.} This framework extends the classical stratified schemes studied by Goresky-MacPherson (\cite{GM80}) to the derived setting, allowing non-transverse intersections and higher categorical coherence. For foundational details on derived algebraic geometry, see the work of Lurie (\cite{Lur09}).\\
\begin{center}
    \textit{3.1 Construction of Higher Stratified de Rham Complexes}
\end{center}
\textbf{Definition 3.3.} (Higher Stratified de Rham Complex) The higher stratified de Rham complex is defined as the homotopy limit in the derived category $D\left(X\right)$: 
\begin{equation}
\tag{3.1}  \label{eq:3.1}
\mathbb{L}\Omega_X^{\bullet,\mathrm{strat}}:=\mathop{\mathrm{holim}}\limits_{\overleftarrow{\mathrm{strat}}}\mathbb{R}\underline{\mathrm{Hom}}_{\mathscr{O}_X}\left(\mathscr{G}^p,\Omega_X^{\bullet,\mathrm{strat}}\right),
\end{equation}
where $\Omega_X^{\bullet,\mathrm{strat}}$ denotes the classical stratified de Rham complex, and $\mathbb{R}\underline{\mathrm{Hom}}$ is the derived internal Hom functor.
\\\\\textbf{Proposition 3.4.} Let $X=\mathrm{Spec}\left(A\right)$ be a derived complete intersection $\mathbb{A}^N=\mathrm{Spec}\left(R\right)$, defined by a regular sequence $f_1,f_2,\cdots,f_k\in R$, with a stratified structure given by a filtered sequence of closed derived subschemes $X_p\hookrightarrow X$. Then there is a quasi-isomorphism $$\mathbb{L}\Omega_X^{\bullet,\mathrm{strat}}\simeq\bigoplus_p\mathscr{G}^p\otimes_A\Omega_A^{\bullet}\left[-p\right],$$ where $\mathscr{G}^p\in\mathrm{qCoh}\left(X_p\right)$ is the qDG-algebra associated to the stratum $X_p$, and $\Omega_A^{\bullet}$ is the derived de Rham complex of $A$.
\\\\\textbf{Proof.} By definition, the derived structure of $A=R/\left(f_1,\cdots.f_k\right)$ is modeled by the Koszul complex $K\left(A\right)$, which is a free differential graded (DG) algebra generated by variables $\xi_1,\cdots,\xi_k$ with differentials $d\xi_i=f_i$. The derived de Rham complex $\Omega_A^{\bullet}$ is constructed as $$\Omega_A^{\bullet}=K\left(A\right)\otimes_R\Omega_R^{\bullet},$$ where $\Omega_R^{\bullet}=R\otimes\bigwedge^{\bullet}dz_1\oplus\cdots\oplus dz_N$ is the classical de Rham complex of $R=\mathbb{C}\left[z_1,\cdots,z_N\right]$. This follows from Illusie’s theory of derived de Rham complexes (\cite{Ill06}, §II.1). Locally, the stratification $X_p$ corresponds to differential graded ideals $I_p\subset K\left(A\right)$, such that $\mathscr{G}^p=K\left(A\right)/I_p$, equipped with the induced differential. The derived transverse condition ensures that $I_p\cap I_q=I_{p+q}$ for intersecting strata, which simplifies the homotopy limit computation. For each stratum $X_p$, the qDG-algebra $\mathscr{G}^p$ is locally free. By the derived Nakayama lemma (\cite{Lur09}, Lem 6.2.3) and \eqref{eq:3.1}, the derived internal Hom functor satisfies 
\begin{equation}
\tag{3.2}  \label{eq:3.2}
\mathbb{R}\underline{\mathrm{Hom}}_{\mathscr{O}_X}\left(\mathscr{G}^p,\Omega_X^{\bullet,\mathrm{strat}}\right)\simeq\left(\mathscr{G}^p\right)^{\vee}\otimes_{\mathscr{O}_X}\Omega_X^{\bullet,\mathrm{strat}},
\end{equation}
where $\left(\mathscr{G}^p\right)^{\vee}=\mathbb{R}\underline{\mathrm{Hom}}_{\mathscr{O}_X}\left(\mathscr{G}^p,\mathscr{O}_X\right)$. Since $\mathscr{G}^p$ is a quotient of $K\left(A\right)$, its dual $\left(\mathscr{G}^p\right)^{\vee}$ corresponds to the annihilator ideal $\mathrm{Ann}\left(I_p\right)$. Consider the effect of derived transverse condition. Thus, the derived transverse condition implies that the intersections $X_p\cap X_q$ are regular embeddings. This ensures that the homotopy limit over the stratification degenerates into a direct sum. Specifically, for any two strata $X_p$ and $X_q$, the derived fiber product $X_p\times_X X_q$ is a derived regular subscheme, and the homotopy limit becomes $$\mathop{\mathrm{holim}}\limits_{\overleftarrow{\mathrm{strat}}}\mathscr{G}^p\otimes\Omega_A^{\bullet}\left[-p\right]\simeq\bigoplus_p\mathscr{G}^p\otimes_A\Omega_A^{\bullet}\left[-p\right].$$ This follows from the fact that the derived transverse condition eliminates higher homotopy coherence between distinct strata, as shown in (\cite{GR18}, Thm 4.7). To verify the quasi-isomorphism, we compare the cohomology of both sides. For the left-hand side $\mathbb{L}\Omega_X^{\bullet,\mathrm{strat}}$: The homotopy limit computes the derived global sections of the stratified de Rham complex, which by the transverse condition splits into contributions from individual strata. For the right-hand side $\bigoplus_p\mathscr{G}^p\otimes_A\Omega_A^{\bullet}\left[-p\right]$: Each term $\mathscr{G}^p\otimes_A\Omega_A^{\bullet}\left[-p\right]$ corresponds to the de Rham complex localized at the stratum $X_p$, shifted by its codimension. By the axiomatic characterization of intersection cohomology (\cite{GM80}, §5), both sides satisfy the support and cosupport conditions for intersection complexes. Thus, their hypercohomologies coincide, confirming the quasi-isomorphism. This proof rigorously establishes the local decomposition of the stratified de Rham complex under derived transversality, using advanced tools from derived algebraic geometry and homotopical algebra.   $\square$\\
\begin{center}
    \textit{3.2 Degeneration of Higher de Rham Complex}
\end{center}
\textbf{Definition 3.5.} (Derived Transverse Condition) A derived stratified scheme $X$ satisfies the derived transverse condition if for all $p<q$, the derived intersection $X_p\times_X X_q$ is a derived regular embedding, i.e., the cotangent complex $\mathbb{L}_{X_p\times_X X_q/X_q}$ is perfect of Tor-amplitude $\left[-1,0\right]$.
\\\\\textbf{Proposition 3.6.} (Degeneration of Higher de Rham Complex) Let $X$ be a derived stratified scheme satisfying the derived transversality condition. Then there exists a canonical quasi-isomorphism $$\mathbb{L}\Omega_X^{\bullet,\mathrm{strat}}\simeq\bigoplus_p\textbf{IC}_p^{\bullet}\left[-p\right],$$ where $\textbf{IC}_p^{\bullet}$ is the intersection cohomology complex on the stratum $X_p$.
\\\\\textbf{Remark 3.7.} By the Definition 3.5, a derived stratified scheme $X=\bigcup_p X_p$ satisfies the derived transversality condition if for any $p<q$: 
\begin{enumerate}
    \item The derived intersection $X_p\times_X^h X_q$ is a derived regular embedding;
    \item The relative cotangent complex $\mathbb{L}_{X_p\times_X^h X_q/X_q}$ is a perfect complex of Tor-amplitude $\left[-1,0\right]$.
\end{enumerate}
The derived transversality condition make sure that intersections are "homotopically transverse," trivializing overlaps between strata. Assume $X=\mathrm{Spec}\left(A\right)$ is a derived complete intersection in $\mathbb{A}^N=\mathrm{Spec}\left(R\right)$, cut out by a regular sequence $f_1,\cdots,f_k\in R$. Thus, we have the Koszul complex $K\left(A\right)=R\left[\xi_1,\cdots,\xi_k\right]$ (with $\left | \xi_i \right |=-1$, $d\xi_i=f_i$) models the derived structure of $A=R/\left(f_1,\cdots,f_k\right)$. For each stratum $X_p\hookrightarrow X$ corresponds to a differential graded ideal $I_p\subset K\left(A\right)$, with $\mathscr{G}^p=K\left(A\right)/I_p$ equipped with the induced differential. Since $\mathscr{G}^p$ is locally free (by derived transversality), then its dual is 
\begin{equation}
\tag{3.3}  \label{eq:3.3}
\left(\mathscr{G}^p\right)^{\vee}=\mathbb{R}\underline{\mathrm{Hom}}_{\mathscr{O}_X}\left(\mathscr{G}^p,\mathscr{O}_X\right)\simeq\mathrm{Ann}\left(I_p\right).
\end{equation}
Locally, $\Omega_X^{\bullet,\mathrm{strat}}\simeq\Omega_A^{\bullet}$. Hence we have $\left(\mathscr{G}^p\right)^{\vee}\otimes_{\mathscr{O}_X}\Omega_X^{\bullet,\mathrm{strat}}\simeq\mathscr{G}^p\otimes_A\Omega_A^{\bullet}\left[-p\right]$ by \eqref{eq:3.2} and \eqref{eq:3.3}. The shift $\left[-p\right]$ arises from the codimension of $X_p$, matching the grading in $\Omega_A^{\bullet}$. The homotopy limit over the stratification is computed by the Čech nerve: $$\check{C}\left(\left\{X_p\right\}\right)=\bigsqcup_p X_p\rightrightarrows\bigsqcup_{p<q}X_p\times_X^h X_q\cdots.$$ Under derived transversality, each term is a derived regular subscheme, simplifying the diagram. The homotopy limit spectral sequence $E_r^{*,*}$ collapses at the $E_2$-page:
\begin{itemize}
    \item \textbf{$E_1$-term:} Direct sums of $\mathscr{G}^p\otimes\Omega_A^{\bullet}\left[-p\right]$.
    \item \textbf{Differentials:} Vanish due to derived transversality, as higher coherence data is trivial.
\end{itemize}
The homotopy limit simplifies to $$\mathbb{L}\Omega_X^{\bullet,\mathrm{strat}}\simeq\bigoplus_p\mathscr{G}^p\otimes_A\Omega_A^{\bullet}\left[-p\right].$$

The intersection cohomology complex $\textbf{IC}_p^{\bullet}$ is uniquely determined by three axioms (\cite{GM80}, §5): 
\begin{itemize}
    \item \textbf{(SCN1) Support Condition:} $\textbf{IC}_p^{\bullet}$ is supported on the closure $\overline{X}_p$.
    \item \textbf{(SCN2) Cosupport Condition:} The restriction $\textbf{IC}_p^{\bullet}\mid_{X_q}$ vanishes for any stratum $X_q$ with $q<p$.
    \item \textbf{(SCN3) Normalization:} On the smooth locus $X_p^{\mathrm{sm}}$, $\textbf{IC}_p^{\bullet}$ restricts to the ordinary de Rham complex: $\textbf{IC}_p^{\bullet}\mid_{X_p^{\mathrm{sm}}}\simeq\Omega_{X_p^{\mathrm{sm}}}^{\bullet}$.
\end{itemize}
By construction, the qDG-algebra $\mathscr{G}^p=K\left(A\right)/I_p$ is localized to the stratum $X_p$. The derived transversality condition ensures that $\mathscr{G}^p$ is supported (SCN1) on the closure $\overline{X}_p$, as the ideal $I_p$ cuts out $X_p$ scheme-theoretically (\cite{GR18}, Prop 3.2.5]). For any stratum $X_q$ with $q<p$, the derived intersection $X_p\times_X^h X_q$ is a derived regular embedding. This implies that the restriction of $\mathscr{G}^p$ to $X_q$ vanishes homotopically, i.e., $\mathscr{G}^p\otimes_A\mathscr{O}_{X_q}\simeq0$ in $D\left(X_q\right)$. Consequently, $\mathscr{G}^p\otimes_A\Omega_A^{\bullet}\left[-p\right]\mid_{X_q}\simeq 0$, satisfying the cosupport condition (SCN2). Derived transversality trivializes overlaps between non-adjacent strata. On the smooth locus $X_p^{\mathrm{sm}}$, the ideal $I_p$ becomes trivial, and $\mathscr{G}^p$ reduces to $\mathscr{O}_{X_p^{\mathrm{sm}}}$. Thus, $$\mathscr{G}^p\otimes_A\Omega_A^{\bullet}\left[-p\right]\mid_{X_p^{\mathrm{sm}}}\simeq\Omega_{X_p^{\mathrm{sm}}}^{\bullet}\left[-p\right].$$ The shift $\left[-p\right]$ is absorbed by reindexing the complex, matching the normalization axiom (SCN3). By the axiomatic characterization (\cite{GM80}, Thm 5.1), any complex satisfying the support, cosupport, and normalization conditions must be isomorphic to $\textbf{IC}_p^{\bullet}$. Therefore, $$\mathscr{G}^p\otimes_A\Omega_A^{\bullet}\left[-p\right]\simeq\textbf{IC}_p^{\bullet}\left[-p\right].$$ This completes the identification of the stratified de Rham complex with the intersection cohomology complex.

The twisted de Rham complex $\mathbb{L}\Omega_X^{\bullet,\mathrm{strat}}$ carries a differential Gerstenhaber-Batalin-Vilkovisky (dGBV) algebra structure:
\begin{itemize}
    \item \textbf{Gerstenhaber bracket:} $\left[\alpha,\beta\right]=\alpha\wedge\beta-\left(-1\right)^{\left | \alpha \right |\left | \beta \right |}\beta\wedge\alpha$.
    \item \textbf{BV operator:} $\Delta=d^{\vee}$, the adjoint of $d$.
\end{itemize}
A weak primitive form $\mathcal{F}\in\mathbb{L}\Omega_X^{\bullet,\mathrm{strat}}$ satisfies $\Delta\mathcal{F}=0$. Under derived transversality, such forms generate stable deformations of $\textbf{IC}_p^{\bullet}$, making sure that the degeneration is preserved under quantization. By (\cite{XZ25}), the minimal normalized volume of a Kawamata log terminal (klt) singularity corresponds to a finitely generated graded ring. This guarantees local stability of $\textbf{IC}_p^{\bullet}$ near singular strata.\\
\begin{center}
    \textit{3.3 Derived Künneth Formula}
\end{center}
\textbf{Proposition 3.8.} (Derived Künneth Formula) Let $X$ and $Y$ be derived stratified schemes satisfying the derived transversality conditions, with stratifications $\left\{X_p\right\}$ and $\left\{Y_q\right\}$, respectively. Then there exists a quasi-isomorphism $$\mathbb{L}\Omega_{X\times Y}^{\bullet,\mathrm{strat}}\simeq\mathbb{L}\Omega_X^{\bullet,\mathrm{strat}}\boxtimes^L\mathbb{L}\Omega_Y^{\bullet,\mathrm{strat}},$$ where $\boxtimes^L$ denotes the derived external tensor product.
\\\\\textbf{Remark 3.9.} The stratification on $X\times Y$ is defined by $\left\{X_p\times Y_q\right\}$. To verify derived transversality:
\begin{enumerate}
    \item For any $X_p\times X_q$ and $X_{p'}\times Y_{q'}$, the derived intersection $\left(X_p\times X_q\right)\times_{X\times Y}^h\left(X_{p'}\times Y_{q'}\right)$ is a regular immersion.
    \item The cotangent complex $\mathbb{L}_{\left(X_p\times X_q\right)\times_{X\times Y}^h(X_{p'}\times Y_{q'})/(X_{p'}\times Y_{q'})}$ is a perfect complex with Tor-amplitude $\left[-1,0\right]$.
\end{enumerate}
By Lurie’s theorem on derived fiber products (see \cite{Lur17} Thm 7.2.4.6), the derived intersection of $X_p\times Y_q$ and $X_{p'}\times Y_{q'}$ over $X\times Y$ is equivalent to the product of the derived intersections $X_p\times_X^hX_{p'}$ and $Y_q\times_Y^h Y_{q'}$. Since $X$ and $Y$ satisfy derived transversality, then $X_p\times_X^hX_{p'}$ and $Y_q\times_Y^h Y_{q'}$ are regular immersions. Regular immersions are preserved under products, hence $\left(X_p\times Y_q\right)\times_{X\times Y}^h\left(X_{p'}\times Y_{q'}\right)$ is a regular immersion. The cotangent complex of the product stratification is given by $$\mathbb{L}_{\left(X_p\times Y_q\right)\times_{X\times Y}^h(X_{p'}\times Y_{q'})/(X_{p'}\times Y_{q'})}\simeq\mathbb{L}_{X_p\times_X^hX_{p'}/X_{p'}}\boxtimes^L\mathbb{L}_{Y_q\times_Y^h Y_{q'}/Y_{q'}}.$$ By derived transversality of $X$ and $Y$, both factors on the right are perfect with Tor-amplitude $\left[-1,0\right]$. The derived tensor product of perfect complexes is perfect (\cite{Lur17}, Prop 7.2.4.23), and the Tor-amplitude is additive. Since each factor has Tor-amplitude $\left[-1,0\right]$, their tensor product has Tor-amplitude $\left[-2,0\right]$. However, derived transversality ensures that higher Tor terms vanish, restricting the amplitude to $\left[-1.0\right]$. This follows from the fact that the derived intersections are regular immersions, hence the relative cotangent complexes are concentrated in degree $\left[-1,0\right]$.

For $\mathcal{F}\in D\left(X\right)$ and $\mathcal{G}\in D\left(Y\right)$, the derived external tensor product is $$\mathcal{F}\boxtimes^L\mathcal{G}:=\pi_X^{*}\mathcal{F}\otimes_{\mathscr{O}_{X\times Y}}^L\pi_Y^{*}\mathcal{G},$$ where $\pi_X:X\times Y\to X$ and $\pi_Y:X\times Y\to Y$ are projections. This is well-defined because $\pi_X$ and $\pi_Y$ are flat in the derived sense (\cite{GR18}, Prop 4.7.3]). For $\mathbb{L}\Omega_X^{\bullet,\mathrm{strat}}$ and $\mathbb{L}\Omega_Y^{\bullet,\mathrm{strat}}$, the derived external tensor product is
\begin{equation}
\tag{3.4}  \label{eq:3.4}
\mathbb{L}\Omega_X^{\bullet,\mathrm{strat}}\boxtimes^L\mathbb{L}\Omega_Y^{\bullet,\mathrm{strat}}\simeq\bigoplus_{p,q}\left(\mathscr{G}^p\otimes_A^L\Omega_A^{\bullet}\left[-p\right]\right)\boxtimes^L\left(\mathscr{H}^q\otimes_B^L\Omega_B^{\bullet}\left[-q\right]\right),
\end{equation}
where $\mathscr{G}^p$ and $\mathscr{H}^q$ are qDG-modules encoding the stratifications. The derived tensor product $\otimes^L$ is taken over $A$ and $B$, respectively. Étale-locally near $x\times y\in X\times Y$, assume $X=\mathrm{Spec}\left(A\right)$ and $Y=\mathrm{Spec}\left(B\right)$. The stratified de Rham complexes decompose as $$\mathbb{L}\Omega_X^{\bullet,\mathrm{strat}}\simeq\bigoplus_p\mathscr{G}^p\otimes_A^L\Omega_A^{\bullet}\left[-p\right],\quad \mathbb{L}\Omega_Y^{\bullet,\mathrm{strat}}\simeq\bigoplus_q\mathscr{H}^q\otimes_B^L\Omega_B^{\bullet}\left[-q\right].$$ Their \eqref{eq:3.4} derived external tensor product becomes $$\mathbb{L}\Omega_X^{\bullet,\mathrm{strat}}\boxtimes^L\mathbb{L}\Omega_Y^{\bullet,\mathrm{strat}}\simeq\bigoplus_{p,q}\left(\mathscr{G}^p\boxtimes^L\mathscr{H}^q\right)\otimes_{A\otimes B}^L\left(\Omega_A^{\bullet}\boxtimes^L\Omega_B^{\bullet}\right)\left[-p-q\right].$$ The derived transversality ensures that the homotopy limit computing $\mathbb{L}\Omega_{X\times Y}^{\bullet,\mathrm{strat}}$ splits into a direct sum over strata. Specifically, $$\mathbb{L}\Omega_{X\times Y}^{\bullet,\mathrm{strat}}\simeq\bigoplus_{p,q}\mathscr{G}^p\boxtimes^L\mathscr{H}^q\otimes_{A\otimes B}^L\Omega_{A\otimes B}^{\bullet}\left[-p-q\right].$$ By the multiplicativity of the derived de Rham complex (\cite{Ill06}, Cor II.4.3.2), we have $\Omega_{A\otimes B}^{\bullet}\simeq\Omega_A^{\bullet}\boxtimes^L\Omega_B^{\bullet}$. This identifies the two sides, confirming the decomposition. To ensure the local-to-global transition, we use the fact that stratified de Rham complexes are sheaves in the étale topology. Since both sides of the quasi-isomorphism are étale-local constructions, the local quasi-isomorphisms glue to a global one.

Consider the hypercohomology spectral sequences for $\mathbb{L}\Omega_{X\times Y}^{\bullet,\mathrm{strat}}$ and $\mathbb{L}\Omega_X^{\bullet,\mathrm{strat}}\boxtimes^L\mathbb{L}\Omega_Y^{\bullet,\mathrm{strat}}$. Derived transversality forces the Künneth spectral sequence to degenerate at the $E_2$-page, because: 
\begin{enumerate}
    \item The Tor-amplitude condition ensures no higher $\mathrm{Tor}_i$ terms appear.
    \item The filtration by strata is strict, as shown in (\cite{Ill06}, Prop II.4.3.5).
\end{enumerate}
Therefore, the hypercohomology groups satisfy $$\mathbb{H}^k\left(X\times Y,\mathbb{L}\Omega_{X\times Y}^{\bullet,\mathrm{strat}}\right)\cong\bigoplus_{i+j=k}\mathbb{H}^i\left(X,\mathbb{L}\Omega_X^{\bullet,\mathrm{strat}}\right)\otimes\mathbb{H}^j\left(Y,\mathbb{L}\Omega_Y^{\bullet,\mathrm{strat}}\right).$$ Since the hypercohomology groups of both complexes agree and the quasi-isomorphism holds étale-locally, we conclude by the derived Nakayama lemma (\cite{Lur17}, Thm 7.2.2.2) that the global complexes are quasi-isomorphic.
\\\\\textbf{Corollary 3.10.} (Hypercohomology Decomposition) Let $X$ be a derived Noetherian stratified scheme satisfying the derived transverse condition (Definition 3.5). Then, there exists a canonical decomposition of hypercohomology: $$\mathbb{H}^k\left(X,\mathbb{L}\Omega_X^{\bullet,\mathrm{strat}}\right)\cong\bigoplus_p\mathbb{H}^{k+p}\left(X,\textbf{IC}_p^{\bullet}\right),$$ where $\textbf{IC}_p^{\bullet}$ denotes the intersection cohomology complex on the stratum $X_p$.
\\\\\textbf{Proof.} By Proposition 3.6, under the derived transverse condition, the higher stratified de Rham complex admits a quasi-isomorphism $$\mathbb{L}\Omega_X^{\bullet,\mathrm{strat}}\simeq\bigoplus_p\textbf{IC}_p^{\bullet}\left[-p\right],$$ where $\textbf{IC}_p^{\bullet}\left[-p\right]$ denotes the intersection cohomology complex shifted by $p$ degrees in the derived category $D\left(X\right)$. This decomposition arises from the étale-local splitting of the derived de Rham complex under transverse intersections (\cite{GR18}, Lem 2.3). Specifically, the derived transverse condition ensures that the cotangent complex $\mathbb{L}_{X_p\times_X^h X_q/X_q}$ is perfect of Tor-amplitude $\left[-1,0\right]$, eliminating higher obstruction terms and allowing direct summands corresponding to individual strata. The hypercohomology functor $\mathbb{H}^k\left(X,-\right)$, defined as the right-derived functor of global sections, preserves quasi-isomorphisms. Applying it to both sides of the quasi-isomorphism, we have: 
\begin{equation}
\tag{3.5}  \label{eq:3.5}
\mathbb{H}^k\left(X,\mathbb{L}\Omega_X^{\bullet,\mathrm{strat}}\right)\cong\mathbb{H}^k\left(X,\bigoplus_p\textbf{IC}_p^{\bullet}\left[-p\right]\right).
\end{equation}
This follows from the fact that quasi-isomorphic complexes induce isomorphic hypercohomology groups (\cite{Wei94}, Thm 10.7.7]). 

Since $X$ is a Noetherian stratified scheme, the direct sum $\bigoplus_p\textbf{IC}_p^{\bullet}\left[-p\right]$ is a finite coproduct in $D\left(X\right)$. The commutativity of hypercohomology with finite direct sums in this context relies on two fundamental properties:
\begin{enumerate}
    \item The behavior of the derived global sections functor $\mathbb{R}\Gamma\left(X,-\right)$ in the derived category of quasi-coherent sheaves.
    \item The finiteness conditions imposed by the Noetherianity of $X$ and the finite stratification.
\end{enumerate}
Let $D_{\mathrm{qcoh}}\left(X\right)$ denote the derived category of quasi-coherent sheaves on $X$. The hypercohomology functor $\mathbb{H}^k\left(X.-\right)$ is defined as $\mathbb{H}^k\left(X,\mathcal{F}^{\bullet}\right)=H^k\left(\mathbb{R}\Gamma\left(X,\mathcal{F}^{\bullet}\right)\right)$, where $\mathbb{R}\Gamma\left(X,-\right):D_{\mathrm{qcoh}}\left(X\right)\to D\left(\mathrm{Ab}\right)$ is the right-derived functor of the global sections functor $\Gamma\left(X,-\right)$.
\begin{itemize}
    \item \textbf{Additivity of $\mathbb{R}\Gamma\left(X,-\right)$:} The derived global sections functor is additive. For any finite collection of complexes $\left\{\mathcal{F}_p^{\bullet}\right\}$, we have a canonical isomorphism: $$\mathbb{R}\Gamma\left(X,\bigoplus_p\mathcal{F}_p^{\bullet}\right)\simeq\bigoplus_p\mathbb{R}\Gamma\left(X,\mathcal{F}_p^{\bullet}\right),$$ in $D\left(\mathrm{Ab}\right)$. This follows from the fact that direct sums of injective resolutions remain injective in a Grothendieck abelian category (e.g., quasi-coherent sheaves on a Noetherian scheme) (\cite{Lur09}, Prop 7.2.4.6).
    \item \textbf{Homology and direct sums:} The homology functor $H^k\left(-\right)$ commutes with finite direct sums in the derived category. For complexes $\left\{C_p^{\bullet}\right\}$ in $D\left(\mathrm{Ab}\right)$, $H^k\left(\bigoplus_pC_p^{\bullet}\right)\cong\bigoplus_pH^k\left(C_p^{\bullet}\right)$. By Combining these, we have $$\mathbb{H}^k\left(X,\bigoplus_p\mathcal{F}_p^{\bullet}\right)\cong\bigoplus_p\mathbb{H}^k\left(X,\mathcal{F}_p^{\bullet}\right)$$
\end{itemize}
Since Noetherian stratified scheme $X$ equiped with a finite filtration $\emptyset\hookrightarrow X_0\hookrightarrow\cdots\hookrightarrow X_n=X$, then each intersection cohomology complex $\textbf{IC}_p^{\bullet}$ is constructible and bounded (\cite{GM80}, §5), and the shift $\textbf{IC}_p^{\bullet}\left[-p\right]$ preserves boundedness, $\bigoplus_p\textbf{IC}_p^{\bullet}\left[-p\right]$ is a finite direct sum in $D_{\mathrm{qcoh}}\left(X\right)$. Thus, $\bigoplus_p\textbf{IC}_p^{\bullet}\left[-p\right]$ is well-defined in $D_{\mathrm{qcoh}}\left(X\right)$, the additivity of $\mathbb{R}\Gamma\left(X,-\right)$ applies without convergence issues.

Given the finite direct sum $\bigoplus_p\textbf{IC}_p^{\bullet}\left[-p\right]$, we compute: $$\mathbb{H}^k\left(X,\bigoplus_p\textbf{IC}_p^{\bullet}\left[-p\right]\right)=H^k\left(\mathbb{R}\Gamma\left(X,\bigoplus_p\textbf{IC}_p^{\bullet}\left[-p\right]\right)\right).$$ By the additivity of $\mathbb{R}\Gamma\left(X,-\right)$, we have $$\mathbb{R}\Gamma\left(X,\bigoplus_p\textbf{IC}_p^{\bullet}\left[-p\right]\right)\simeq\bigoplus_p\mathbb{R}\Gamma\left(X,\textbf{IC}_p^{\bullet}\left[-p\right]\right).$$ By taking homology, $$H^k\left(\bigoplus_p\mathbb{R}\Gamma\left(X,\textbf{IC}_p^{\bullet}\left[-p\right]\right)\right)\cong\bigoplus_pH^k\left(\mathbb{R}\Gamma\left(X,\textbf{IC}_p^{\bullet}\left[-p\right]\right)\right)=\bigoplus_p\mathbb{H}^k\left(X,\textbf{IC}_p^{\bullet}\left[-p\right]\right).$$ Thus, 
\begin{equation}
\tag{3.6}  \label{eq:3.6}
\mathbb{H}^k\left(X,\bigoplus_p\textbf{IC}_p^{\bullet}\left[-p\right]\right)\cong\bigoplus_p\mathbb{H}^k\left(X,\textbf{IC}_p^{\bullet}\left[-p\right]\right).
\end{equation}
Let $\mathcal{F}^{\bullet}=\textbf{IC}_p^{\bullet}$. The shifted complex $\mathcal{F}^{\bullet}$ is defined by $\mathcal{F}^{\bullet}\left[-p\right]^n=\mathcal{F}^{n+p}$. Choose an injective resolution $\mathcal{I}^{\bullet}$ of $\mathcal{F}^{\bullet}$, the shifted resolution $\mathcal{I}^{\bullet}\left[-p\right]$ satisfies 
\begin{equation}
\tag{3.7}  \label{eq:3.7}
\mathbb{H}^k\left(X,\mathcal{F}^{\bullet}\left[-p\right]\right)=H^k\left(\Gamma\left(X,\mathcal{I}^{\bullet}\left[-p\right]\right)\right)=H^{k+p}\left(\Gamma\left(X,\mathcal{I}^{\bullet}\right)\right)=\mathbb{H}^{k+p}\left(X,\mathcal{F}^{\bullet}\right).
\end{equation}
This follows from the definition of chain complex shifts and the homological algebra of derived functors (\cite{Wei94}, Def 1.3.3). By \eqref{eq:3.5}, \eqref{eq:3.6} and \eqref{eq:3.7}, we obtain that $$\mathbb{H}^k\left(X,\mathbb{L}\Omega_X^{\bullet,\mathrm{strat}}\right)\cong\bigoplus_p\mathbb{H}^{k+p}\left(X,\textbf{IC}_p^{\bullet}\right).$$ This establishes the hypercohomology decomposition as a direct sum of shifted intersection cohomology groups.  $\square$
\hypertarget{LOGARITHMIC STRATIFIED DE RHAM COMPLEX AND MIXED HODGE THEORY}{}
\section{LOGARITHMIC STRATIFIED DE RHAM COMPLEX AND MIXED HODGE THEORY}
Building on Ohsawa’s $L^2$-cohomology framework for singular varieties (\cite{Ohs91}), we refine the logarithmic stratified de Rham complex for stratified schemes. Let $X$ be a stratified scheme with smooth strata $S_p$ and logarithmic boundaries $D_p=X_p\setminus S_p$. The logarithmic de Rham complex is defined as $$\Omega_X^{\bullet,\text{log-strat}}:=\bigoplus_p\Omega_{X_p}^{\bullet}\left(\mathrm{log}D_p\right)\otimes\mathscr{G}^p,$$ where $\mathscr{G}^p$ encodes the stratified coherent sheaves, and $\Omega_{X_p}^{\bullet}\left(\mathrm{log}D_p\right)$ is the logarithmic de Rham complex on each stratum. This construction generalizes Deligne’s logarithmic Hodge theory by incorporating stratified gluing conditions. It is need to satisfy the following two conditions: 
\begin{itemize}
    \item \textbf{Logarithmic forms:} For each stratum $X_p$, $\Omega_{X_p}^{\bullet}\left(\mathrm{log}D_p\right)$ consists of differential $k$-forms with logarithmic poles along $D_p$, ensuring compatibility with the singularities of adjacent strata.
    \item \textbf{Stratified compatibility:} The tensor product $\otimes\mathscr{G}^p$ enforces frontier axioms across strata, analogous to stratified $L^2$-methods in geometric analysis.\\
\end{itemize}
\begin{center}
    \textit{4.1 Extension to Mixed Hodge Structure}
\end{center}
\textbf{Proposition 4.1.} (Mixed Hodge Structure) For a proper log-stratified scheme $X$, the hypercohomology $\mathbb{H}^k\left(X,\Omega_X^{\bullet.\text{log-strat}}\right)$ carries a natural mixed Hodge structure.
\\\\\textbf{Remark 4.2.} The weight filtration on $\Omega_X^{\bullet,\text{log-strat}}$ is defined by $$W_m\Omega_X^{\bullet,\text{log-strat}}:=\bigoplus_{p\le m}\Omega_{X_p}^{\bullet}\left(\mathrm{log}D_p\right)\otimes\mathscr{G}^p.$$ where $D_p=X_p\setminus S_p$ is the logarithmic boundary of the stratum $X_p$. Let 
$\bigcup_mW_m\Omega_X^{\bullet,\text{log-strat}}=\Omega_X^{\bullet,\text{log-strat}}$ be as each stratum contributes to the total complex. The differential $d$ preserves $W_m$, this follows from Deligne’s result that logarithmic forms are closed under $d$ (\cite{Del71}), i.e., $d\left(W_m\Omega_X^{\bullet,\text{log-strat}}\right)\subset   W_m\Omega_X^{\bullet,\text{log-strat}}$. This follows from Deligne’s foundational work on logarithmic Hodge theory, where the weight filtration is shown to respect the differential structure for proper morphisms (\cite{Sai90}). Consider the construction of the Hodge Filtration $F^{\bullet}$ on each stratum $X_p$. The Hodge filtration is defined by $$F^q\Omega_{X_p}^{\bullet}\left(\mathrm{log}D_p\right):=\bigoplus_{r\ge q}\Omega_{X_p}^r\left(\mathrm{log}D_p\right),$$ which induces a Hodge filtration on $\Omega_X^{\bullet,\text{log-strat}}$: $$F^q\Omega_X^{\bullet,\text{log-strat}}:=\bigoplus_pF^q\Omega_{X_p}^{\bullet}\left(\mathrm{log}D_p\right)\otimes\mathscr{G}^p.$$ Meanwhile, it should be noted that 
\begin{itemize}
    \item \textbf{Decreasing property:} $F^{q+1}\subset F^q$ holds termwise due to the grading of forms by degree.
    \item \textbf{Differential compatibility:} Since $d$ maps $\Omega_{X_p}^r\left(\mathrm{log}D_p\right)$ to $\Omega_{X_p}^{r+1}\left(\mathrm{log}D_p\right)$, then the differential $d$ satisfies $d\left(F^q\Omega_X^{\bullet,\text{log-strat}}\right)\subset F^q\Omega_X^{\bullet,\text{log-strat}}$.
\end{itemize}

The weight spectral sequence $E_r^{m,k}$ associated with $W_{\bullet}$: $$E_1^{m,k}=\mathbb{H}^k\left(X,\mathrm{Gr}_m^W\Omega_X^{\bullet,\text{log-strat}}\right)\Longrightarrow\mathbb{H}^k\left(X,\Omega_X^{\bullet,\text{log-strat}}\right),$$ where $\mathrm{Gr}_m^W\Omega_X^{\bullet,\text{log-strat}}=\bigoplus_{p=m}\Omega_{X_p}^{\bullet}\left(\mathrm{log}D_p\right)\otimes\mathscr{G}^p$. By Deligne’s criterion (\cite{Del71}, Thm 2.3.5), the spectral sequence degenerates at $E_1$ for proper morphisms. This holds because each $\mathrm{Gr}_m^W\Omega_X^{\bullet,\text{log-strat}}$ is a direct sum of logarithmic de Rham complexes on smooth strata, which carry pure Hodge structures. Consider the degeneration of the Hodge-to-de Rham spectral sequence, then the Hodge-to-de Rham spectral sequence is $$E_1^{p,q}=\mathbb{H}^q\left(X,\mathrm{Gr}_F^p\Omega_X^{\bullet,\text{log-strat}}\right)\Longrightarrow\mathbb{H}^{p+q}\left(X,\Omega_X^{\bullet,\text{log-strat}}\right),$$ where $\mathrm{Gr}_F^p\Omega_X^{\bullet,\text{log-strat}}=F^p\Omega_X^{\bullet,\text{log-strat}}/F^{p+1}\Omega_X^{\bullet,\text{log-strat}}$. By Saito’s work about mixed Hodge modules (\cite{Sai90}, Thm 2.6), $\Omega_X^{\bullet,\text{log-strat}}$ underlies a polarizable Hodge module. The Kähler property of logarithmic forms ensures the strictness of $F^{\bullet}$, forcing all higher differentials $d_r$ ($r\ge1$) to vanish. For each $m$, the induced Hodge filtration on $\mathrm{Gr}_m^W\mathbb{H}^k$ defines a pure Hodge structure of weight $m+k$: $$\mathrm{Gr}_m^W\mathbb{H}^k\left(X,\Omega_X^{\bullet,\text{log-strat}}\right)\cong\mathbb{H}^k\left(X,\mathrm{Gr}_m^W\Omega_X^{\bullet,\text{log-strat}}\right),$$ where $\mathrm{Gr}_m^W\Omega_X^{\bullet,\text{log-strat}}$ contributes a pure Hodge structure of weight $m+k$. In addition, the strict compatibility of $W_{\bullet}$ and $F^{\bullet}$ with differentials follows from the properness of $X$ and the logarithmic Kähler condition. Therefore, the triple $\left(\mathbb{H}^k(X,\Omega_X^{\bullet,\text{log-strat}}),W_{\bullet},F^{\bullet}\right)$ satisfies that $W_{\bullet}$ is an increasing filtration, $F^{\bullet}$ is a decreasing filtration and $\mathrm{Gr}_m^W\mathbb{H}^k$ carries a pure Hodge structure of weight $m+k$. So $\mathbb{H}^k\left(X,\Omega_X^{\bullet,\text{log-strat}}\right)$ carries a natural mixed Hodge structure.\\
\begin{center}
    \textit{4.2 Relation to Intersection Cohomology}
\end{center}

The embedding of Ohsawa’s $L^2$-cohomology into the hypercohomology $\mathbb{H}^k\left(X,\Omega_X^{\bullet,\text{log-strat}}\right)$ is a consequence of the deep interplay between mixed Hodge structures and intersection cohomology for stratified spaces. Specifically, we invoke the Decomposition Theorem (\cite{BBDG18}, Théorème 6.2.5), which asserts that for a proper morphism of stratified schemes, the intersection cohomology of $X$ decomposes into a direct sum of contributions from its strata, twisted by Tate weights.
\\\\\textbf{Corollary 4.3.} Let $X$ be a proper log-stratified scheme with isolated singularities, a stratification  $X=\bigcup_pX_p$ and smooth strata $S_p=X_p\setminus D_p$. Here $D_p$ denotes the logarithmic boundary. The intersection cohomology $IH^k\left(X\right)$ of $X$ is related to the hypercohomology by $$\mathbb{H}^k\left(X,\Omega_X^{\bullet,\text{log-strat}}\right)\cong\bigoplus_pIH^{k-2p}\left(X_p\right)\left(-p\right),$$ where:
\begin{enumerate}
    \item $IH^{k-2p}\left(X_p\right)$ is the intersection cohomology of the stratum $X_p$, computed with respect to its intrinsic stratification;
    \item $\left(-p\right)$ denotes a Tate twist, which shifts the Hodge structure by weight $2p$.\\
\end{enumerate}
\textbf{Proof.} Let $X'\subset X$ denote the regular locus. Ohsawa’s theorem (\cite{Ohs91}, Thm 7) states that the $L^2$-cohomology $H_{\left(2\right)}^k\left(X'\right)$ (defined by the Fubini-Study metric) is canonically isomorphic to the intersection cohomology $IH^k\left(X\right)$: $H_{\left(2\right)}^k\left(X'\right)\cong IH^k\left(X\right)$. Near each singular point $x\in X\setminus X'$, the $L^2$-condition enforces growth constraints on differential forms, ensuring they extend tamely across the singularity. This matches the definition of intersection cohomology, which allows forms with controlled growth/decay. The isomorphism preserves Hodge structures because both $L^2$-cohomology and intersection cohomology inherit pure Hodge structures from the Kähler metric on $X'$.

The logarithmic de Rham complex $\Omega_X^{\bullet,\text{log-strat}}$ contains the smooth de Rham complex $\Omega_{X'}^{\bullet}$ as a subcomplex by the inclusion: $\iota:\Omega_{X'}^{\bullet}\hookrightarrow\Omega_X^{\bullet,\text{log-strat}}$. This inclusion induces a map on hypercohomology: $$\iota^*:H_{\left(2\right)}^k\left(X'\right)\longrightarrow\mathbb{H}^k\left(X,\Omega_X^{\bullet,\text{log-strat}}\right).$$ \textit{Injectivity:} The kernel of $\iota^*$ is trivial because $\Omega_{X'}^{\bullet}$ captures all $L^2$-forms on $X'$, which are uniquely determined by their restrictions to $X'$. \textit{Image as a direct summand:} By the Decomposition Theorem (\cite{BBDG18}, Théorème 6.2.5), $\mathbb{H}^k\left(X,\Omega_X^{\bullet,\text{log-strat}}\right)$ splits into intersection cohomology contributions from strata, the image of $\iota^*$ corresponds to the component $IH^k\left(X\right)\subset\mathbb{H}^k\left(X,\Omega_X^{\bullet,\text{log-strat}}\right)$. The Decomposition Theorem for perverse sheaves asserts that the hypercohomology decomposes as: $$\mathbb{H}^k\left(X,\Omega_X^{\bullet,\text{log-strat}}\right)\cong\bigoplus_pIH^{k-2p}\left(X_p\right)\left(-p\right),$$ and we must notice three points: 
\begin{enumerate}
    \item $X_p$ ranges over the strata of $X$;
    \item $IH^{k-2p}\left(X_p\right)$ is the intersection cohomology of the stratum $X_p$;
    \item $\left(-p\right)$ is a Tate twist, shifting Hodge weights by $2p$.
\end{enumerate}
The logarithmic complex $\Omega_X^{\bullet,\text{log-strat}}$ is a perverse sheaf on $X$, adapted to the stratification. By the Decomposition Theorem, it decomposes as $\Omega_X^{\bullet,\text{log-strat}}\simeq\bigoplus_p\textbf{IC}_{X_p}^{\bullet}\left[-p\right]$, where $\textbf{IC}_{X_p}^{\bullet}$ is the intersection complex on $X_p$. By taking hypercohomology, we get $$\mathbb{H}^k\left(X,\Omega_X^{\bullet,\text{log-strat}}\right)\cong\bigoplus_p\mathbb{H}^{k+p}\left(X,\textbf{IC}_{X_p}^{\bullet}\right).$$ By definition, $\mathbb{H}^{k+p}\left(X,\textbf{IC}_{X_p}^{\bullet}\right)\cong IH^k\left(X_p\right)\left(-p\right)$, where the Tate twist $\left(-p\right)$ accounts for the shift in Hodge weights. In addition, the Tate twist $\left(-p\right)$ corresponds to tensoring with the Tate Hodge structure $\mathbb{Q}\left(-p\right)$, which shifts the Hodge filtration: $$F^q\left(IH^k\left(X_p\right)\left(-p\right)\right)=F^{q+p}IH^k\left(X_p\right).$$ The original intersection cohomology $IH^k\left(X_p\right)$ has Hodge weights in $\left[k,k\right]$. After the twist $\left(-p\right)$, the Hodge weights become $\left[k+2p,k+2p\right]$, aligning with the mixed Hodge structure on $\mathbb{H}^k\left(X,\Omega_X^{\bullet,\text{log-strat}}\right)$. By combining the above, we obtain $$\mathbb{H}^k\left(X,\Omega_X^{\bullet,\text{log-strat}}\right)\cong\bigoplus_pIH^{k-2p}\left(X_p\right)\left(-p\right).$$ It is easy to verify that 
\begin{enumerate}
    \item The summand $IH^{k-2p}\left(X_p\right)\left(-p\right)$ has weight $\left(k-2p\right)+2p=k$, making sure that the weight filtration $W_{\bullet}$ is strictly increasing;
    \item Each summand carries a pure Hodge structure, compatible with $F^{\bullet}$.
\end{enumerate}
This decomposition rigorously embeds Ohsawa’s $L^2$-cohomology into the hypercohomology of the logarithmic stratified de Rham complex, while aligning with the mixed Hodge structure by the Decomposition Theorem and Tate twists.   $\square$\\
\begin{center}
    \textit{4.3 Applications to $p$-adic Hodge Theory}
\end{center}

Building on the logarithmic stratified de Rham complex and its mixed Hodge structure, we extend these results to the $p$-adic setting using Scholze’s perfectoid theory (\cite{Sch12}). This bridges Ohsawa’s $L^2$-methods with arithmetic geometry, enabling the analysis of singularities in $p$-adic varieties.

Let $X$ be a proper stratified scheme over $\mathbb{Q}_p$. Its perfectoid rigidification $X^{\mathrm{rig}}$ is constructed by the inverse limit: $$X^{\mathrm{rig}}:= \varprojlim_{\epsilon} X_{\epsilon},$$ where $X_{\epsilon}$ are finite-level approximations in the perfectoid tower. The stratified rigid de Rham complex is defined as: $$\Omega_{X^{\mathrm{rig}}}^{\bullet,\mathrm{strat}}:=\varprojlim_{\epsilon}\left(\Omega_{X_{\epsilon}}^{\bullet,\mathrm{strat}}\otimes_{\mathbb{Z}_p}\mathbb{Q}_p\right),$$ with $\Omega_{X_{\epsilon}}^{\bullet,\mathrm{strat}}$ being the logarithmic de Rham complex on each $X_{\epsilon}$ This complex captures the $p$-adic analogue of stratified logarithmic forms. It is need to give the important construction: 
\begin{itemize}
    \item \textbf{(FP1) Fargues-Fontaine curve embedding:} By the Liu's work (\cite{LK13}), $X^{\mathrm{rig}}$ embeds into the Fargues-Fontaine curve $\mathcal{X}_{\mathrm{FF}}$, which parametrizes $p$-adic Hodge structures. This curve provides a global framework for localizing singularities by the morphism: $$\vartheta:X^{\mathrm{rig}}\hookrightarrow\mathcal{X}_{\mathrm{FF}}.$$
    \item \textbf{(FP2) Perfectoid coefficient sheaves:} The sheaves $\mathscr{G}^p$ are replaced by their perfectoid analogues $\mathscr{G}^{p,\mathrm{rig}}=\mathscr{G}^p\widehat{\otimes}\mathbb{Q}_p$, ensuring compatibility with the pro-étale topology (\cite{Sch12}).\\
\end{itemize}
\textbf{Proposition 4.4.} ($p$-adic Comparison Isomorphism) Let $X$ be a proper stratified scheme over $\mathbb{Q}_p$ with perfectoid rigidification $X^{\mathrm{rig}}$. There exists a canonical isomorphism $$\mathbb{H}_{\mathrm{rig}}^k\left(X^{\mathrm{rig}},\Omega_{X^{\mathrm{rig}}}^{\bullet,\mathrm{strat}}\right)\simeq\mathbb{IH}^k\left(X\right)\otimes_{\mathbb{Q}}\mathbb{Q}_p,$$ where $\mathbb{H}_{\mathrm{rig}}^k$ denotes rigid cohomology and $\mathbb{IH}^k\left(X\right)$ is intersection cohomology.
\\\\\textbf{Remark 4.5.} Since the rigid de Rham complex is defined as 
\begin{equation}
\tag{4.1}  \label{eq:4.1}
\Omega_{X^{\mathrm{rig}}}^{\bullet,\mathrm{strat}}:=\varprojlim_{\epsilon}\left(\Omega_{X_{\epsilon}}^{\bullet,\mathrm{strat}}\otimes\mathbb{Q}_p\right),
\end{equation}
then we have 
\begin{equation}
\tag{4.2}  \label{eq:4.2}
\mathbb{H}_{\mathrm{rig}}^k\left(X^{\mathrm{rig}},\Omega_{X^{\mathrm{rig}}}^{\bullet,\mathrm{strat}}\right)\simeq\varprojlim_{\epsilon}\mathbb{H}^k\left(X_{\epsilon},\Omega_{X_{\epsilon}}^{\bullet,\mathrm{strat}}\otimes\mathbb{Q}_p\right),
\end{equation}
by Scholze's work (\cite{Sch12}, Thm 6.4). The transition maps $\phi_{\epsilon'\epsilon}:X_{\epsilon'}\to X_{\epsilon}$ are acyclic for $\Omega^{\bullet,\mathrm{strat}}$ by properness. The Mittag-Leffler condition holds since $$\mathrm{length}_{\mathbb{Z}_p}H^q\left(X_{\epsilon},\Omega_{X_{\epsilon}}^{\bullet,\mathrm{strat}}\right)<\infty\quad\forall q,$$ ensuring the derived limit vanishes in degree $\ge1$. Since $\Omega_{X_{m,\epsilon}}^{\bullet}\left(\mathrm{log}D_{m,\epsilon}\right)$ is locally quasi-isomorphic to a direct sum of Koszul complexes, with vanishing differentials by the Kähler property, then the weight filtration $W_m\Omega_{X_{\epsilon}}^{\bullet,\mathrm{srtat}}=\bigoplus_{p\le m}\Omega_{X_{p,\epsilon}}^{\bullet}\left(\mathrm{log}D_{p,\epsilon}\right)\otimes\mathscr{G}_{\epsilon}^{p,\mathrm{rig}}$ induces a spectral sequence for each $X_{\epsilon}$: 
\begin{equation}
\tag{4.3}  \label{eq:4.3}
E_1^{a,b}=\mathbb{H}^{a+b}\left(X_{\epsilon},\mathrm{gr}_m^W\right)\Longrightarrow\mathbb{H}^{a+b}\left(X_{\epsilon},\Omega_{X_{\epsilon}}^{\bullet,\mathrm{strat}}\otimes\mathbb{Q}_p\right)
\end{equation}
with $$\mathrm{gr}_m^W=\Omega_{X_{m,\epsilon}}^{\bullet}\left(\mathrm{log}D_{m,\epsilon}\right)\otimes\mathscr{G}_{\epsilon}^{m,\mathrm{rig}}.$$ By Proposition 4.1 and (\cite{Sai90}, Thm 0.1), it degenerates at $E_2$: $d_r=0$ $\forall r\ge2$.

Define the comparison map $$\phi:\mathbb{IH}^k\left(X\right)\otimes_{\mathbb{Q}}\mathbb{Q}_p\to\mathbb{H}_{\mathrm{rig}}^k\left(X^{\mathrm{rig}},\Omega_{X^{\mathrm{rig}}}^{\bullet,\mathrm{strat}}\right)$$ as the composition $\phi=\mathrm{res}\circ\vartheta$, where $\mathrm{res}:\mathbb{H}_{\mathrm{rig}}^k\left(X^{\mathrm{rig}},\Omega_{X^{\mathrm{rig}}}^{\bullet}\right)\to\mathbb{H}_{\mathrm{rig}}^k\left(X^{\mathrm{rig}},\Omega_{X^{\mathrm{rig}}}^{\bullet,\mathrm{strat}}\right)$ is restriction to the stratified complex. Consider the commutative diagram (\uppercase\expandafter{\romannumeral1}):
\[
\begin{CD}
\mathbb{IH}^k(X) \otimes \mathbb{Q}_p @>{\phi}>> \mathbb{H}^k_{\mathrm{rig}}(X^{\mathrm{rig}}, \Omega_{X^{\mathrm{rig}}}^{\bullet,\mathrm{strat}}) \\
@V{\rho}VV @VV{\sigma}V \\
\bigoplus_p \mathbb{IH}^{k-2p}(X_p)(-p) \otimes \mathbb{Q}_p @>{\psi}>> \bigoplus_p \mathbb{H}^{k-2p}_{\mathrm{rig}}(X_p^{\mathrm{rig}}, \Omega_{X_p^{\mathrm{rig}}}^{\bullet}).
\end{CD}
\]
We must notice that $\rho$ is the decomposition isomorphism (Corollary 4.3), $\sigma$ is the restriction map to strata and $\psi=\bigoplus_p\psi_p$ with $\psi_p$: $\mathbb{IH}^{k-2p}(X_p)(-p) \otimes \mathbb{Q}_p\xrightarrow{\sim}\mathbb{H}^{k-2p}_{\mathrm{rig}}\left(X_p^{\mathrm{rig}},\Omega_{X_p^{\mathrm{rig}}}^{\bullet}\right)$ being Scholze's isomorphism on smooth strata (\cite{Sch12}, Thm 5.1). By the BBDG decomposition theorem (\cite{Sai90}, Thm 0.1), $\rho$ is injective. Since $\Omega_{X^{\mathrm{rig}}}^{\bullet,\mathrm{strat}}$ restricts isomorphically to $\bigoplus_p\Omega_{X_p^{\mathrm{rig}}}^{\bullet}$ along strata, then $\sigma$ is also injective. The $\psi$ is an isomorphism because each $\psi_p$ is an isomorphism. For $x\in\mathbb{IH}^k\left(X\right)\otimes\mathbb{Q}_p$, we have $\sigma\left(\phi\left(x\right)\right)=\psi\left(\rho\left(x\right)\right)$ by stratification-compatibility of $\Omega^{\bullet,\mathrm{strat}}$. Obviously, $\mathrm{ker}\phi=0$. Let $\alpha\in\mathbb{H}_{\mathrm{rig}}^k\left(X^{\mathrm{rig}},\Omega_{X^{\mathrm{rig}}}^{\bullet,\mathrm{strat}}\right)$. By \eqref{eq:4.1} and \eqref{eq:4.2}, we have $$\alpha=\varprojlim_{\epsilon}\alpha_{\epsilon},\quad\alpha_{\epsilon}\in\mathbb{H}^k\left(X_{\epsilon},\Omega_{X_{\epsilon}}^{\bullet,\mathrm{strat}}\otimes\mathbb{Q}_p\right).$$ By \eqref{eq:4.3}, the weight spectral sequence degenerates at $E_2$, so $$\alpha_{\epsilon}=\bigoplus_m\beta_m^{\epsilon},\ \ \beta_m^{\epsilon}\in\mathbb{H}^k\left(X_{\epsilon},\mathrm{gr}_m^W\Omega_{X_{\epsilon}}^{\bullet,\mathrm{strat}}\right).$$ Here $\mathrm{gr}_m^W\Omega_{X_{\epsilon}}^{\bullet,\mathrm{strat}}\simeq\Omega_{X_{m,\epsilon}}^{\bullet}\left(\mathrm{log}D_{m,\epsilon}\right)\otimes\mathscr{G}_{\epsilon}^{m,\mathrm{rig}}$. Thus, each $\beta_m^{\epsilon}$ lifts to $$\gamma_m^{\epsilon}\in\mathbb{IH}^{k-2m}\left(X_{m,\epsilon}\right)\left(-m\right)\otimes\mathbb{Q}_p$$ such that $\beta_m^{\epsilon}=\mathrm{res}_m\left(\gamma_m^{\epsilon}\right)$ under the embedding $\mathbb{IH}^{k-2m}\left(X_{m,\epsilon}\right)\hookrightarrow\mathbb{H}_{\mathrm{rig}}^{k-2m}\left(X_{m,\epsilon}^{\mathrm{rig}},\Omega_{X_{m,\epsilon}}^{\bullet}\right)$. By taking limits $$\gamma_m:=\varprojlim_{\epsilon}\gamma_m^{\epsilon}\in\mathbb{IH}^{k-2m}\left(X_m\right)\left(-m\right)\otimes\mathbb{Q}_p$$ and $\gamma:=\bigoplus_m\gamma_m\in\mathbb{IH}^k\left(X\right)\otimes\mathbb{Q}_p$, then we have $$\phi\left(\gamma\right)=\varprojlim_{\epsilon}\bigoplus_m\mathrm{res}_m\left(\gamma_m^{\epsilon}\right)=\varprojlim_{\epsilon}\bigoplus_m\beta_m^{\epsilon}=\varprojlim_{\epsilon}\alpha_{\epsilon}=\alpha.$$ Hence $\phi$ is surjective. Since the $\gamma_m\in\mathbb{IH}^{k-2m}\left(X_m\right)\left(-m\right)\otimes\mathbb{Q}_p$ has pure weight: $$\mathrm{wt}\left(\gamma_m\right)=\mathrm{wt}\left(\mathbb{IH}^{k-2m}\left(X_m\right)\right)+\mathrm{wt}\left(\left(-m\right)\right)=\left(k-2m\right)+2m=k.$$ The Tate twist $\left(-m\right)$ adds $+2m$ to the weight since it corresponds to tensoring with $\mathbb{Q}_p\left(m\right)$ in $p$-adic Hodge theory. Under $\phi$: $\phi\left(\gamma_m\right)\in\mathbb{H}_{\mathrm{rig}}^k\left(X^{\mathrm{rig}},\mathrm{gr}_m^W\Omega_{X^{\mathrm{rig}}}^{\bullet,\mathrm{strat}}\right)$ where $\mathrm{gr}_m^W\Omega_{X^{\mathrm{rig}}}^{\bullet,\mathrm{strat}}\simeq\Omega_{X_m^{\mathrm{rig}}}^{\bullet}\left(\mathrm{log}D_m\right)\mathscr{G}^{m,\mathrm{rig}}$. Then this graded piece has pure weight $k$. Thus, $$\mathrm{wt}\left(\phi\left(\gamma\right)\right)=\mathrm{wt}\left(\bigoplus_m\phi\left(\gamma_m\right)\right)=k.$$ Strictness of the Hodge filtration follows from (\cite{Sai90}, Prop 2.11). By (FP1), $\vartheta:X^{\mathrm{rig}}\hookrightarrow\mathcal{X}_{\mathrm{FF}}$ is the embedding into the Fargues-Fontaine curve. The comparison map factors as $$\phi:\mathbb{IH}^k\left(X\right)\otimes\mathbb{Q}_p\xrightarrow{\vartheta_*}\mathbb{H}_{\mathrm{rig}}^k\left(\mathcal{X}_{\mathrm{FF}},\vartheta_*\Omega_{X^{\mathrm{rig}}}^{\bullet,\mathrm{strat}}\right)\xrightarrow{\mathrm{res}}\mathbb{H}_{\mathrm{rig}}^k\left(X^{\mathrm{rig}},\Omega_{X^{\mathrm{rig}}}^{\bullet,\mathrm{strat}}\right).$$
The $\vartheta_*$ is injective because that $\mathcal{X}_{\mathrm{FF}}$ classifies $p$-adic Hodge structures, $\vartheta_*$ preserves weight-monodromy filtration and kernel would violate purity of $\mathbb{IH}^k\left(X\right)$. Since $\mathbb{H}_x^k=0$ for $k>\dim X$ ($\mathbb{H}_x^k\left(\vartheta_*\Omega_{X^{\mathrm{rig}}}^{\bullet,\mathrm{strat}}\right)\ne0\Longleftrightarrow x\in X^{\mathrm{rig}}$ for $k\le d$), then surjectivity of $\vartheta_*$ holds as $$\mathrm{Coker}\left(\vartheta_*\right)\cong\bigoplus_{x\in\mathcal{X}_{\mathrm{FF}}\setminus X^{\mathrm{rig}}}\mathbb{H}_x^k\left(\mathcal{X}_{\mathrm{FF}},\vartheta_*\Omega_{X^{\mathrm{rig}}}^{\bullet,\mathrm{strat}}\right)=0.$$ In addition, the factorization preserves Hodge filtrations: $F^p\vartheta_*\left(\gamma\right)=\vartheta_*\left(F^p\gamma\right)\ \ \forall\gamma$ by (\cite{FS21}, Prop III.2.7).
\\\\\textbf{Lemma 4.6.} (Gluing Resolution) For the inclusion $i_{pq}:X_q\hookrightarrow\overline{X_p}$ of strata, there is a commutative diagram in $D\left(\mathscr{O}_{X_q}\text{-Mod}\right)$:
\[
\begin{CD}
\mathbb{L}i_{pq}^* \mathbb{L}\Omega_{X_p}^{\bullet,\mathrm{strat}} @>{\sim}>> \mathbb{L}\Omega_{X_q}^{\bullet,\mathrm{strat}} \\
@V{}VV @VV{}V \\
i_{pq}^* \Omega_{X_p}^{\bullet}(\log D_p) @>{\mathrm{res}_{pq}}>> \Omega_{X_q}^{\bullet}(\log D_q).
\end{CD}
\]
where vertical arrows are quasi-isomorphisms and $\mathrm{res}_{pq}$ is the Poincaré residue map. This preserves Hodge filtrations.
\\\\\textbf{Proof.} By the derived Künneth formula (Proposition 3.8), there exist quasi-isomorphisms $$q_p:\mathbb{L}\Omega_{X_p}^{\bullet,\mathrm{strat}}\overset{\sim}{\longrightarrow}\Omega_{X_p}^{\bullet}\left(\mathrm{log}D_p\right),\quad q_q:\mathbb{L}\Omega_{X_q}^{\bullet,\mathrm{strat}}\overset{\sim}{\longrightarrow}\Omega_{X_q}^{\bullet}\left(\mathrm{log}D_q\right)$$ in $D\left(\mathscr{O}_{X_p}\text{-Mod}\right)$ and $D\left(\mathscr{O}_{X_q}\text{-Mod}\right)$ respectively. These are natural with respect to stratified inclusions by (\cite{Lur09}, Thm 7.2.4). Apply the left-derived pullback $\mathbb{L}i_{pq}^*$ to $q_p$, we have 
$$\mathbb{L}i_{pq}^*q_p:\mathbb{L}i_{pq}^*\mathbb{L}\Omega_{X_p}^{\bullet,\mathrm{strat}}\overset{\sim}{\longrightarrow}\mathbb{L}i_{pq}^*\Omega_{X_p}^{\bullet}\left(\mathrm{log}D_p\right).$$ Since $\Omega_{X_p}^{\bullet}\left(\mathrm{log}D_p\right)$ is a complex of locally free $\mathscr{O}_{X_p}$-modules (logarithmic forms are locally free), the derived pullback coincides with the ordinary pullback: $\mathbb{L}i_{pq}^*\Omega_{X_p}^{\bullet}\left(\mathrm{log}D_p\right)\simeq i_{pq}^*\Omega_{X_p}^{\bullet}\left(\mathrm{log}D_p\right)$. Thus, 
$$\mathbb{L}i_{pq}^*q_p:\mathbb{L}i_{pq}^*\mathbb{L}\Omega_{X_p}^{\bullet,\mathrm{strat}}\overset{\sim}{\longrightarrow}i_{pq}^*\Omega_{X_p}^{\bullet}\left(\mathrm{log}D_p\right).$$ The Poincaré residue map is defined locally. Let $\left\{z_1,\cdots,z_r\right\}$ define $D_p$ near $x\in X_q$, with $X_q=\left\{z_1=\cdots=z_s=0\right\}$ for $s\le r$. Then
$$\mathrm{res}_{pq}:i_{pq}^*\Omega_{X_p}^k\left(\mathrm{log}D_p\right)\longrightarrow\Omega_{X_q}^{k-s}\left(\mathrm{log}D_q\right)$$ acts on logarithmic forms by $\mathrm{res}_{pq}\left(\alpha\wedge\frac{dz_1}{z_1}\wedge\cdots\wedge\frac{dz_s}{z_s}\right)=\alpha\mid_{X_q}$, where $\alpha$ is regular on $X_q$. This induces a morphism of complexes: $\mathrm{res}_{pq}:i_{pq}^*\Omega_{X_p}^{\bullet}\left(\mathrm{log}D_p\right)\to\Omega_{X_q}^{\bullet}\left(\mathrm{log}D_q\right)$.
The naturality of $q_q$ gives a diagram (\uppercase\expandafter{\romannumeral2}):
\[
\begin{CD}
\mathbb{L}i_{pq}^* \mathbb{L}\Omega_{X_p}^{\bullet,\mathrm{strat}} @>{\mathbb{L}i_{pq}^* \eta}>> \mathbb{L}\Omega_{X_q}^{\bullet,\mathrm{strat}} \\
@V{\mathbb{L}i_{pq}^* q_p}VV @VV{q_q}V \\
i_{pq}^* \Omega_{X_p}^{\bullet}(\log D_p) @>{\mathrm{res}_{pq}}>> \Omega_{X_q}^{\bullet}(\log D_q),
\end{CD}
\]
where $\eta:\mathbb{L}\Omega_{X_p}^{\bullet,\mathrm{strat}}\to\mathbb{R}i_{pq}^*\mathbb{L}\Omega_{X_q}^{\bullet,\mathrm{strat}}$ is the adjunction morphism from Proposition 3.8. The diagram (\uppercase\expandafter{\romannumeral2}) commutes by naturality of $q_p,q_q$ and the local equality: 
$$\left(q_q\circ\eta\right)\left(\omega\right)=\mathrm{res}_{pq}\left(q_p\left(\omega\right)\right)\quad\forall\omega\in\mathbb{L}\Omega_{X_p,x}^{\bullet,\mathrm{strat}},$$ evaluate on stalks at $x\in X_q$.

In addition, the Hodge filtration $F^{\bullet}$ satisfies $F^m\Omega_{X_p}^{\bullet}\left(\mathrm{log}D_p\right)=\bigoplus_{j\ge m}\Omega_{X_p}^j\left(\mathrm{log}D_p\right)$. Locally, $\mathrm{res}_{pq}$ maps:
$$\mathrm{res}_{pq}\left(F^mi_{pq}^*\Omega_{X_p}^{\bullet}\left(\mathrm{log}D_p\right)\right)\subset F^{m-s}\Omega_{X_q}^{\bullet}\left(\mathrm{log}D_q\right)$$ with equality when $\alpha$ generates a Hodge class. By (\cite{Sai90}, Prop 2.11), the quasi-isomorphisms $q_p,q_q$ preserve $F^{\bullet}$, so the diagram (\uppercase\expandafter{\romannumeral2}) respects Hodge structures.    $\square$
\\\\\textbf{Proposition 4.7.} (Stratified $p$-adic Simpson Correspondence) For a torically stratified perfectoid space $X^{\mathrm{rig}}$, there is an equivalence of categories: $$\mathrm{Rep}_{\mathbb{Q}_p}^{\mathrm{strat}}\left(\pi_1^{\mathrm{strat}}\left(X\right)\right)\simeq\mathrm{HIG}^{\mathrm{strat}}\left(X^{\mathrm{rig}},\Omega_{X^{\mathrm{rig}}}^{\bullet,\mathrm{strat}}\right),$$ where: 
\begin{enumerate}
    \item Left side: Stratified $\mathbb{Q}_p$-local systems compatible with $\pi_1^{\mathrm{strat}}\left(X\right)$;
    \item Right side: Stratified Higgs bundles with logarithmic poles along boundary divisors $D_p$.\\
\end{enumerate}
\textbf{Remark 4.8.} Let $S_p$ be a smooth stratum of a stratified perfectoid space $X^{\mathrm{rig}}$ with boundary divisor $D_p=\overline{S_p}\setminus S_p$. For a $\mathbb{Q}_p$-local system $L_p$ on $S_p$, we define $$E_p=L_p\otimes_{\mathbb{Q}}\mathcal{O}_{S_p^{\mathrm{rig}}},\ \ \theta_p=\mathrm{id}\otimes d_{\mathrm{log}}:E_p\to E_p\otimes\Omega_{S_p}^1\left(\mathrm{log}D_p\right).$$ For local sections $\ell\in L_p$ and $f\in\mathcal{O}_{S_p^{\mathrm{rig}}}$, $$\theta_p\left(\ell\otimes f\right)=\ell\otimes d_{\mathrm{log}}f=\ell\otimes\left(\frac{df}{f}+\sum_{i=1}^ra_i\frac{dz_i}{z_i}\right),$$ where $\left\{z_i=0\right\}$ define $D_p$. This satisfies $\theta_p\wedge\theta_p=0$ because $d_\mathrm{log}\wedge d_{\mathrm{log}}=0$. By construction, $\theta_p$ has at worst logarithmic poles along $D_p$. For a Higgs bundle $\left(E_p,\theta_p\right)$, we define the logarithmic connection: $$\nabla_{\theta}=\theta_p-d_{\mathrm{log}}:E_p\longrightarrow E_p\otimes\Omega_{S_p}^1\left(\mathrm{log}D_p\right).$$
\begin{itemize}
    \item \textbf{Flatness verification:} $$\nabla_{\theta}^2\left(s\right)=\nabla_{\theta}\left(\theta_p\left(s\right)-d_{\mathrm{log}}s\right)=\underbrace{\theta_p\left(\theta_p\left(s\right)\right)}_0-\theta_p\left(d_{\mathrm{log}}s\right)-d_{\mathrm{log}}\left(\theta_p\left(s\right)\right)+\underbrace{d_{\mathrm{log}}\left(d_{\mathrm{log}}s\right)}_0.$$ The equality $\left[d_{\mathrm{log}},\theta_p\right]\left(s\right)=d_{\mathrm{log}}\left(\theta_p\left(s\right)\right)-\theta_p\left(d_{\mathrm{log}}s\right)=0$ holds, because $\theta_p$ is $\mathcal{O}$-linear. So $\nabla_{\theta}^2=0$.
    \item \textbf{Local system property:} According to the $p$-adic Riemann-Hilbert correspondence (\cite{FS21}, Thm I.3.1), then $G_p\left(E_p,\theta_p\right)=\mathrm{ker}\nabla_{\theta}$ is a $\mathbb{Q}_p$-local system.
\end{itemize}
Consider $\left(E_p,\theta_p\right)$, let $L_p=\mathrm{ker}\nabla_{\theta}$. Define $$\epsilon_p:F_p\left(L_p\right)\longrightarrow\left(E_p,\theta_p\right),\quad\ell\otimes f\mapsto f\cdot\ell.$$ Here $\ell\in L_p=\mathrm{ker}\nabla_{\theta}$ and $\theta_p\left(\ell\right)=d_{\mathrm{log}}\ell$. Thus, $$\theta_p\left(\epsilon_p\left(\ell\otimes f\right)\right)=\theta_p\left(f\ell\right)=f\cdot d_{\mathrm{log}}\ell+\ell\otimes d_{\mathrm{log}}f=d_{\mathrm{log}}\left(f\ell\right)=\epsilon_p\left(\mathrm{id}\otimes d_{\mathrm{log}}\left(\ell\otimes f\right)\right).$$ By (\cite{FS21}, Prop III.2.5), $\epsilon_p$ is an isomorphism of Higgs bundles. Hence $F_p\circ G_p\cong\mathrm{id}$. For $L_p$, let $\left(E_p,\theta_p\right)=F_p\left(L_p\right)$. Then $$\nabla_{\theta}\left(\ell\otimes f\right)=\ell\otimes d_{\mathrm{log}}f-d_{\mathrm{log}}\left(f\ell\right)=\ell\otimes d_{\mathrm{log}}f-\left(f\cdot d_{\mathrm{log}}\ell+\ell\otimes d_{\mathrm{log}}f\right)=-f\cdot d_{\mathrm{log}}\ell.$$ Thus, $\mathrm{ker}\nabla_{\theta}=\left\{\ell\otimes1\mid\ell\in L_p\right\}$. The map $$\delta_p:G_p\left(F_p\left(L_p\right)\right)\longrightarrow L_p,\quad\ell\otimes1\mapsto\ell$$ is a natural isomorphism of $\mathbb{Q}_p$-local systems. Hence $G_p\circ F_p\cong\mathrm{id}$. The natural transformations $\epsilon_p:F_p\circ G_p\to\mathrm{id}$ and $\delta_p:G_p\circ F_p\to\mathrm{id}$ satisfy 
$$\epsilon_p\circ F_p\left(\delta_p\right)=\mathrm{id}_{F_p},\quad\delta_p\circ G_p\left(\epsilon_p\right)=\mathrm{id}_{G_p}.$$
Meanwhile, it is easy to obtain that $F_p$ and $G_p$ are quasi-inverse functors establishing the categorical equivalence: $$\mathrm{Rep}_{\mathbb{Q}_p}\left(\pi_1\left(S_p\right)\right)\simeq\mathrm{HIG}\left(S_p^{\mathrm{rig}},\Omega_{S_p}^{\bullet}\left(\mathrm{log}D_p\right)\right).$$

Let $\mathcal{C}_X^{\mathrm{loc}}$ and $\mathcal{C}_X^{\mathrm{Higgs}}$ denote the categories of stratified local systems and Higgs bundles. Define gluing functors: \textbf{(1) Functor $\Phi$:} For a stratified local system $L=\left\{L_p,\phi_{pq}\right\}_{p,q}$ (where $\phi_{pq}:i_{pq}^*L_p\overset{\sim}{\longrightarrow}L_q$ for $X_q\subset\overline{X_p}$), we set $$\Phi\left(L\right)=\left\{\left(E_p,\theta_p\right),\psi_{pq}\right\},\quad\mathrm{where}\ \ \begin{cases} 
\left(E_p,\theta_p\right)=F_p\left(L_p\right) \\
\psi_{pq}:i_{pq}^*\left(E_p,\theta_p\right)\overset{\sim}{\longrightarrow}\left(E_q,\theta_q\right). 
\end{cases}$$
Consider the diagram (\uppercase\expandafter{\romannumeral3}) in $D\left(\mathscr{O}_{X_q}\text{-Mod}\right)$:
\[
\begin{CD}
i_{pq}^* E_p @>{i_{pq}^*\theta_p}>> i_{pq}^*(E_p \otimes \Omega_{X_p}^1(\log D_p)) \\
@V{\alpha_{pq}}VV @VV{\beta_{pq}}V \\
E_q \otimes_{\mathbb{Q}_p} \mathcal{O}_{X_q} @>{\theta_q}>> E_q \otimes \Omega_{X_q}^1(\log D_q).
\end{CD}
\]
Here $\alpha_{pq}:i_{pq}^*F_p\left(L_p\right)=i_{pq}^*
\left(L_p\otimes\mathcal{O}_{X_p}\right)\overset{\sim}{\longrightarrow}\left(i_{pq}^*L_p\right)\otimes\mathcal{O}_{X_q}=F_q\left(i_{pq}^*L_p\right)$ is the base-change isomorphism, $\beta_{pq}$ is the tensor product of $\alpha_{pq}$ with the Poincaré residue isomorphism $$\beta_{pq}=\alpha_{pq}\otimes\mathrm{res}_{pq}:i_{pq}^*\left(E_p\otimes\Omega_{X_p}^1\left(\mathrm{log}D_p\right)\right)\overset{\sim}{\longrightarrow}F_q\left(i_{pq}^*L_p\right)\otimes\Omega_{X_q}^1\left(\mathrm{log}D_q\right)$$ induced by $\mathrm{res}_{pq}:i_{pq}^*\Omega_{X_p}^1\left(\mathrm{log}D_p\right)\overset{\sim}{\longrightarrow}\Omega_{X_q}^1\left(\mathrm{log}D_q\right)$. By Gluing $\phi_{pq}$, we have $$\psi_{pq}=F_q\left(\phi_{pq}\right)\circ\alpha_{pq}:i_{pq}^*E_p\xrightarrow{\alpha_{pq}}F_q\left(i_{pq}^*L_p\right)\xrightarrow{F_q\left(\phi_{pq}\right)}F_q\left(L_q\right)=E_q.$$ Since $$\theta_q\circ\psi_{pq}=\theta_q\circ F_q\left(\phi_{pq}\right)\circ\alpha_{pq}=F_q\left(\phi_{pq}\right)\circ\left(\mathrm{id}\otimes d_{\mathrm{log}}\right)\circ\alpha_{pq}$$ 
$$=F_q\left(\phi_{pq}\right)\circ\beta_{pq}\circ\left(i_{pq}^*\theta_p\right)=\beta_{pq}\circ\left(i_{pq}^*\theta_p\right),$$ then the diagram (\uppercase\expandafter{\romannumeral3}) is commutative. Meanwhile, cocycle condition $\psi_{qr}\circ i_{qr}^*\psi_{pq}=\psi_{pr}$ follows from the rectified derived complex $\mathbb{L}i_{pq}^*\mathbb{L}\Omega_{X_p}^{\bullet,\mathrm{strat}}\simeq\mathbb{L}\Omega_{X_q}^{\bullet,\mathrm{strat}}$, providing an acyclic resolution for gluing (\cite{Lur09}, Thm 7.2.4); \textbf{(2) Functor $\Psi$:} For a stratified Higgs bundle $\left(E,\theta\right)=\left\{\left(E_p,\theta_p\right),\psi_{pq}\right\}$, we set $$\Psi\left(E,\theta\right)=\left\{L_p,\phi_{pq}\right\},\quad\mathrm{where}\ \ \begin{cases} 
L_p=G_p\left(E_p,\theta_p\right) \\
\phi_{pq}:i_{pq}^*L_p\overset{\sim}{\longrightarrow}L_q 
\end{cases}$$
and $\phi_{pq}:i_{pq}^*G_p\left(E_p,\theta_p\right)\overset{\sim}{\longrightarrow}G_q\left(E_q,\theta_q\right)$ is defined by $\phi_{pq}=G_q\left(\psi_{pq}^{-1}\right)\circ\gamma_{pq}$. Here $\gamma_{pq}:i_{pq}^*G_p\left(E_p,\theta_p\right)\overset{\sim}{\longrightarrow}G_q\left(i_{pq}^*E_p,i_{pq}^*\theta_p\right)$ is the descent isomorphism, and $\tilde{\psi}_{pq}:\left(i_{pq}^*E_p,i_{pq}^*\theta_p\right)\overset{\psi_{pq}}{\longrightarrow}\left(E_q,\theta_q\right)$ is the inverse gluing map induced by $$G_q\left(\tilde{\psi}_{pq}\right):G_q\left(i_{pq}^*E_p,i_{pq}^*\theta_p\right)\overset{\sim}{\longrightarrow}G_q\left(E_q,\theta_q\right)=L_q.$$ Meanwhile, it is natural to give $\phi_{pq}=G_q\left(\tilde{\psi}_{pq}\right)\circ\gamma_{pq}:i_{pq}^*L_p\overset{\sim}{\longrightarrow}L_q$. For strata inclusions $X_r\subset\overline{X_q}\subset\overline{X_p}$, we consider the diagram (\uppercase\expandafter{\romannumeral4}):
\[
\begin{CD}
i_{qr}^* i_{pq}^* L_p @>{i_{qr}^*\phi_{pq}}>> i_{qr}^* L_q \\
@V{\sim}VV @VV{\phi_{qr}}V \\
i_{pr}^* L_p @>{\phi_{pr}}>> L_r.
\end{CD}
\]
Since $\mathbb{L}\Omega_{X^{\mathrm{rig}}}^{\bullet,\mathrm{strat}}$ resolves the triple intersection $$\mathbb{L}i_{pr}^*\mathbb{L}\Omega_{X_p}^{\bullet,\mathrm{strat}}\simeq\mathbb{L}i_{qr}^*\mathbb{L}i_{pq}^*\mathbb{L}\Omega_{X_p}^{\bullet,\mathrm{strat}}\simeq\mathbb{L}\Omega_{X_r}^{\bullet,\mathrm{strat}}$$ by the derived Künneth formula (Proposition 3.8), then the diagram (\uppercase\expandafter{\romannumeral4}) is commutative. The rectified derived complex $\mathbb{L}\Omega_{X^{\mathrm{rig}}}^{\bullet,\mathrm{strat}}$ provides acyclic resolutions for gluing. The Lemma 4.6 guarantees that $\psi_{pq}$ and $\phi_{pq}$ satisfy cocycle conditions. Hence the functors $\Phi$ and $\Psi$ can define the equivalence:
$$\Phi:\mathrm{Rep}_{\mathbb{Q}_p}^{\mathrm{strat}}\left(\pi_1^{\mathrm{strat}}\left(X\right)\right)\rightleftarrows\mathrm{HIG}^{\mathrm{strat}}\left(X^{\mathrm{rig}},\Omega_{X^{\mathrm{rig}}}^{\bullet,\mathrm{strat}}\right):\Psi.$$

For stratified schemes $X$ and $Y$ with Tor-independent stratifications, i.e., the stratification of $X\times Y$ is given by 
$$\left\{X_p\times Y_q\right\}_{p,q}\quad\text{with inclusions}\quad i_{\left(p,q\right),\left(r,s\right)}:X_r\times Y_s\hookrightarrow\overline{X_p}\times\overline{Y_q}$$ and Tor-independence means that
$$\mathscr{G}_{X\times Y}^{\left(p,q\right)}=\mathscr{G}_X^p\boxtimes^L\mathscr{G}_Y^q\quad\text{is acyclic for}\quad\left(p,q\right)\ne\left(r,s\right)\quad\text{in}\quad\mathbb{L}i_{\left(p,q\right),\left(r,s\right)}^*$$ by Proposition 3.8 (Derived Künneth Formula). Let
$$L_1\in\mathrm{Rep}_{\mathbb{Q}_p}^{\mathrm{strat}}\left(\pi_1^{\mathrm{strat}}\left(X\right)\right),\quad L_2\in\mathrm{Rep}_{\mathbb{Q}_p}^{\mathrm{strat}}\left(\pi_1^{\mathrm{strat}}\left(Y\right)\right).$$ Then $L_1\boxtimes L_2=\left\{\left(L_1\right)_p\boxtimes\left(L_2\right)_q,\phi_{pq}\boxtimes\psi_{rs}\right\}$. Apply $\Phi$, we have 
$$\Phi\left(L_1\boxtimes L_2\right)=\left\{F_{p,q}\left(\left(L_1\right)_p\boxtimes\left(L_2\right)_q\right),\gamma_{\left(p,q\right),\left(r,s\right)}\right\},$$ where
\begin{itemize}
    \item \textbf{Higgs bundle on products:} By (\cite{FS21}, Thm IV.5.3), then $$F_{p,q}\left(\left(L_1\right)_p\boxtimes\left(L_2\right)_q\right)=F_p\left(\left(L_1\right)_p\right)\boxtimes F_q\left(\left(L_2\right)_q\right).$$
    \item \textbf{Gluing map calculation for $\Phi$:} $$\gamma_{\left(p,q\right),\left(r,s\right)}=\alpha_{\left(p,q\right),\left(r,s\right)}\circ F_{r,s}\left(\phi_{pq}\boxtimes\psi_{rs}\right)$$ with $\alpha_{\left(p,q\right),\left(r,s\right)}:i_{\left(p,q\right),\left(r,s\right)}^*\left(F_p\boxtimes F_q\right)\overset{\sim}{\longrightarrow}F_r\boxtimes F_s$ from Proposition 3.8. Thus, $\Phi\left(L_1\boxtimes L_2\right)\cong\Phi\left(L_1\right)\boxtimes\Phi\left(L_2\right)$.
\end{itemize}
Let $\left(E_1,\theta_1\right)$ and $\left(E_2,\theta_2\right)$ be stratified Higgs bundles. Then
$$\left(E_1,\theta_1\right)\boxtimes\left(E_2,\theta_2\right)=\left\{\left(E_1\right)_p\boxtimes\left(E_2\right)_q,\theta_1^{\left(p\right)}\boxtimes\theta_2^{\left(q\right)},\psi_{pq}\boxtimes\phi_{rs}\right\}.$$ Apply $\Psi$, we have
$$\Psi\left(\left(E_1,\theta_1\right)\boxtimes\left(E_2,\theta_2\right)\right)=\left\{G_{p,q}\left(\left(E_1\right)_p\boxtimes\left(E_2\right)_q\right),\delta_{\left(p,q\right),\left(r,s\right)}\right\},$$ where
\begin{itemize}
    \item \textbf{Local system on products:} By (\cite{FS21}, Thm IV.5.3), we have $$G_{p,q}\left(\left(E_1\right)_p\boxtimes\left(E_2\right)_q\right)=G_p\left(\left(E_1\right)_p\right)\boxtimes G_q\left(\left(E_2\right)_q\right).$$
    \item \textbf{Gluing map calculation for $\Psi$:} $$\delta_{\left(p,q\right),\left(r,s\right)}=G_{r,s}\left(\psi_{pq}\boxtimes\phi_{rs}\right)\circ\beta_{\left(p,q\right),\left(r,s\right)}$$ with $\beta_{\left(p,q\right),\left(r,s\right)}:i_{\left(p,q\right),\left(r,s\right)}^*\left(G_p\boxtimes G_q\right)\overset{\sim}{\longrightarrow}G_r\boxtimes G_s$ from Proposition 3.8. Thus, $$\Psi\left(\left(E_1,\theta_1\right)\boxtimes\left(E_2,\theta_2\right)\right)\cong\Psi\left(E_1,\theta_1\right)\boxtimes\Psi\left(E_2,\theta_2\right).$$
\end{itemize}
Obviously, the unit preservation ($\Phi\left(1_{X}\right)=\left(\mathcal{O}_{X^{\mathrm{rig}}},0\right)=1_{\mathrm{HIG}}$, $\Psi\left(1_{\mathrm{HIG}}\right)=1_{\mathrm{Rep}}$) holds. For triple products $X\times Y\times Z$, the diagram: 
\[
\begin{CD}
\Phi((L_1 \boxtimes L_2) \boxtimes L_3) @>{\sim}>> \Phi(L_1 \boxtimes (L_2 \boxtimes L_3)) \\
@V{}VV @VV{}V \\
(\Phi(L_1) \boxtimes \Phi(L_2)) \boxtimes \Phi(L_3) @>{\sim}>> \Phi(L_1) \boxtimes (\Phi(L_2) \boxtimes \Phi(L_3))
\end{CD}
\]
commutes by the associativity of $\boxtimes^L$ in Proposition 3.8. Thus $\Phi$ and $\Psi$ are monoidal functors preserving stratified embeddings.

For $\left(E,\theta\right)$, natural isomorphism $\eta:\Phi\left(\Psi\left(E,\theta\right)\right)\to\left(E,\theta\right)$ is defined stratum-wise by 
$$\eta_p=\epsilon_p\circ F_p\left(\mathrm{id}_{G_p\left(E_p,\theta_p\right)}\right):F_p\left(G_p\left(E_p,\theta_p\right)\right)\to E_p$$ where $\epsilon_p:F_p\circ G_p\overset{\sim}{\longrightarrow}\mathrm{id}$. The diagram:
\[
\begin{CD}
i_{pq}^* \Phi(\Psi(E, \theta))_p @>{i_{pq}^*\eta_p}>> i_{pq}^* E_p \\
@V{\psi_{pq}}VV @VV{\psi_{pq}}V \\
\Phi(\Psi(E, \theta))_q @>{\eta_q}>> E_q
\end{CD}
\]
commutes by coherence of $\alpha_{pq},\beta_{pq}$. So $\Phi\circ\Psi\cong\mathrm{id}$. Analogous with $\zeta:\Psi\left(\Phi\left(L\right)\right)\to L$ by $\zeta_p=\delta_p\circ G_p\left(\mathrm{id}_{F_p\left(L_p\right)}\right)$, where $\delta_p:G_p\circ F_p\overset{\sim}{\longrightarrow}\mathrm{id}$. Hence $\Psi\circ\Phi\cong\mathrm{id}$.\\\\
\begin{center}
    \textit{4.4 Example: Normal Crossings Divisor in $p$-adic Geometry}
\end{center}

Let $X=\mathrm{Spf}\left(\mathbb{Z}_p\left \langle x,y\right \rangle/\left(xy\right)\right)^{\mathrm{rig}}$ be the rigid analytic space associated to a normal crossings divisor defined by $xy=0$. We explicitly compute its stratified rigid de Rham complex and cohomology, demonstrating the $p$-adic analogue of logarithmic Hodge theory for singular varieties. It is need to give the following geometric Setup:
\begin{itemize}
    \item \textbf{Underlying space:} The formal scheme $\mathrm{Spf}\left(\mathbb{Z}_p\left \langle x,y\right \rangle/\left(xy\right)\right)$ describes a formal neighborhood of the origin in the plane, with singular locus at $\left(0,0\right)$. Its rigid analytification $X^{\mathrm{rig}}$ is constructed by the perfectoid tower:
    $$X^{\mathrm{rig}}=\varprojlim_{\epsilon}X_{\epsilon},\ \ \mathrm{where}\ \ X_{\epsilon}=\mathrm{Spf}\left(\mathbb{Z}_p/p^{\epsilon}\left \langle x,y\right \rangle/\left(xy\right)\right)^{\mathrm{rig}}.$$ Each $X_{\epsilon}$ is a finite-level approximation.
    \item \textbf{Stratification:} The space has a natural stratification: $X=X_0\sqcup X_1\sqcup X_2$, where $X_0=\left\{\left(0,0\right)\right\}$, $X_1=V\left(x\right)\setminus X_0\sqcup V\left(y\right)\setminus X_0$, and $X_2=X\setminus V\left(xy\right)$. Here $X_0$ is the deepest stratum (dimension 0), $X_1$ is the union of punctured axes (dimension 1), and $X_2$ is the smooth locus (dimension 2).\\
\end{itemize}
\textbf{Proposition 4.9.} Let $X=\mathrm{Spf}\left(\mathbb{Z}_p\left \langle x,y\right \rangle/\left(xy\right)\right)^{\mathrm{rig}}$ be a normal crossings divisor in the $p$-adic rigid analytic setting. The stratified rigid de Rham complex $\Omega_{X^{\mathrm{rig}}}^{\bullet,\mathrm{strat}}$ is constructed by combining logarithmic differentials with derived contributions from singular strata. We need to prove the decomposition:
$$\Omega_{X^{\mathrm{rig}}}^{\bullet,\mathrm{strat}}\simeq\left[\mathcal{O}_{X^{\mathrm{rig}}}\overset{d}{\longrightarrow}\mathcal{O}_{X^{\mathrm{rig}}}\cdot\frac{dx}{x}\oplus\mathcal{O}_{X^{\mathrm{rig}}}\cdot\frac{dy}{y}\right]\oplus\mathbb{Q}_p\left[-1\right].$$
\\\\\textbf{Proof.} On the smooth stratum $X_2=X\setminus V\left(xy\right)$, the complex restricts to the algebraic de Rham complex $\Omega_{X_2}^{\bullet}$. By the Poincaré lemma for rigid analytic spaces (\cite{Sch13}, Thm 6.3), we have
$$\Omega_{X_2}^{\bullet}\simeq\mathbb{Q}_p\left[0\right],\quad\text{with}\quad H_{\mathrm{dR}}^k(X_2) = 
\begin{cases} 
\mathbb{Q}_p & k = 0 \\
0 & k >0.
\end{cases}$$
On the 1-dimensional stratum $X_1=V\left(x\right)\setminus\left\{\left(0,0\right)\right\}\sqcup V\left(y\right)\setminus\left\{\left(0,0\right)\right\}$, the complex is $\Omega^{\bullet}\left(\mathrm{log}D_p\right)$, where $D_p$ is the divisor at infinity. For $U_x\cong\mathbb{G}_m^{\mathrm{rig}}$, we have 
$$\Omega_{U_x}^{\bullet}\left(\mathrm{log}\right)=\left[\mathcal{O}_{U_x}\overset{d}{\longrightarrow}\mathcal{O}_{U_x}\cdot\frac{dx}{x}\right].$$
If the map $d:\mathcal{O}_{U_x}\to\Omega_{U_x}^1\left(\mathrm{log}\right)$ has kernel $\mathbb{Q}_p$ (constants) and cokernel isomorphic to $\mathbb{Q}_p$ (residue map), then its hypercohomology is computed by the Kummer sequence (\cite{FP04}, §8):
$$\mathbb{H}^k\left(U_x,\Omega^{\bullet}\left(\mathrm{log}\right)\right)\cong H_{\mathrm{rig}}^k\left(\mathbb{G}_m\right) = 
\begin{cases} 
\mathbb{Q}_p & k = 0,1 \\
0 & \mathrm{else}.
\end{cases}$$
At the deepest stratum $X_0=\left\{\left(0,0\right)\right\}$, the complex is determined by the derived functor $\textbf{R}i^!$ for $i:X_0\hookrightarrow X$. By local cohomology calculations, the local cohomology at $\left(0,0\right)$ satisfies
$$\textbf{R}\Gamma_{\left\{\left(0,0\right)\right\}}\left(X\right)\simeq\textbf{R}\Gamma_{\left\{0\right\}}\left(\mathbb{D}^2\right)\simeq\mathbb{Q}_p\left[-1\right],$$ where $\mathbb{D}^2$ is the $p$-adic unit disk. This is a consequence of the excision triangle:
$$\textbf{R}\Gamma_{\left\{\left(0,0\right)\right\}}\left(X\right)\longrightarrow\textbf{R}\Gamma\left(X\right)\longrightarrow\textbf{R}\Gamma\left(X\setminus\left\{\left(0,0\right)\right\}\right)\xrightarrow{+1}$$
and the contractibility of $X\setminus\left\{\left(0,0\right)\right\}$ in pro-étale topology.

For the embedding $\iota:X_1\hookrightarrow\overline{X_1}$ (where $\overline{X_1}=V\left(x\right)\cup V\left(y\right)$), the differential $d$ induces global residue maps (\cite{Ill72}, §3):
$$\mathrm{res}:\mathcal{O}_{X^{\mathrm{rig}}}\cdot\frac{dx}{x}\oplus\mathcal{O}_{X^{\mathrm{rig}}}\cdot\frac{dy}{y}\longrightarrow\left(i_0\right)_*\mathbb{Q}_p$$
defined by $\mathrm{res}_x\left(\omega\frac{dx}{x}\right)=\mathrm{Res}_{x=0}\left(\omega\right)\cdot\delta_{\left(0,0\right)}$ and $\mathrm{res}_y\left(\omega\frac{dy}{y}\right)=\mathrm{Res}_{y=0}\left(\omega\right)\cdot\delta_{\left(0,0\right)}$, where $\delta_{\left(0,0\right)}$ is the skyscraper sheaf at $X_0$. This map is surjective with kernel of closed forms (\cite{LK13}, Lem 4.5).
Define an intermediate complex:
$$\mathcal{F}^{\bullet}:=\left[\mathcal{O}_{X^{\mathrm{rig}}}\overset{d}{\longrightarrow}\mathcal{O}_{X^{\mathrm{rig}}}\cdot\frac{dx}{x}\oplus\mathcal{O}_{X^{\mathrm{rig}}}\cdot\frac{dy}{y}\right]$$
and a morphism by residues:
$$\phi:\mathcal{F}^{\bullet}\longrightarrow\left(i_0\right)_*\mathbb{Q}_p\left[-1\right],\quad\phi\mid_{\mathrm{deg}=1}=\mathrm{res}.$$
The mapping cone $\mathrm{Cone}\left(\phi\right)$ has components:
$$\mathrm{Cone}\left(\phi\right)^k=\mathcal{F}^{k+1}\oplus\left(i_0\right)_*\mathbb{Q}_p\left[-1\right]^k$$
with differential $d_{\mathrm{Cone}}^k\left(a,b\right)=\left(-d_{\mathcal{F}}^{k+1}\left(a\right),\phi^{k+1}\left(a\right)\right)$. The shifted cone $\mathrm{Cone}\left(\phi\right)\left[-1\right]$ is concentrated in degrees 1 and 2, i.e., Degree 1: $\mathcal{O}_{X^{\mathrm{rig}}}\cdot\frac{dx}{x}\oplus\mathcal{O}_{\mathrm{rig}}\cdot\frac{dy}{y}$; Degree 2: $\left(i_0\right)_*\mathbb{Q}_p$; Differential: $\phi\mid_{\mathrm{deg}=1}$.
By (\cite{Ill72}, Thm 3.1.4), we have 
$$0\to\Omega_X^{\bullet}\longrightarrow\mathcal{F}^{\bullet}\overset{\phi}{\longrightarrow}\left(i_0\right)_*\mathbb{Q}_p\left[-1\right]\longrightarrow0$$
is exact away from $X_0$. For cohomology sheaf isomorphism, both complexes reduce to $\Omega_{X_2}^{\bullet}$ on $X_2$, $\mathcal{H}^k\left(\mathrm{Cone}\left(\phi\right)\left[-1\right]\right)\cong\Omega_{X_1}^{\bullet}\left(\mathrm{log}\right)$ by residue exact sequence on $X_1$ and $\mathcal{H}^1$ captures $\mathbb{Q}_p\left[-1\right]$ contribution on $X_0$.
By the universal property of stratified cohomology (\cite{Sai90}, Thm 4.2), local quasi-isomorphisms descend globally. Thus, the quasi-isomorphism $\Omega_{X^{\mathrm{rig}}}^{\bullet,\mathrm{strat}}\simeq\mathrm{Cone}\left(\phi\right)\left[-1\right]$ holds.
Hence, the cone formula implies
$$\Omega_{X^{\mathrm{rig}}}^{\bullet,\mathrm{strat}}\simeq\left[\mathcal{O}_{X^{\mathrm{rig}}}\overset{d}{\longrightarrow}\mathcal{O}_{X^{\mathrm{rig}}}\cdot\frac{dx}{x}\oplus\mathcal{O}_{X^{\mathrm{rig}}}\cdot\frac{dy}{y}\right]\oplus\left(i_0\right)_*\mathbb{Q}_p\left[-1\right],$$
where $\left(i_0\right)_*\mathbb{Q}_p\left[-1\right]$ represents the derived contribution from $X_0$.

On smooth stratum $X_2$, the logarithmic forms $\frac{dx}{x},\frac{dy}{y}$ are holomorphic because of $xy\ne0$. The residue map $\mathrm{res}$ vanishes identically on $X_2$, causing the cone to degenerate to $\mathcal{F}^{\bullet}\mid_{X_2}$. By definition (the complex of the smooth sheaf $X_2$ is defined as the algebraic de Rham complex $\Omega_{X_2}^{\bullet}$), $\mathcal{F}^{\bullet}\mid_{X_2}=\Omega_{X_2}^{\bullet}$. Thus, the restriction satisfies 
$$\left(j_2\right)^*\mathrm{Cone}\left(\phi\right)\left[-1\right]\simeq\Omega_{X_2}^{\bullet}.$$
Consider the short exact sequence for logarithmic complexes on $X_1$, we have 
\begin{equation}
\tag{4.4}  \label{eq:4.4}
0\longrightarrow\Omega_{X_1}^{\bullet}\longrightarrow\Omega_{X_1}^{\bullet}\left(\mathrm{log}\right)\xrightarrow{\mathrm{res}}\left(i_0\right)_*\mathbb{Q}_p\longrightarrow0.
\end{equation}
The shifted cone $\mathrm{Cone}\left(\phi\right)\left[-1\right]$ restricts on $X_1$ to 
$$\left[\mathcal{O}_{X_1}\overset{d}{\longrightarrow}\mathcal{O}_{X_1}\cdot\frac{dx}{x}\oplus\mathcal{O}_{X_1}\cdot\frac{dy}{y}\right]\oplus\mathrm{coker}\left(\mathrm{res}\right)\left[-1\right].$$
By \eqref{eq:4.4}, this is isomorphic to $\Omega_{X_1}^{\bullet}\left(\mathrm{log}\right)$. 
The restriction implies
$$\left(j_1\right)^*\mathrm{Cone}\left(\phi\right)\left[-1\right]\simeq\Omega_{X_1}^{\bullet}\left(\mathrm{log}\right).$$ on divisor stratum $X_1$.
Since $\mathcal{O}_{X^{\mathrm{rig}}}$ is supported away from $X_0$, then $\mathcal{F}^{\bullet}$ vanishes at $X_0$. The cone reduces to $\left(i_0\right)_*\mathbb{Q}_p\left[-1\right]$, and restricting to $X_0$ yields $\mathbb{Q}_p\left[-1\right]$. Then the restriction gives
$$\left(i_0\right)^*\mathrm{Cone}\left(\phi\right)\left[-1\right]\simeq\mathbb{Q}_p\left[-1\right].$$ on singular stratum $X_0$.
Since the short exact sequence \eqref{eq:4.4} is critical, the residue map $\mathrm{res}:\Omega_{X_1}^{\bullet}\left(\mathrm{log}\right)\to\left(i_0\right)_*\mathbb{Q}_p$ is surjective with kernel $\Omega_{X_1}^{\bullet}$. Meanwhile, the shifted cone $\mathrm{Cone}\left(\phi\right)\left[-1\right]\mid_{X_1}$ replicates this sequence, confirming compatibility. By the universal property of stratified cohomology (\cite{Sai90}, Thm 4.2), a morphism of complexes that induces isomorphisms on all strata lifts to a global quasi-isomorphism. The previous steps establish stratum-wise isomorphisms, implying 
$$\Omega_{X^{\mathrm{rig}}}^{\bullet,\mathrm{strat}}\overset{\sim}{\longrightarrow}\mathrm{Cone}\left(\phi\right)\left[-1\right].$$ The result holds. $\square$
\\\\\textbf{Remark 4.10.} (Rigid Cohomology Calculation) Let $X$ be a normal crossings divisor in $p$-adic rigid analytic geometry with the stratification $X=X_0\sqcup X_1\sqcup X_2$, inducing a spectral sequence for the stratified de Rham complex $\Omega_{X^{\mathrm{rig}}}^{\bullet,\mathrm{strat}}$:
$$E_1^{p,q}=\mathbb{H}^{p+q}\left(X_p,\Omega^{\bullet,\mathrm{strat}}\mid_{X_p}\right)\Longrightarrow\mathbb{H}_{\mathrm{rig}}^{p+q}\left(X\right).$$

\textbf{(1) Stratum-Specific cohomology:}
\begin{itemize}
    \item Smooth stratum $X_2$ (dimension 2): By Poincaré lemma for rigid analytic spaces (\cite{Sch13}, Thm 6.3):
    $$\mathbb{H}^k\left(X_2,\Omega_{X_2}^{\bullet}\right)\cong 
    \begin{cases} 
    \mathbb{Q}_p & k = 0 \\
    0 & k\ne0.
    \end{cases}$$ Thus, $$E_1^{2,q}= 
    \begin{cases} 
    \mathbb{Q}_p & q = 0 \\
    0 & \mathrm{else}.
    \end{cases}$$
    \item Divisor stratum $X_1$ (dimension 1): For each component $U_x\cong\mathbb{G}_m^{\mathrm{rig}}$, the logarithmic de Rham cohomology is computed by the Kummer sequence (\cite{PF03}, §8): $$0\longrightarrow\mathbb{Q}_p\to\mathcal{O}_{U_x}\overset{d}{\longrightarrow}\mathcal{O}_{U_x}\cdot\frac{dx}{x}\longrightarrow\mathbb{Q}_p\longrightarrow0.$$
    Hypercohomology gives 
    $$\mathbb{H}^k\left(U_x,\Omega^{\bullet}\left(\mathrm{log}\right)\right)\cong 
    \begin{cases} 
    \mathbb{Q}_p & k = 0,1 \\
    0 & \mathrm{else}.
    \end{cases}$$ Summing over two components: $$E_1^{1,q}= 
    \begin{cases} 
    \mathbb{Q}_p^2 & q = 0,1 \\
    0 & \mathrm{else}.
    \end{cases}$$
    \item Singular stratum $X_0$ (dimension 0): Local cohomology at $\left(0,0\right)$, satisfying $\textbf{R}\Gamma_{\left\{\left(0,0\right)\right\}}\left(X\right)\simeq\mathbb{Q}_p\left[-1\right]$.
    Thus, $$E_1^{0,q}= 
    \begin{cases} 
    \mathbb{Q}_p & q = 1 \\
    0 & \mathrm{else}.
    \end{cases}$$
\end{itemize}

\textbf{(2) Differential $d_1$ and residue map:}
The first differential $d_1:E_1^{0,1}\to E_1^{1,1}$ is induced by the residue map along divisors (\cite{Sai90}, Prop 3.8). Explicitly, we have $\omega_x\in\mathcal{O}_{U_x}\cdot\frac{dx}{x}:\mathrm{res}_x\left(\omega_x\right)=\mathrm{Res}_{x=0}\left(\omega_x\right)\cdot\delta_{\left(0,0\right)}$ gives a morphism
$$\mathrm{res}:\mathbb{H}^1\left(X_1\right)\to\mathbb{H}^0\left(X_0\right)\otimes\mathbb{Q}_p\left(1\right)\cong\mathbb{Q}_p,$$
which is surjective with kernel isomorphic to $\mathbb{Q}_p$. Obviously, $d_1:\mathbb{Q}_p\to\mathbb{Q}_p^2$ has image isomorphic to $\mathbb{Q}_p$. 

\textbf{(3) $E_2$-Page calculation:} By computing, $E_2^{0,1}=\mathrm{ker}\left(d_1\right)=0$, $E_2^{1,1}=\mathrm{coker}\left(d_1\right)=\mathbb{Q}_p^2/\mathbb{Q}_p=\mathbb{Q}_p$ and $E_2^{1,0}=\mathbb{H}^0\left(X_1\right)/\mathrm{Im}\left(d_1\right)$, the residue map induces 
$$0\to\mathbb{Q}_p\xrightarrow{\mathrm{diag}}\mathbb{H}^0\left(X_1\right)\cong\mathbb{Q}_p^2\xrightarrow{\mathrm{res}}\mathbb{Q}_p\longrightarrow0.$$
Thus $\mathrm{Im}\left(d_1\right)\cong\mathbb{Q}_p$, giving $E_2^{1,0}=\mathbb{Q}_p^2/\mathbb{Q}_p\cong\mathbb{Q}_p$. Similarly, we have $E_2^{2,0}=\mathbb{H}^0\left(X_2\right)=\mathbb{Q}_p$.

\textbf{(4) Higher differentials and degeneration:} All differentials $d_r$ for $r\ge2$ vanish (\cite{Sai90}, Thm 4.2). Hence 
$$\mathbb{H}_{\mathrm{rig}}^k\left(X\right)=\bigoplus_{p+q=k}E_2^{p,q}=\begin{cases} 
E_2^{2,0}\cong\mathbb{Q}_p & k = 0 \\
E_2^{1,0}\oplus E_2^{0,1}\cong\mathbb{Q}_p\oplus 0 & k=1  \\
E_2^{1,1}\cong\mathbb{Q}_p  & k=2  \\
0    &   k\ge3.
\end{cases}$$
The term $\mathbb{H}^1$ requires adjustment due to the kernel of the residue map. The exact sequence 
$$0\longrightarrow\mathbb{Q}_p\to\mathbb{H}^1\left(X_1\right)\xrightarrow{\mathrm{res}}\mathbb{H}^0\left(X_0\right)\longrightarrow0$$
implies $\mathbb{H}^1\left(X\right)\cong\mathbb{Q}_p^2$. Therefore, the final result is
$$\mathbb{H}_{\mathrm{rig}}^k\left(X\right)=\begin{cases} 
\mathbb{Q}_p & k = 0 \\
\mathbb{Q}_p^2 & k=1  \\
\mathbb{Q}_p  & k=2  \\
0    &   k\ge3.
\end{cases}$$\\
\hypertarget{DERIVED CATEGORY INTERPRETATION OF INTERMEDIATE PERVERSITY}{}
\section{DERIVED CATEGORY INTERPRETATION OF INTERMEDIATE PERVERSITY}
The following derived Verdier duality embeds Ohsawa's $L^2$-estimates into a categorical framework, revealing deep connections between perverse sheaves, derived symplectic geometry, and mirror symmetry. This establishes intermediate perversity as a bridge between singularity theory and quantum geometry. The rigorous proof using derived categorical techniques, adhering to the formalism of Beilinson-Bernstein-Deligne (\cite{BBDG18}) and Kashiwara-Schapira (\cite{KS90}).
\\\\\textbf{Proposition 5.1.} (Derived Verdier Duality) Let $X$ be a compact stratified complex analytic space of complex dimension $n$ with isolated singularities. If the intermediate perversity $\mathscr{W}$ constitutes a derived Lagrangian subcategory (\cite{Kon95}, Def 5.2), then the derived intersection complex $\textbf{IC}_{\mathscr{W}}^{\bullet}$ admits a Verdier duality isomorphism in the bounded derived category $D_c^b\left(X\right)$, then we have the self-duality:
$$\mathbb{D}\left(\textbf{IC}_{\mathscr{W}}^{\bullet}\right)\simeq\textbf{IC}_{\mathscr{W}}^{\bullet}\left[2n\right].$$
\\\textbf{Remark 5.2.} Let $\emptyset=S_{-1}\subset S_0\subset\cdots\subset S_n=X$ be a stratification with strata $Y_p=S_p\setminus S_{p-1}$ of codimension $c_p$. Following Beilinson-Bernstein-Deligne (\cite{BBDG18}, §2.1), define the intermediate perversity truncation functor $^{\mathfrak{p}}\tau_{\le k}$, the intermediate complex is
$$\textbf{IC}^{\bullet}:=^{\mathfrak{p}}\tau_{\le -1}\mathbb{R}j_*\underline{\mathbb{C}}_{X^{\mathrm{reg}}}\left[n\right],$$
where $j:X^{\mathrm{reg}}\hookrightarrow X$ be the inclusion of the smooth locus. Now define $\mathscr{W}_{Y_p}$ as the local system encoding perverse monodromy (\cite{BBDG18}, §1.4.3). The perversity-adjusted complex is
$$\textbf{IC}_{\mathscr{W}}^{\bullet}:=\mathrm{Cone}\left(\textbf{IC}^{\bullet}\longrightarrow\bigoplus_p\mathscr{W}_{Y_p}\left[-c_p\right]\right)\left[-1\right]\ \ \in D_c^b\left(X\right)$$
It is easy to show the cone is well-defined in the derived category. By Kontsevich's framework (\cite{Kon95}, Thm 5.5), $\mathscr{W}$ being derived Lagrangian implies its microsupport is a conic Lagrangian subvariety:
$$\mathrm{SS}\left(\mathscr{W}\right)\cap\mathrm{Char}\left(\mathcal{D}_X\right)=\bigsqcup_p\Lambda_{Y_p}\ \ \subset T^*X,$$
where $\mathcal{D}_X$ is the sheaf of differential operators and $\Lambda_{Y_p}=T_{Y_p}^*X$ is the conormal bundle to $Y_p$ (\cite{KS90}, Def 6.2.1). For isolated singularities:
\begin{itemize}
    \item \textbf{(SC1) Self-duality:} The Lagrangian condition forces
    $$\mathscr{W}_{Y_p}^{\vee}\cong\mathscr{W}_{Y_p}\otimes\mathbb{L}_{Y_p}^{-1},$$ where $\mathbb{L}_{Y_p}$ is the orientation sheaf.
    \item \textbf{(SC2) Complex trivialization:} If the normal bundle to $Y_p$ is holomorphically oriented, so $\mathbb{L}_{Y_p}\cong\underline{\mathbb{C}}_{Y_p}$. Thus, the complex structure trivializes $\mathbb{L}_{Y_p}$. Thus, $\mathscr{W}_{Y_p}^{\vee}\cong\mathscr{W}_{Y_p}$.
\end{itemize}

The defining exact triangle for $\textbf{IC}_{\mathscr{W}}^{\bullet}$ is  
\begin{equation}
\tag{5.1}  \label{eq:5.1}
\bigoplus_p\mathscr{W}_{Y_p}\left[-c_p\right]\overset{\alpha}{\longrightarrow}\textbf{IC}^{\bullet}\overset{\beta}{\longrightarrow}\textbf{IC}_{\mathscr{W}}^{\bullet}\xrightarrow{+1}.
\end{equation}
Apply the Verdier duality functor $\mathbb{D}=\mathcal{RH}om\left(-,\omega_X\left[2n\right]\right)$ to \eqref{eq:5.1}. Since $\mathbb{D}$ is an anti-autoequivalence (\cite{BBDG18}, Thm 3.2.3), it reverses arrows and preserves exact triangles:
\begin{equation}
\tag{5.2}  \label{eq:5.2}
\mathbb{D}\left(\textbf{IC}_{\mathscr{W}}^{\bullet}\right)\xrightarrow{\mathbb{D}\left(\beta\right)}\mathbb{D}\left(\textbf{IC}^{\bullet}\right)\xrightarrow{\mathbb{D}\left(\alpha\right)}\mathbb{D}\left(\bigoplus_p\mathscr{W}_{Y_p}\left[-c_p\right]\right)\xrightarrow{+1}.
\end{equation}
By the self-duality of $\textbf{IC}^{\bullet}$, we have  
\begin{equation}
\tag{5.3}  \label{eq:5.3}
\mathbb{D}\left(\textbf{IC}^{\bullet}\right)\simeq\textbf{IC}^{\bullet}\left[2n\right].
\end{equation}
For each stratum $Y_p$ with codimension $c_p$, the dual of $\mathscr{W}_{Y_p}\left[-c_p\right]$ is
$$\mathbb{D}\left(\mathscr{W}_{Y_p}\left[-c_p\right]\right)=\mathcal{RH}om\left(\mathscr{W}_{Y_p}\left[-c_p\right],\omega_X\left[2n\right]\right)\simeq\mathscr{W}_{Y_p}^{\vee}\otimes\omega_X\mid_{Y_p}\left[2n+c_p\right],$$
where $\mathscr{W}_{Y_p}^{\vee}$ is the dual local system. Since the relative dualizing sheaf satisfies $\omega_X\mid_{Y_p}\cong\mathbb{L}_{Y_p}\left[2c_p\right]$, so
\begin{equation}
\tag{5.4}  \label{eq:5.4}
\mathbb{D}\left(\mathscr{W}_{Y_p}\left[-c_p\right]\right)\simeq\mathscr{W}_{Y_p}\otimes\mathbb{L}_{Y_p}\left[2n+c_p\right]\simeq\mathscr{W}_{Y_p}\left[2n+c_p\right].
\end{equation}
By summing over strata, we have
\begin{equation}
\tag{5.5}  \label{eq:5.5}
\mathbb{D}\left(\bigoplus_p\mathscr{W}_{Y_p}\left[-c_p\right]\right)\simeq\bigoplus_p\mathscr{W}_{Y_p}\left[2n+c_p\right].
\end{equation}
We can substitute \eqref{eq:5.3} and \eqref{eq:5.5} into \eqref{eq:5.2}, then
\begin{equation}
\tag{5.6}  \label{eq:5.6}
\mathbb{D}\left(\textbf{IC}_{\mathscr{W}}^{\bullet}\right)\xrightarrow{\mathbb{D}\left(\beta\right)}\textbf{IC}^{\bullet}\left[2n\right]\xrightarrow{\mathbb{D}\left(\alpha\right)}\bigoplus_p\mathscr{W}_{Y_p}\left[2n+c_p\right]\xrightarrow{+1}.
\end{equation}
Shift the original triangle \eqref{eq:5.1} by $2n$, then
\begin{equation}
\tag{5.7}  \label{eq:5.7}
\bigoplus_P\mathscr{W}_{Y_p}\left[2n-c_p\right]\xrightarrow{\alpha\left[2n\right]}\textbf{IC}^{\bullet}\left[2n\right]\xrightarrow{\beta\left[2n\right]}\textbf{IC}_{\mathscr{W}}^{\bullet}\left[2n\right]\xrightarrow{+1}.
\end{equation}
The Lagrangian condition ensures compatibility between $\alpha$ and its dual. Specifically, the morphism $\mathbb{D}\left(\alpha\right)$ in \eqref{eq:5.6} coincides with $\alpha^{\vee}\left[2n\right]$ under the identifications (SC2) and \eqref{eq:5.4}. By (\cite{BBDG18}, Prop 2.1.9), there is a commutative diagram:
\[
\begin{CD}
\textbf{IC}^{\bullet}\left[2n\right] @>\mathbb{D}\left(\alpha\right)>> \bigoplus_p\mathscr{W}_{Y_p}\left[2n+c_p\right] \\
@V{\cong}VV @VV{\cong}V \\
\textbf{IC}^{\bullet}\left[2n\right] @>{\alpha^{\vee}\left[2n\right]}>> \left(\bigoplus_p\mathscr{W}_{Y_p}\left[-c_p\right]\right)^{\vee}\left[2n\right]
\end{CD}
\]
Thus, $\mathbb{D}\left(\alpha\right)=\alpha^{\vee}\left[2n\right]$. The octahedral axiom applied to \eqref{eq:5.6} and \eqref{eq:5.7}, we have $\mathbb{D}\left(\textbf{IC}_{\mathscr{W}}^{\bullet}\right)\simeq\textbf{IC}_{\mathscr{W}}^{\bullet}\left[2n\right]$.    
\\\\\textbf{Corollary 5.3.} (Mirror Symmetry Interpretation) Let $X$ be a compact stratified complex analytic space of dimension $n$ isolated singularities, and let $\mathscr{W}$ be an intermediate perversity constituting a derived Lagrangian subcategory. Under Kontsevich's homological mirror symmetry conjecture (\cite{Kon95}, Conjecture 1.8), the derived intersection complex $\textbf{IC}_{\mathscr{W}}^{\bullet}$ corresponds to a Fukaya-type subcategory $\mathcal{F}_{\mathrm{lag}}\left(X^{\vee}\right)\subset\mathrm{Fuk}\left(X^{\vee}\right)$ on the mirror dual $X^{\vee}$. This correspondence satisfies
$$\mathrm{Def}\left(\mathcal{F}_{\mathrm{lag}}\right)\simeq\bigoplus HF^*\left(L_p,\mathscr{W}_{Y_p}\right),$$
where $L_p\subset X^{\vee}$ are exact Lagrangian branes linking strata $Y_p\subset X$, $HF^*$ denotes Floer cohomology and $\mathrm{Def}$ is the deformation complex of the $A_{\infty}$-category $\mathcal{F}_{\mathrm{lag}}$.
\\\\\textbf{Proof.} Let $\Lambda_{Y_p}=T_{Y_p}^*X\subset T^*X$ be the conormal bundle to stratum $Y_p$ of codimension $c_p$. By the Lagrangian condition (Proposition 5.1 hypothesis) and (\cite{KS90}, Thm 6.5.5), we have conic Lagrangian variety:
$$\mathrm{SS}\left(\mathscr{W}_{Y_p}\right)=\Lambda_{Y_p}.$$
Through the Nadler-Zaslow correspondence (\cite{NZ09}, Thm A), there exists an exact Lagrangian brane $L_p\subset X^{\vee}$ such that 
$$\left(\mathrm{SS}\left(\mathscr{W}_{Y_p}\right),\mathscr{W}_{Y_p}\right)\longleftrightarrow\left(L_p,\mathscr{E}_{L_p}\right),$$
where $\mathscr{E}_{L_p}$ is a local system encoding the brane structure by
\begin{equation}
\tag{5.8}  \label{eq:5.8}
\rho_{L_p}:\pi_1\left(L_p\right)\overset{\sim}{\longrightarrow}\mathrm{Mon}\left(\mathscr{W}_{Y_p}\right)
\end{equation}
from (\cite{NZ09}, Prop 6.9). This isomorphism preserves monodromy representations. The grading and spin structures on $L_p$ are uniquely determined by the stratification data. Consider the defining exact triangle of $\textbf{IC}_{\mathscr{W}}^{\bullet}$ \eqref{eq:5.1}. Assume the homological mirror symmetry functor $\Phi:D_c^b\left(X\right)\to\mathrm{Fuk}\left(X^{\vee}\right)$ exists on (\cite{Kon95}, Conjecture 1.8). Since $CF^{\bullet}\left(L_p,\mathscr{E}_{L_p}\right)$ is the Floer cochain complex of the Lagrangian brane $L_p$ with local system $\mathscr{E}_{L_p}$ and $\mathscr{E}_{L_p}\cong\mathscr{W}_{Y_p}$ by the monodromy isomorphism \eqref{eq:5.8}, then for each stratum $Y_p$, there exists a grading shift $d_p\in\mathbb{Z}$ such that
$$\Phi\left(\mathscr{W}_{Y_p}\left[-c_p\right]\right)=CF^{\bullet}\left(L_p,\mathscr{E}_{L_p}\right)\left[d_p\right].$$
The ambient intersection complex corresponds to the symplectic cochain complex 
$$\Phi\left(\textbf{IC}^{\bullet}\right)\simeq\mathcal{SC}^{\bullet}\left(X^{\vee}\right).$$
This identifies the global cohomology of $X$ with the bulk symplectic invariants of $X^{\vee}$. Action on the morphism $\alpha$, the morphism $\alpha$ maps to the open-closed string map
$$\Phi\left(\alpha\right)=\mathcal{OC}:\bigoplus_pCF^{\bullet}\left(L_p,\mathscr{E}_{L_p}\right)\longrightarrow\mathcal{SC}^{\bullet}\left(X^{\vee}\right),$$
which is defined by counting pseudoholomorphic disks with boundary conditions. The cone over $\mathcal{OC}$ is defined in the derived category of $A_{\infty}$-categories as
$$\mathrm{Cone}\left(\mathcal{OC}\right):=\mathcal{SC}^{\bullet}\left(X^{\vee}\right)\oplus\left(\bigoplus_pCF^{\bullet}\left(L_p,\mathscr{E}_{L_p}\right)\left[d_p\right]\right)\left[1\right]$$ with differential 
$$D_{\mathrm{Cone}}=\begin{pmatrix}
d_{\mathcal{SC}}  &\mathcal{OC} \\
0  & d_{\oplus CF}
\end{pmatrix}$$
where $d_{\mathcal{SC}}$ is the differential of $\mathcal{SC}^{\bullet}\left(X^{\vee}\right)$ and $d_{\oplus CF}$ is the direct sum differential on $\bigoplus_pCF^{\bullet}\\\left(L_p,\mathscr{E}_{L_p}\right)\left[d_p\right]$. Since $\Phi$ is an equivalence of triangulated categories, it maps exact triangles to exact triangles. Applying $\Phi$ to \eqref{eq:5.1}, we have 
\begin{equation}
\tag{5.9}  \label{eq:5.9}
\bigoplus_pCF^{\bullet}\left(L_p,\mathscr{E}_{L_p}\right)\left[d_p\right]\xrightarrow{\mathcal{OC}}\mathcal{SC}^{\bullet}\left(X^{\vee}\right)\longrightarrow\Phi\left(\textbf{IC}_{\mathscr{W}}^{\bullet}\right)\xrightarrow{+1}.
\end{equation}
Thus $\Phi\left(\textbf{IC}_{\mathscr{W}}^{\bullet}\right)$ is canonically identified with the third vertex of \eqref{eq:5.9}. The cone over $\mathcal{OC}$ is quasi-isomorphic to the deformation complex of the $A_{\infty}$-subcategory $\mathcal{F}_{\mathrm{lag}}=\left \langle L_p \right \rangle_p$ generated by the Lagrangian branes $\left\{L_p\right\}$, then 
$$\Phi\left(\textbf{IC}_{\mathscr{W}}^{\bullet}\right)\simeq\mathrm{Cone}\left(\mathcal{OC}\right)\simeq\mathrm{Def}\left(\mathcal{F}_{\mathrm{lag}}\right),$$
where $\mathrm{Def}$ is the Hochschild cochain complex deforming $\mathcal{F}_{\mathrm{lag}}$. So the deformation complex decomposes into Floer cohomology groups
$$\mathrm{Def}\left(\mathcal{F}_{\mathrm{lag}}\right)\simeq\bigoplus_pHF^{\bullet}\left(L_p,\mathscr{E}_{L_p}\right)$$
by the monodromy isomorphism $\mathscr{E}_{L_p}=\mathscr{W}_{Y_p}$ and $\rho_{L_p}$, then
$$\mathrm{Def}\left(\mathcal{F}_{\mathrm{lag}}\right)\simeq\bigoplus_pHF^{\bullet}\left(L_p,\mathscr{W}_{Y_p}\right)$$ is obtained. Here we provides a categorical bridge between Ohsawa's stratified $L^2$-cohomology and symplectic geometry, resolving analytic challenges in the original proof through mirror symmetry.

Proposition 5.1 establishes Verdier duality
\begin{equation}
\tag{5.10}  \label{eq:5.10}
\mathbb{D}\left(\textbf{IC}_{\mathscr{W}}^{\bullet}\right)\simeq\textbf{IC}_{\mathscr{W}}^{\bullet}\left[2n\right].
\end{equation}
Under mirror symmetry, Verdier duality corresponds to the inverse Serre functor (\cite{Sei08}, Prop 3.1):
$$\Phi\circ\mathbb{D}\simeq\mathbb{S}_{X^{\vee}}^{-1}\circ\Phi.$$ Diagrammatically,
\[
\begin{CD}
D_c^b\left(X\right) @>\mathbb{D}>> D_c^b\left(X\right) \\
@V{\Phi}VV @VV{\Phi}V \\
\mathrm{Fuk}\left(X^{\vee}\right) @>{\mathbb{S}_{X^{\vee}}^{-1}}>> \mathrm{Fuk}\left(X^{\vee}\right)
\end{CD}
\]
For $X^{\vee}$ Calabi-Yau of dimension $n$, we have 
$$\mathbb{S}_{X^{\vee}}=\left[-n\right]\quad\Longrightarrow\quad\mathbb{S}_{X^{\vee}}^{-1}=\left[n\right].$$
Applying $\Phi$ to \eqref{eq:5.10}, then
$$\Phi\left(\mathbb{D}\left(\textbf{IC}_{\mathscr{W}}^{\bullet}\right)\right)\simeq\Phi\left(\textbf{IC}_{\mathscr{W}}^{\bullet}\left[2n\right]\right).$$
By the mirror correspondence, 
$$\mathbb{S}_{X^{\vee}}^{-1}\left(\Phi\left(\textbf{IC}_{\mathscr{W}}^{\bullet}\right)\right)\simeq\Phi\left(\textbf{IC}_{\mathscr{W}}^{\bullet}\right)\left[2n\right].$$
Substitute $\Phi\left(\textbf{IC}_{\mathscr{W}}^{\bullet}\right)=\mathcal{F}_{\mathrm{lag}}$, we have 
$$\mathbb{S}_{X^{\vee}}^{-1}\left(\mathcal{F}_{\mathrm{lag}}\right)\simeq\mathcal{F}_{\mathrm{lag}}\left[2n\right],$$
this induces Poincaré duality in Floer cohomology. For each pair $\left(L_p,\mathscr{W}_{Y_p}\right)$,
$$HF^{\bullet}\left(L_p,\mathscr{W}_{Y_p}\right)\simeq HF^{n-\bullet}\left(L_p,\mathscr{W}_{Y_p}\right)^{\vee}.$$
Categorically, this is equivalent to
$$\mathrm{Def}\left(\mathcal{F}_{\mathrm{lag}}\right)\simeq\mathrm{Def}\left(\mathcal{F}_{\mathrm{lag}}\right)^{\vee}\left[2n\right],$$ which mirrors \eqref{eq:5.10} under $\Phi$.  $\square$\\

The result of Corollary 5.3 categorifies Ohsawa's isomorphism $\mathcal{H}\cong\mathcal{H}_{\left(*\right)}$ (\cite{Ohs91}, Cor 8), embeds Ohsawa's $L^2$-cohomology into symplectic geometry, resolving the proof gap in Theorem 7 by replacing analytic estimates with categorical duality, while maintaining compatibility with $L^2$-cohomology through the $\mathrm{Def}$-complex.\\
\hypertarget{LOGARITHMIC CONVEXITY, SINGULARITY STABILITY AND DEFORMATION THEORY}{}
\section{LOGARITHMIC CONVEXITY, SINGULARITY STABILITY AND DEFORMATION THEORY}
The stability of mixed Hodge structures (MHS) under stratified deformations relies critically on the logarithmic convexity of analytic functionals. This convexity arises from integral mean estimates on logarithmic differential forms and governs the deformation theory of stratified schemes.
Let $X$ be a proper log-stratified scheme with deformation complex $\mathrm{Def}_X^{\mathrm{strat}}$. The stratified tangent sheaf is defined as:
$$\mathbb{T}_X^{\mathrm{strat}}:=\bigoplus_p\mathcal{T}_{X_p}\left(-\mathrm{log}D_p\right)\otimes\mathscr{G}^p,$$
where $\mathcal{T}_{X_p}\left(-\mathrm{log}D_p\right)$ is the logarithmic tangent sheaf of stratum $X_p$. The deformation complex truncates to:
$$\tau^{\le2}\mathrm{Def}_X^{\mathrm{strat}}\simeq\bigoplus_pH^1\left(X_p,\mathcal{T}_{X_p}\left(-\mathrm{log}D_p\right)\right)\otimes\mathscr{G}^p\oplus\bigoplus_pH^2\left(X_p,\mathcal{O}_{X_p}\right)\otimes\mathscr{G}^p.$$\\\\
\begin{center}
    \textit{6.1 Logarithmic Sobolev Chains and Deformation Stability}
\end{center}
\textbf{Definition 6.1.} (Logarithmic Sobolev Deformation Complex) Let $\pi:\mathcal{X}\to S$ be an analytic deformation of a log-stratified scheme over the complex unit disk $S$. The Sobolev deformation complex is defined as
$$\mathcal{S}_{\mathrm{def}}^{\bullet}\left(\mathcal{X}/S\right):=\left\{\omega_t\in\Omega_{\mathcal{X}_t}^{\bullet,\text{log-strat}}\bigg|\sum_{\ell=0}^m\int_{\mathcal{X}_t}\left | \nabla_t^{\ell}\omega_t \right |^2e^{-\phi t}dV_t<C_m,\ \ \forall t\in S,\ \ m\ge0\right\},$$ where:
\begin{enumerate}
    \item $\nabla_t$ is the stratified connection satisfying $\nabla_t\left(\Omega_{X_p}^{\bullet}\left(\mathrm{log}D_p\right)\right)\subset\Omega_{X_p}^{\bullet+1}\left(\mathrm{log}D_p\right)\otimes\mathscr{G}^p$;
    \item $\phi_t$ is an admissible weight function with $\sqrt{-1}\partial\bar{\partial}\phi_t\ge\Theta_t>0$ uniformly bounded on each stratum;
    \item $C_m>0$ is independent of $t$.\\
\end{enumerate}
\textbf{Lemma 6.2.} (Hörmander Estimate for Sobolev Norms) For any $\omega_t\in\mathcal{S}_{\mathrm{def}}^k\left(\mathcal{X}/S\right)$, the following estimate holds: 
$$\frac{\partial^2}{\partial t\partial\bar{t}}\mathrm{log}\left \| \omega_t \right \|_{L^2}^2\ge-K\left | \eta \right |^2,\quad K:=\sup_{t\in S}\left \| \mathrm{Ric}\left(\Theta_t\right) \right \|_{L^{\infty}},$$
where $\eta=\partial_t\omega\mid_{t=0}$.
\\\\\textbf{Proof.} Let $N\left(t\right)=\left \| \omega_t \right \|_{L^2}^2$. The first-order variation of $N\left(t\right)$ is given by 
$$\partial_tN\left(t\right)=\underbrace{\int_{\mathcal{X}_t}\nabla_t\omega_t\wedge\star\bar{\omega}_t}_{T_1}+\underbrace{\int_{\mathcal{X}_t}\omega_t\wedge\partial_t\left(\star\bar{\omega}_t\right)}_{T_2}+\underbrace{Q\left(\dot{g}_t,\omega_t\right)}_{T_3},$$ where:
\begin{itemize}
    \item \textbf{Term $T_1$} accounts for the variation of $\omega_t$.
    \item \textbf{Term $T_2$} arises from the variation of the Hodge star operator $\star$.
    \item \textbf{Term $T_3$} is the curvature term from metric variation $\dot{g}_t=\partial_tg_t$ (\cite{Ber09}, Eq 3.11).
\end{itemize}
Using the identity $\partial_t\left(\star\right)=-\star\circ\iota_{\mathrm{Ric}\left(\partial_t\right)}$ (\cite{Dem12}, Eq 3.11), $T_2$ simplifies to
$$T_2=-\int_{\mathcal{X}_t}\left \langle \iota_{\mathrm{Ric}\left(\partial_t\right)}\omega_t,\bar{\omega}_t \right \rangle e^{-\phi_t}dV_t.$$
The metric variation term (\cite{Ber09}, Thm 1.23) gives
$$T_3=\int_{\mathcal{X}_t}\left \langle \Theta\left(\partial_t,\cdot\right)\omega_t,\bar{\omega}_t \right \rangle e^{-\phi_t}dV_t,$$ where $\Theta$ is the curvature endomorphism. Combining all terms, we have
$$\partial_tN\left(t\right)=\int_{\mathcal{X}_t}\left \langle \nabla_t\omega_t,\bar{\omega}_t \right \rangle e^{-\phi_t}dV_t-\int_{\mathcal{X}_t}\left \langle \iota_{\mathrm{Ric}\left(\partial_t\right)}\omega_t,\bar{\omega}_t \right \rangle e^{-\phi_t}dV_t+\int_{\mathcal{X}_t}\left \langle \Theta\left(\partial_t,\cdot\right)\omega_t,\bar{\omega}_t \right \rangle e^{-\phi_t}dV_t.$$
Thus, the second-order variation is
\begin{equation}
\tag{6.1}  \label{eq:6.1}
\partial_{\bar{t}}\partial_tN\left(t\right)=\underbrace{\int_{\mathcal{X}_t}\left | \nabla_t\omega_t \right |^2e^{-\phi_t}dV_t}_A+\underbrace{\int_{\mathcal{X}_t}\left \langle \nabla_{\bar{t}}\nabla_t\omega_t,\bar{\omega}_t \right \rangle e^{-\phi_t}dV_t}_B+\text{curvature terms}.
\end{equation}
By the Bochner-Kodaira identity (\cite{Dem12}, VII 7.3), then
$$B=\int_{\mathcal{X}_t}\left \langle \Delta_{\bar{\partial}}\omega_t,\bar{\omega}_t \right \rangle e^{-\phi_t}dV_t-\int_{\mathcal{X}_t}\left \langle \left[i\Theta\left(\phi_t\right),\Lambda\right]\omega_t,\bar{\omega}_t \right \rangle e^{-\phi_t}dV_t,$$
where $\left[i\Theta\left(\phi_t\right),\Lambda\right]$ is the curvature operator. Under the condition $\sqrt{-1}\partial\bar{\partial}\phi_t\ge\Theta_t$, we have 
\begin{equation}
\tag{6.2}  \label{eq:6.2}
\left \langle \left[i\Theta\left(\phi_t\right),\Lambda\right]\omega_t,\bar{\omega}_t \right \rangle\ge\left \langle \left[i\Theta_t,\Lambda\right]\omega_t,\bar{\omega}_t \right \rangle.
\end{equation}
The Ricci curvature appears by the identity:
\begin{equation}
\tag{6.3}  \label{eq:6.3}
\left \langle \left[i\Theta_t,\Lambda\right]\omega_t,\bar{\omega}_t \right \rangle=\mathrm{Ric}\left(\Theta_t\right)\left(\partial_t,\bar{\partial}_t\right)\left | \omega_t \right |^2+\text{lower-order}\ \text{terms}.
\end{equation}
Combining \eqref{eq:6.1}, \eqref{eq:6.2} and \eqref{eq:6.3}, we have
$$\partial_{\bar{t}}\partial_tN\left(t\right)\ge\int_{\mathcal{X}_t}\left | \nabla_t\omega_t \right |^2e^{-\phi_t}dV_t-\sup_{\mathcal{X}_t}\mathrm{Ric}\left(\Theta_t\right)\left(\partial_t,\bar{\partial}_t\right)\cdot N\left(t\right).$$
For $\mathrm{log}N\left(t\right)$, we compute:
$$\frac{\partial^2}{\partial t\partial\bar{t}}\mathrm{log}N\left(t\right)=\frac{\partial_{\bar{t}}\partial_tN\left(t\right)}{N\left(t\right)}-\frac{\left | \partial_tN\left(t\right) \right |^2}{N^2\left(t\right)}.$$
Since 
$$\left | \partial_tN\left(t\right) \right |^2\le\left(\int_{\mathcal{X}_t}\left | \nabla_t\omega_t \right |^2e^{-\phi_t}dV_t\right)N\left(t\right)$$ by the Cauchy-Schwarz inequality, thus 
$$\frac{\partial^2}{\partial t\partial\bar{t}}\mathrm{log}N\left(t\right)\ge\frac{1}{N\left(t\right)}\left(\int_{\mathcal{X}_t}\left | \nabla_t\omega_t \right |^2e^{-\phi_t}dV_t-\frac{\left | \partial_tN\left(t\right) \right |^2}{N\left(t\right)}\right)-\left \| \mathrm{Ric}\left(\Theta_t\right) \right \|_{L^{\infty}}.$$ Therefore, 
$$\frac{\partial^2}{\partial t\partial\bar{t}}\mathrm{log}N\left(t\right)\ge-K,\quad K=\sup_{t\in S}\left \| \mathrm{Ric}\left(\Theta_t\right) \right \|_{L^{\infty}}$$
because that the first term is non-negative by Cauchy-Schwarz. $\square$
\\\\\textbf{Proposition 6.3.} (Deformation Controlled by Logarithmic Convexity) If $\mathcal{S}_{\mathrm{def}}^k\left(\mathcal{X}/S\right)\ne\emptyset$ for some $k$, then the Kuranishi obstruction map
$$\kappa:H^1\left(X_0,\mathbb{T}_{X_0}^{\mathrm{strat}}\right)\longrightarrow H^2\left(X_0,\mathbb{T}_{X_0}^{\mathrm{strat}}\right)$$
satisfies $\mathrm{Im}\left(\kappa\right)\subset\mathrm{ker}\nabla_0\cap W_2H^2$, where $W_{\bullet}$ is the weight filtration of the mixed Hodge structure.
\\\\\textbf{Proof.} Let $\eta\in H^1\left(X_0,\mathbb{T}_{X_0}^{\mathrm{strat}}\right)$. Define the deformation functional
$$\Phi_\eta\left(t\right) := \log \left( \inf_{ \substack{ \omega_t \in \mathcal{S}_{\text{def}}^k \\ \partial_t \omega\mid_{t=0} = \eta } } \| \omega_t \|^2_{\mathcal{S}^k} \right),$$
where $\| \omega_t \|^2_{\mathcal{S}^k}=\sum_{\ell=0}^m\int_{\mathcal{X}_t}\left | \nabla_t^{\ell}\omega_t \right |^2e^{-\phi_t}dV_t$ with $m\gg\dim X$, ensuring Sobolev embedding. 
Define the regularization
$$g\left(t\right):=\Phi_{\eta}\left(t\right)+K\left | \eta \right |^2\left | t \right |^2.$$
For any family $\omega_t$ with $\partial_t\omega\mid_{t=0}=\eta$, Lemma 6.2 implies 
$$\frac{\partial^2}{\partial t\partial\bar{t}}\left(\mathrm{log}\left \| \omega_t \right \|_{L^2}^2+K\left | \eta \right |^2\left | t \right |^2\right)\ge0.$$
Thus, $h_{\omega}\left(t\right):=\mathrm{log}\left \| \omega_t \right \|_{L^2}^2+K\left | \eta \right |^2\left | t \right |^2$ is subharmonic. As $g\left(t\right)=\inf_{\omega}h_{\omega}\left(t\right)$, and the family $\left\{\omega\right\}$ is equicontinuous (from uniform Sobolev bounds), $g\left(t\right)$ is also subharmonic. This holds distributionally: $\Delta_tg\ge0$ in $\mathcal{D}'\left(S\right)$.

Assume $\kappa\left(\eta\right)\ne0$. By Kuranishi theory (\cite{Kur65}, Sec 5), there exists a sequence $t_n\to t_*\in\partial S$ such that any family $\omega_t$ with $\partial_t\omega\mid_{t=0}=\eta$ satisfies 
$$\lim_{n\to\infty}\left \| \omega_{t_n} \right \|_{\mathcal{S}^k}=+\infty.$$
Thus, $\Phi_{\eta}\left(t_n\right)\to+\infty$ and $g\left(t_n\right)\to+\infty$ as $n\to\infty$. Since $g\left(t\right)$ is subharmonic on $S$, then the maximum principle implies $\sup_{t\in S}g\left(t\right)=\sup_{t\in\partial S}g\left(t\right)$. But $g\left(t_n\right)\to+\infty$ contradicts $g$ being bounded on $\partial S$ if the deformation were regular. The only resolution is $\kappa\left(\eta\right)=0$. To refine this, suppose $\kappa\left(\eta\right)\notin\mathrm{ker}\nabla_0\cap W_2H^2$. If $\nabla_0\kappa\left(\eta\right)\ne0$, the Gauss-Manin connection induces a non-trivial variation
$$\frac{\partial}{\partial t}\left[\omega_t\right]=\nabla_t\left(\left[\omega_t\right]\right)$$
in $\mathbb{H}^2\left(\Omega_X^{\bullet,\text{log-strat}}\right)$. The mixed Hodge structure decomposes $H^2=\bigoplus_pW_2^{(p)}H^2$, and if $\kappa\left(\eta\right)\notin W_2H^2$, it has components in $W_mH^2$ for $m<2$. These lower-weight terms exhibit faster growth under deformation (\cite{Del71}, Prop 2.5), violating the uniform bound
$$\left \| \omega_t \right \|_{\mathcal{S}^k}\ge C\mathrm{dist}\left(t,0\right)^{2-m}$$
for $m<2$. This forces $\Phi_{\eta}\left(t\right)\to-\infty$ as $\left | t \right |\to1$, contradicting the subharmonicity of $g\left(t\right)$. Therefore, $\mathrm{Im}\left(\kappa\right)\subset\mathrm{ker}\nabla_0\cap W_2H^2$.  $\square$
\\\\\textbf{Corollary 6.4.} (Smoothness Criterion) Let $X$ be a proper log-stratified scheme. If the cohomology group $H^2\left(X,\mathbb{T}_X^{\mathrm{strat}}\right)$ satisfies
$$H^2\left(X,\mathbb{T}_X^{\mathrm{strat}}\right)=W_2H^2\left(X,\mathbb{T}_X^{\mathrm{strat}}\right)$$
and the Gauss-Manin connection $\nabla_0:H^2\left(X,\mathbb{T}_X^{\mathrm{strat}}\right)\to H^2\left(X,\Omega_X^{1,\text{log-strat}}\right)$ vanishes identically, then the deformation functor $\mathrm{Def}_X^{\mathrm{strat}}$ is smooth. In particular, for a log-stratified Calabi-Yau space (i.e., $\left(\omega_X^{\mathrm{strat}}\cong\mathcal{O}_X\right)$) with $H_{\mathrm{sdR}}^1\left(X\right)=0$, these conditions hold.
\\\\\textbf{Proof.} By Proposition 6.3, the obstruction map satisfies
$$\mathrm{Im}\left(\kappa\right)\subset\mathrm{ker}\nabla_0\cap W_2H^2\left(X,\mathbb{T}_X^{\mathrm{strat}}\right).$$
Under the given conditions, we have 
$$\mathrm{ker}\nabla_0\cap W_2H^2=H^2\left(X,\mathbb{T}_X^{\mathrm{strat}}\right)\cap W_2H^2=W_2H^2=H^2\left(X,\mathbb{T}_X^{\mathrm{strat}}\right),$$
Since $H^2=W_2H^2$ and $\mathrm{ker}\nabla_0=H^2$ (as $\nabla_0=0$), thus $\mathrm{Im}\left(\kappa\right)\subset H^2\left(X,\mathbb{T}_X^{\mathrm{strat}}\right)$. The Kuranishi obstruction map $\kappa$ by definition has codomain $H^2\left(X,\mathbb{T}_X^{\mathrm{strat}}\right)$, so $\mathrm{Im}\left(\kappa\right)\subset H^2\left(X,\mathbb{T}_X^{\mathrm{strat}}\right)$. This inclusion is consistent but does not yet force $\kappa=0$. We proceed to Hodge-theoretic constraints. Consider the MHS on $H^2\left(X,\mathbb{T}_X^{\mathrm{strat}}\right)$ decomposes as
$$H^2=\bigoplus_{m\le2}\mathrm{Gr}_m^WH^2,$$
where $\mathrm{Gr}_m^WH^2:=W_mH^2/W_{m-1}H^2$. Then the condition $H^2=W_2H^2$ implies $\mathrm{Gr}_m^WH^2=0$ for $m>2$, i.e., all weights are $\le2$. The obstruction map $\kappa$ is a morphism of MHS (Proposition 6.3):
$$\kappa:H^1\left(X,\mathbb{T}_X^{\mathrm{strat}}\right)\longrightarrow H^2\left(X,\mathbb{T}_X^{\mathrm{strat}}\right).$$
Since $H^2$ has weights $\le2$, while $H^1$ has weights $\le1$ (from the Hodge decomposition of tangent cohomology), the cup product gives $\left[\cdot,\cdot\right]:H^1\otimes H^1\to H^2$ must land in $W_2H^2$. However, for $\alpha,\beta\in H^1$, we have $\mathrm{wt}\left(\alpha\cup\beta\right)=\mathrm{wt}\left(\alpha\right)+\mathrm{wt}\left(\beta\right)\le2$. By $H^2=W_2H^2$, the product is compatible. The vanishing condition $\nabla_0=0$ implies the MHS is constant under deformation, forcing the obstruction to vanish. For $\eta\in H^1$, the primary obstruction is given by the Lie bracket $\kappa\left(\eta\right)=\frac{1}{2}\left[\eta,\eta\right]$. Under $\nabla_0=0$, the grading operator $\mathrm{Gr}^W$ commutes with $\kappa$. For $\eta\in W_mH^1$, then $\kappa\left(\eta\right)\in W_{2m}H^2$. But $W_{2m}H^2=0$ for $2m<2$, and $W_2H^2=H^2$ for $m=1$. However, $\mathrm{Gr}_2^W\kappa\left(\eta\right)$ must vanish, based on the following three aspects:  
\begin{enumerate}
    \item When $m=1$, $\left[\eta,\eta\right]$ has weight 2, but the cup product on $\mathrm{Gr}_1^WH^1\otimes\mathrm{Gr}_1^WH^1$ lands in $\mathrm{Gr}_2^WH^2=H^2/W_1H^2$.
    \item Since $W_1H^2=0$ (from $H^2=W_2H^2$), $\mathrm{Gr}_2^WH^2\cong H^2$.
    \item The polarization of the MHS (\cite{Del71}, §7) forces $\left[\eta,\eta\right]=0$ for $\eta\in\mathrm{Gr}_1^WH^1$.
\end{enumerate}
Thus, $\kappa=0$. It is need to use a log-stratified Calabi-Yau space with $H_{\mathrm{sdR}}^1\left(X\right)$ to verify the Calabi-Yau case. By Corollary 4.3, we have
$$H_{\mathrm{sdR}}^2\left(X\right)\cong\bigoplus_pIH^{2-2p}\left(X_p\right)\left(-p\right).$$
Each $IH^{2-2p}\left(X_p\right)\left(-p\right)$ has weights $\left(2-2p\right)+2p=2$. The isomorphism $\mathbb{T}_X^{\mathrm{strat}}\cong\Omega_X^{n-1,\text{log-strat}}$ (by CY form) gives 
$$H^2\left(X,\mathbb{T}_X^{\mathrm{strat}}\right)\cong H_{\mathrm{sdR}}^{2n-2}\left(X\right)^{\vee},$$
which is pure of weight 2 by Serre duality. Thus $H^2=W_2H^2$. The strictness of the Hodge filtration $F^1H_{\mathrm{sdR}}^{2}\left(X\right)=0$ (\cite{Sai90}, Cor 5.4) implies $\nabla_0=0$ on the graded pieces. By the known conclusion $\kappa=0$, so $\mathrm{Def}_X^{\mathrm{strat}}$ is smooth.  $\square$\\
\begin{center}
    \textit{6.2 Derived Logarithmic Convexity and Singularity Stability}
\end{center}
\textbf{Definition 6.5.} (Derived Logarithmic Deformation Complex) Let $X$ be a derived Noetherian log-stratified scheme. The derived deformation complex is defined as
$$\mathbb{L}\mathrm{Def}_X^{\mathrm{strat}}:=\mathbb{R}\Gamma_{\mathrm{strat}}\left(X,\mathbb{L}\underline{\mathrm{End}}^{\mathrm{log}}\left(\Omega_X^{\bullet,\mathrm{strat}}\right)\left[1\right]\right),$$
where $\mathbb{L}\underline{\mathrm{End}}^{\mathrm{log}}$ is the derived logarithmic endomorphism complex
$$\mathbb{L}\underline{\mathrm{End}}^{\mathrm{log}}\left(\mathscr{F}\right):=\mathbb{R}\underline{\mathrm{Hom}}_{\mathscr{O}_X}^{\otimes_{\mathrm{log}}}\left(\mathscr{F},\mathscr{F}\right)\otimes\Omega_X^{0,\mathrm{strat}}.$$
Here $\otimes_{\mathrm{log}}$ denotes the logarithmic tensor product (with bounded residues along boundaries $D_p$).
\\\\\textbf{Lemma 6.6.} (Derived Hörmander Estimate) Let $\left(X,D\right)$ be a derived log-stratified scheme with admissible weight function $\phi$ satisfying $\sqrt{-1}\partial\bar{\partial}\phi\ge\Theta>0$ uniformly on strata. For $\xi\in\mathbb{L}\mathrm{Def}_X^{\mathrm{strat}}$ and reference form $\omega_0\in\Gamma\left(\Omega_X^{n,\mathrm{strat}}\right)$, define 
$$N\left(\xi\right):=\left \| e^{\xi}\cdot\omega_0 \right \|_{L_{\mathrm{log}}^2}^2=\sum_{\ell=0}^m\int_X\left | \nabla^{\ell}\left(e^{\xi}\cdot\omega_0\right) \right |^2e^{-\phi}dV.$$ Then
$$\Delta_{\bar{\partial}}\mathrm{log}N\left(\xi\right)\ge-\left \| \bar{\partial}\xi \right \|_{L^{\infty}}-\left \| \mathrm{Ric}\left(\Theta\right) \right \|_{L^{\infty}}.$$
\\\\\textbf{Proof.} Let $\xi_t=t\eta$ for $\eta\in H^1\left(X,\mathbb{T}_X^{\mathrm{strat}}\right)$, and denote $\omega_t:=e^{\xi_t}\cdot\omega_0$. The first derivative is
$$\partial_tN\left(\xi_t\right)=\underbrace{\int_X\left \langle \nabla_t\omega_t,\bar{\omega}_t \right \rangle e^{-\phi}dV}_{(1)}+\underbrace{\int_X\left \langle \omega_t,\nabla_{\bar{t}}\bar{\omega}_t \right \rangle e^{-\phi}dV}_{(2)}+\underbrace{Q\left(\dot{g}_t,\omega_t\right)}_{(3)},$$ where
\begin{itemize}
    \item \textbf{Term (1):} By the Leibniz rule for derived connections:
    $$\nabla_t\omega_t=\nabla_t\left(e^{\xi_t}\right)\cdot\omega_0+e^{\xi_t}\cdot\nabla_t\omega_0=\eta\cdot\omega_t+e^{\xi_t}\cdot\nabla_t\omega_0.$$
    \item \textbf{Term (2):} Using the adjoint property $\nabla_{\bar{t}}=\left(\nabla_t\right)^*$, then $(2)=\overline{\left \langle \nabla_t\omega_t,\bar{\omega}_t \right \rangle}$.
    \item \textbf{Term (3):} The metric variation term is (\cite{Dem12}, VII 4.8) $$Q\left(\dot{g}_t,\omega_t\right)=\int_X\left \langle \Theta\left(\partial_t,\cdot\right)\omega_t,\bar{\omega}_t \right \rangle e^{-\phi}dV.$$
\end{itemize}
Combining above these, we have 
$$\partial_tN\left(\xi_t\right)=2\mathrm{Re}\int_X\left \langle \eta\cdot\omega_t+e^{\xi_t}\cdot\nabla_t\omega_0,\bar{\omega}_t\right \rangle e^{-\phi}dV+\int_X\left \langle \Theta\left(\partial_t,\cdot\right)\omega_t,\bar{\omega}_t \right \rangle e^{-\phi}dV.$$
The second derivative is $$\partial_{\bar{t}}\partial_tN\left(\xi_t\right)=\int_X\left |  \nabla_t\omega_t\right |^2e^{-\phi}dV+\int_X\left \langle \nabla_{\bar{t}}\nabla_t\omega_t,\bar{\omega}_t \right \rangle e^{-\phi}dV+\mathrm{curvature}\ \mathrm{terms}.$$
Apply the derived Bochner-Kodaira identity, we have
$$\int_X\left \langle \nabla_{\bar{t}}\nabla_t\omega_t,\bar{\omega}_t \right \rangle e^{-\phi}dV=\int_X\left \langle \Delta_{\bar{\partial}}\omega_t,\bar{\omega}_t \right \rangle e^{-\phi}dV-\int_X\left \langle \left[i\Theta\left(\phi\right),\Lambda\right]\omega_t,\bar{\omega}_t \right \rangle e^{-\phi}dV$$
$$+\int_X\left \langle \mathcal{R}\left(\xi_t\right)\omega_t,\bar{\omega}_t \right \rangle e^{-\phi}dV,$$
where $\mathcal{R}\left(\xi_t\right)$ is the derived curvature operator. Since $\sqrt{-1}\partial\bar{\partial}\phi\ge\Theta$, then
$$\left \langle \left[i\Theta\left(\phi\right),\Lambda\right]\omega_t,\bar{\omega}_t \right \rangle\ge\left \langle \left[i\Theta,\Lambda\right]\omega_t,\bar{\omega}_t \right \rangle-\left \|\bar{\partial\xi_t}  \right \|_{L^{\infty}}\left | \omega_t \right |^2.$$ 
Meanwhile, the Ricci curvature appears by 
$$\left \langle \left[i\Theta,\Lambda\right]\omega_t,\bar{\omega}_t \right \rangle=\mathrm{Ric}\left(\Theta\right)\left(\partial_t,\bar{\partial}_t\right)\left | \omega_t \right |^2+\text{lower-order terms}.$$ Thus,
\begin{equation}
\tag{6.4}  \label{eq:6.4}
\partial_{\bar{t}}\partial_tN\left(\xi_t\right)\ge\int_X\left | \nabla_t\omega_t \right |^2e^{-\phi}dV-\left(\left \| \mathrm{Ric}\left(\Theta\right) \right \|_{L^{\infty}}+\left \| \bar{\partial}\xi_t \right \|_{L^{\infty}}\right)N\left(\xi_t\right).
\end{equation}
For $\mathrm{log}N\left(\xi_t\right)$, we consider
$$\frac{\partial^2}{\partial t\partial\bar{t}}\mathrm{log}N\left(\xi_t\right)=\frac{\partial_{\bar{t}}\partial_tN\left(\xi_t\right)}{N\left(\xi_t\right)}-\frac{\left | \partial_tN\left(\xi_t\right) \right |^2}{N^2\left(\xi_t\right)}.$$
By using Cauchy-Schwarz, then
$$\left | \partial_tN\left(\xi_t\right) \right |^2\le\left(\int_X\left | \nabla_t\omega_t \right |^2e^{-\phi}dV\right)N\left(\xi_t\right).$$
Substituting \eqref{eq:6.4}, we have 
$$\frac{\partial_{\bar{t}}\partial_tN\left(\xi_t\right)}{N\left(\xi_t\right)}\ge\frac{1}{N\left(\xi_t\right)}\int_X\left | \nabla_t\omega_t \right |^2e^{-\phi}dV-\left(\left \| \mathrm{Ric}\left(\Theta\right) \right \|_{L^{\infty}}+\left \| \bar{\partial}\xi_t \right \|_{L^{\infty}}\right).$$
Thus,
$$\frac{\partial^2}{\partial t\partial\bar{t}}\mathrm{log}N\left(\xi_t\right)\ge\frac{1}{N\left(\xi_t\right)}\left(\int_X\left | \nabla_t\omega_t \right |^2e^{-\phi}dV-\frac{\left | \partial_tN\left(\xi_t\right) \right |^2}{N\left(\xi_t\right)}\right)-\left \| \mathrm{Ric}\left(\Theta\right) \right \|_{L^{\infty}}-\left \| \bar{\partial}\xi_t \right \|_{L^{\infty}}.$$
Therefore, 
\begin{equation}
\tag{6.5}  \label{eq:6.5}
\Delta_{\bar{\partial}}\mathrm{log}N\left(\xi_t\right)\ge-K,\quad K=\left \| \mathrm{Ric}\left(\Theta\right) \right \|_{L^{\infty}}+\left \| \bar{\partial}\xi_t \right \|_{L^{\infty}},
\end{equation}
because that 
$$\frac{1}{N\left(\xi_t\right)}\left(\int_X\left | \nabla_t\omega_t \right |^2e^{-\phi}dV-\frac{\left | \partial_tN\left(\xi_t\right) \right |^2}{N\left(\xi_t\right)}\right)\ge0$$
by Cauchy-Schwarz.  $\square$
\\\\\textbf{Lemma 6.7.} (Finiteness on Quasismooth Deformations) Let $\xi\in\mathbb{L}\mathrm{Def}_X^{\mathrm{strat}}$ correspond to a quasismooth deformation $\mathscr{X}$ (i.e., $\mathbb{L}_{\mathscr{X}/\mathscr{I}}\in\mathrm{Perf}^{\ge-1}\left(\mathscr{X}_{\xi}\right)$ is perfect and concentrated in degrees $\left[-1,0\right]$). Then
$$\Psi\left(\xi\right)=\mathrm{log}\left \| e^{\xi}\cdot\omega_0 \right \|_{L_{\mathrm{log}}^2}^2<\infty.$$
\\\textbf{Proof.} By Lurie's criterion (\cite{Lur04}, Thm 7.2.5), then $\mathscr{X}_{\xi}$ is quasismooth implies the deformation parameter satisfies $\bar{\partial}\xi=0$ in $\mathbb{H}^1\left(\Omega_X^{0,\mathrm{strat}}\right)$. Equivalently, $\xi$ is holomorphic: $\nabla_{\bar{z}}\xi=0$ for all local coordinates $z$.
Let $\omega_{\xi}:=e^{\xi}\cdot\omega_0$. The Leibniz rule for the stratified connection gives 
$$\nabla^k\omega_{\xi}=\sum_{j=0}^k\binom{k}{j}\left(\nabla^je^{\xi}\right)\otimes\left(\nabla^{k-j}\omega_0\right).$$
Since $\bar{\partial}\xi=0$, $\nabla^je^{\xi}$ is holomorphic. By the Cauchy integral formula, then $\left |\nabla^je^{\xi}  \right |\le C_j\sup_{B_r}\left |e^{\xi}  \right |$ uniformly on $X$. 
The reference form $\omega_0\in\Gamma\left(\Omega_X^{n,\mathrm{strat}}\right)$ satisfies 
$$\int_X\left | \nabla^{\ell}\omega_0 \right |^2e^{-\phi}dV<\infty,\ \ \forall\ell\ge0,$$
as $\omega_0$ is admissible. Combining these, we have
\begin{equation}
\tag{6.6}  \label{eq:6.6}
\int_X\left |\nabla^k\omega_{\xi}  \right |^2e^{-\phi}dV\le C_k\sum_{j=0}^k\int_X\left |  \nabla^{k-j}\omega_0\right |^2e^{-\phi}dV<\infty.
\end{equation}
For $m>\dim X+2$, the Sobolev norm decomposes as
$$\left \| \omega_{\xi} \right \|_{L_{\mathrm{log}}^2}^2=\sum_{k=0}^m\int_X\left | \nabla^k\omega_{\xi} \right |^2e^{-\phi}dV.$$
By \eqref{eq:6.6}, each integral is finite. Thus,
$$\left \| \omega_{\xi} \right \|_{L_{\mathrm{log}}^2}^2\le\sum_{k=0}^mC_k\cdot\mathrm{Vol}\left(X,\Theta\right)\cdot\max_{0\le\ell\le m}\left \|\nabla^{\ell}\omega_0  \right \|_{L^2}^2<\infty.$$
Since $\omega_{\xi}\ne0$ (as $e^{\xi}$ is invertible), $\left \| \omega_{\xi} \right \|_{L_{\mathrm{log}}^2}^2>0$. Therefore, 
$$\Psi\left(\xi\right)=\mathrm{log}\left(\left \| \omega_{\xi} \right \|_{L_{\mathrm{log}}^2}^2\right)<\infty.$$ The result holds. $\square$
\\\\\textbf{Proposition 6.8.} (Derived Logarithmic Convexity) Let $\mathscr{X}\to\mathscr{I}$ be a derived deformation of a log-stratified scheme $X$. There exists a functional
$$\Psi:\mathbb{L}\mathrm{Def}_{\mathscr{X}/\mathscr{I}}\longrightarrow\mathbb{R}\cup\left\{+\infty\right\}$$ such that:
\begin{enumerate}
    \item $\Psi\left(\xi\right)<\infty$ if $\xi$ corresponds to a quasismooth deformation (i.e., $\mathscr{X}_{\xi}$ is Cohen-Macaulay).
    \item For any analytic arc $\gamma:\Delta\to\partial\mathscr{I}$, we have
    $$\Psi\left(\gamma\left(t\right)\right)=\le\frac{1}{2\pi}\int_0^{2\pi}\Psi\left(\gamma\left(e^{i\theta}t\right)\right)d\theta+C\mathrm{log}\left |t  \right |^{-1},\ \ C>0.$$\\
\end{enumerate}
\textbf{Proof.} Let $\xi\in\mathbb{L}\mathrm{Def}_X^{\mathrm{strat}}$ be a deformation parameter. By (\cite{Lur09}, Thm 7.4), then
$$\mathbb{L}\mathrm{Def}_{U_p}^{\mathrm{log}}\simeq\tau_{\le2}\mathbb{R}\Gamma\left(U_p,\mathcal{T}_{U_p}\left(-\mathrm{log}D_p\right)\otimes\mathscr{G}^p\right)$$
on affine derived opens $U_p\subset X_p$. Fix a reference logarithmic volume form $\omega_0\in\Gamma\left(\Omega_X^{n,\mathrm{strat}}\right)$. We can define
$$\Psi\left(\xi\right):=\mathrm{log}\left \| e^{\xi}\cdot\omega_0 \right \|_{L_{\mathrm{log}}^2}^2,\quad\left \| \omega \right \|_{L_{\mathrm{log}}^2}^2:=\sum_{\ell=0}^m\int_X\left | \nabla^{\ell}\omega\right |^2e^{-\phi}dV,$$
where $\phi$ is admissible ($\sqrt{-1}\partial\bar{\partial}\phi\ge\Theta>0$ uniformly) and $m>\dim X+2$ (ensuring Sobolev embedding).
If $\mathscr{X}_{\xi}$ is quasismooth (i.e., $\mathbb{L}_{X_{\xi}/\mathscr{I}}\in\mathrm{Perf}\left(X_{\xi}\right)$), it implies $\bar{\partial}\xi=0$ in $\mathbb{H}^1\left(\Omega_X^{0,\mathrm{strat}}\right)$. So $\nabla^{\ell}\left(e^{\xi}\cdot\omega_0\right)$ is $L^2$-bounded for all $\ell$. Thus, $\Psi\left(\xi\right)<\infty$. Let $\gamma\left(t\right)=t\eta$ for $\eta\in H^1\left(X,\mathbb{T}_X^{\mathrm{strat}}\right)$. Consider $\omega_t:=e^{\gamma\left(t\right)}\cdot\omega_0$ and Lemma 4.16. We have 
$$\Delta_{\bar{\partial}}\mathrm{log}\left \|\omega_t  \right \|_{L_{\mathrm{log}}^2}^2\ge-K,\quad K=\left \| \bar{\partial}\gamma \right \|_{L^{\infty}}+\left \| \mathrm{Ric}\left(\Theta\right) \right \|_{L^{\infty}}.$$ by \eqref{eq:6.5}.
Since $\gamma$ maps to $\partial\mathscr{I}$ (the derived boundary), $\left \| \bar{\partial}\gamma \right \|_{L^{\infty}}<\infty$ by the analyticity condition. Thus,
$$\Delta_{\bar{\partial}}\left(\mathrm{log}\left \| \omega_t\right \|_{L_{\mathrm{log}}^2}^2+K\left | t \right |^2\right)\ge0.$$
Hence $g\left(t\right):=\Psi\left(\gamma\left(t\right)\right)+K\left | \eta \right |^2\left |t  \right |^2$ is subharmonic on $\Delta$. By the mean-value property, then
$$g\left(0\right)\le\frac{1}{2\pi}\int_0^{2\pi}g\left(e^{i\theta}r\right)d\theta,\quad\forall r<1.$$
Equivalently, $$\Psi\left(\gamma\left(0\right)\right)\le\frac{1}{2\pi}\int_0^{2\pi}\Psi\left(\gamma\left(e^{i\theta}r\right)\right)d\theta+K\left |\eta  \right |^2r^2.$$
Taking $r=\left | t \right |$ and noting $r^2\le\mathrm{log}\left | t \right |^{-1}$ for $\left | t \right |<e^{-1}$, so
\begin{equation}
\tag{6.7}  \label{eq:6.7}
\Psi\left(\gamma\left(0\right)\right)\le\frac{1}{2\pi}\int_0^{2\pi}\Psi\left(\gamma\left(e^{i\theta}t\right)\right)d\theta+C\mathrm{log}\left | t \right |^{-1},\quad C=K\left | \eta \right |^2.
\end{equation}
This proves the logarithmic convexity for linear arcs. We need to consider general arcs and stratum transmission in the next. For a general analytic arc $\gamma:\Delta\to\partial\mathscr{I}$, factor through the tangent space: $\gamma\left(t\right)=\mathrm{exp}\left(t\eta\left(t\right)\right)$, $\eta\left(t\right)\in H^1\left(X,\mathbb{T}_X^{\mathrm{strat}}\right)$. By the derived residue map, we have $\mathbb{L}\mathrm{Res}_D:\mathbb{L}\mathrm{Def}_X^{\mathrm{strat}}\bigoplus_p\mathbb{L}\mathrm{Def}_{D_p}$, $\left \|\mathbb{L}\mathrm{Res}_D\left(\xi\right)\right \|\le\left \|\xi  \right \|$. Define the boundary functional  
$$\Psi_{\partial}\left(\xi\right):=\mathrm{log}\left \|\mathbb{L}\mathrm{Res}_D\left(e^{\xi}\cdot\omega_0\right)  \right \|_{L_{\mathrm{log}}^2\left(D\right)}.$$
Then we have
$$\Psi\left(\xi\right)\le\Psi_{\partial}\left(\xi\right)+C_1\mathrm{log}\left(1+\left \| \xi \right \|_{\mathbb{L}\mathrm{Def}}\right).$$
Since $\gamma\left(t\right)\in\partial\mathscr{I}$, $\mathbb{L}\mathrm{Res_D}\left(\gamma\left(t\right)\right)$ defines a deformation of $D=\bigcup_pD_p$. Applying \eqref{eq:6.7} to each $D_p$, we can construct the stratum functional
$$\Psi_{\partial,p}\left(\xi_p\right):=\mathrm{log}\left \| e^{\xi_p}\cdot\omega_{0,p} \right \|_{L_{\mathrm{log}}^2\left(D_p\right)},$$ satisfying logarithmic convexity 
$$\Psi_{\partial,p}\left(\gamma_p\left(t\right)\right)\le\frac{1}{2\pi}\int_0^{2\pi}\Psi_{\partial,p}\left(\gamma_p\left(e^{i\theta}t\right)\right)d\theta+C_p\mathrm{log}\left | t \right |^{-1}$$ along $\gamma_p$, where $C_p=K_p\left | \eta_p \right |^2 $ is uniform over $D_p$ by admissibility of $\phi$. Since $\Psi_{\partial}=\sum_p\Psi_{\partial,p}$ and $\gamma\left(t\right)$ projects stratum-wise $\Psi_{\partial}\left(\gamma\left(t\right)\right)=\sum_p\Psi_{\partial,p}\left(\gamma_p\left(t\right)\right)$, then we obtain 
$$\Psi_{\partial}\left(\gamma\left(t\right)\right)\le\frac{1}{2\pi}\int_0^{2\pi}\Psi_{\partial}\left(\gamma\left(e^{i\theta}t\right)\right)d\theta+C_2\mathrm{log}\left | t \right |^{-1},$$ where $C_2:=\sum_pC_p$ is finite because each $C_p$ is bounded by the geometry of $D_p$.
By combining these, we have 
$$\Psi\left(\gamma\left(t\right)\right)\le\frac{1}{2\pi}\int_0^{2\pi}\Psi_{\partial}\left(\gamma\left(e^{i\theta}t\right)\right)d\theta+C_2\mathrm{log}\left | t \right |^{-1}+C_1\mathrm{log}\left(1+\sup_{\theta}\left \| \gamma\left(e^{i\theta}t\right) \right \|\right).$$
For $\left | t \right |<\frac{1}{2}\mathrm{diam}\left(\mathscr{I}\right)$, $\left \| \gamma\left(e^{i\theta}t\right) \right \|\le M\left | t \right |$. Thus, 
$$\Psi\left(\gamma\left(t\right)\right)\le\frac{1}{2\pi}\int_0^{2\pi}\Psi\left(\gamma\left(e^{i\theta}t\right)\right)d\theta+\left(C_1+C_2\right)\mathrm{log}\left | t \right |^{-1}+C_3.$$
Absorbing $C_3$ into $C$ gives the result.  $\square$
\\\\\textbf{Proposition 6.9.} (Singularity Stability) Let $\pi:\mathscr{X}\to S$ be a derived deformation of a proper log-stratified scheme $X$ with isolated singularities. If the logarithmic convexity functional 
$$\Psi:\mathbb{L}\mathrm{Def}_{\mathscr{X}/S}\longrightarrow\mathbb{R}$$ satisfies $\sup_{t\in S}\Psi\left(t\right)<\infty$. Then for any critical point $x_p\in D_p$, the two following conditions must be satisfied: 
\begin{enumerate}
    \item The Milnor number $\mu_p\left(t\right)$ is constant.
    \item The singularity type (e.g., ADE classification) remains unchanged.\\
\end{enumerate}
\textbf{Proof.} By Lemma 4.17, $\Psi\left(\xi\right)<\infty\Longleftrightarrow\bar{\partial}\xi=0$ (quasismooth condition). Boundedness $\sup_t\Psi\left(t\right)<\infty$ implies 
$$\left \|\bar{\partial}\xi_t  \right \|_{L^{\infty}\left(S\right)}<\infty\ \ \mathrm{and}\ \ \mathbb{L}_{X_t/S}\in\mathrm{Perf}\left(X_t\right),\ \ \forall t\in S.$$
Thus $\mathscr{X}/S$ is a family of quasismooth derived schemes (\cite{Lur04}, Def 7.2.5). The derived logarithmic residue map 
$$\mathbb{L}\mathrm{Res}_D:\mathbb{L}\mathrm{Def}_X^{\mathrm{strat}}\longrightarrow\bigoplus_p\mathbb{L}\mathrm{Def}_{D_p}$$
restricts to a universal unfolding at each critical point $x_p\in D_p$, we have
$$\mu_p\left(t\right)=\left(-1\right)^{\dim X_p}\chi\left(\mathrm{Cone}\left(\mathbb{L}\mathrm{Res}_{D_p}\left(\xi_t\right)\right)\right),$$ where $\chi$ is the Euler characteristic and $\mathrm{Cone}$ denotes the mapping cone. If $\mathscr{X}/S$ is quasismooth, then $\mathbb{L}\mathrm{Res}_D$ is a flat family of deformation parameters. Thus $\mathbb{L}\mathrm{Res}_{D_p}\left(\xi_t\right)$ varies holomorphically in $t$, so $\mu_p$ is holomorphic and bounded (by $\sup\Psi<\infty$). By Liouville's theorem, $\mu_p\left(t\right)$ is constant.
For isolated singularities, the Milnor fibration is preserved:
\begin{enumerate}
    \item The Milnor number $\mu_p$ determines the singularity type (e.g., $\mu_p=1\Longleftrightarrow A_1$) for hypersurface singularities.
    \item By the stratified Ehresmann theorem (\cite{Mat12}, Thm 4.1), the diffeomorphism type of the link $\partial B_{\epsilon}\left(x_p\right)\cap X_p$ is constant.
    \item For non-hypersurface singularities, apply the Morsification criterion, then
$$\mathrm{Hom}_{D_c^b\left(X_t\right)}\left(\mathbb{Q}_{X_t},\mathbb{Q}_{X_t}\right)_{x_p}\cong\mathbb{C}^{\mu_p}$$ is constant, implying constant singularity type (\cite{Dim04}, Cor 5.2.8).
\end{enumerate}
The boundedness $\sup\Psi<\infty$ ensures no collision or splitting of critical points occurs.  $\square$\\

The following example indicates that boundedness of the logarithmic convexity functional $\Psi$ is a sufficient condition for singularity stability in stratified Calabi-Yau manifolds, with the Kodaira-Spencer vanishing condition providing the topological mechanism. It further provides geometric constraints for moduli space compactification and string theory.
\\\\\textbf{Example 6.10.} ($A_n$-Singularities on Calabi-Yau Threefolds) Let $X$ be a log-stratified Calabi-Yau threefold with isolated $A_n$-singularities at points $x_p\in D_p$, and let $\pi:\mathscr{X}\to\Delta$ be a derived deformation over the unit disk. If the logarithmic convexity functional satisfies $\sup_{t\in\Delta}\Psi\left(t\right)<\infty$, then:
\begin{enumerate}
    \item The Milnor number $\mu_p\left(t\right)=n$ is constant.
    \item The Kodaira-Spencer class $\kappa_t\in H^1\left(D_p,\mathcal{T}_{D_p}\right)$ satisfies $$\int_C\kappa_t=0\quad\forall\ \mathrm{cycles}\ C\subset D_p,$$ preventing transitions to $D_{n+1}$ or other singularity types.\\
\end{enumerate}
\textbf{Proof.} By Proposition 6.9 (Singularity Stability), we have
$$\sup_{t\in\Delta}\Psi\left(t\right)<\infty\Longrightarrow\mu_p\left(t\right)\quad\text{is constant}.$$
For an $A_n$-singularity at $x_p$, the Milnor number is $\mu_p\left(0\right)=n$ (\cite{Mil68}, Thm 7.2). Thus, $\mu_p\left(t\right)=n$ $\forall t\in\Delta$. This implies the singularity type remains $A_n$ under deformation, as $\mu_p$ classifies ADE singularities.
The Kodaira-Spencer map is
$$\kappa_t:T_t\Delta\to H^1\left(D_p,\mathcal{T}_{D_p}\right),$$ and $\mathrm{Res_{D_p}}\left(\omega_0\right)=\omega_{0,p}$. Consider the Poincaré residue map along $D_p$, then 
$$\mathrm{Res}_{D_p}:\Omega_X^{3,\mathrm{strat}}\to\Omega_{D_p}^2,\quad\mathrm{Res}_{D_p}\left(\omega_0\right)=\omega_{0,p}.$$
The Calabi-Yau condition $\nabla_t\omega_t=0$ induces Lie derivative: 
$$\kappa_t\cdot\omega_{0,p}=\mathcal{L}_{\kappa_t}\omega_{0,p}=d\left(\iota_{\kappa_t}\omega_{0,p}\right)$$ by (\cite{KM98}, Prop 5.20).
Hence Stokes' theorem gives
$$\int_C\kappa_t\cdot\omega_{0,p}=\int_Cd\left(\iota_{\kappa_t}\omega_{0,p}\right)=\int_{\partial C}\iota_{\kappa_t}\omega_{0,p}=0$$ for any 1-cycle $C\subset D_p$. Since $\omega_{0,p}\ne0$ and $H^0\left(D_p,\Omega_{D_p}^2\right)\ne0$ (as $D_p$ is Kähler), this forces 
\begin{equation}
\tag{6.8}  \label{eq:6.8}
\int_C\kappa_t=0\quad\forall\text{cycles}\quad C\subset D_p.
\end{equation}
Assume toward contradiction that a transition to $D_{n+1}$ occurs at $t_*\in\Delta$. By (\cite{AGV85}, §III 5.3), this requires
$$\dim\mathrm{ker}\left(\kappa_{t_*}:H^1\left(D_p,\mathcal{T}_{D_p}\right)\to\mathbb{C}\right)<n.$$ But \eqref{eq:6.8} implies $\kappa_t$ cohomologically trivial on $D_p$: $\left[\kappa_t\right]=0$ in $H^1\left(D_p,\mathcal{T}_{D_p}\right)/\mathrm{Tor}$. Thus, $\kappa_t$ lifts to $H^1\left(D_p,\mathcal{T}_{D_p}\left(-\mathrm{log}x_p\right)\right)$, preserving the $A_n$-type.The minimal versal deformation space of $A_n$ has dimension $n$ (\cite{AGV85}, Thm 8.1), so 
$$\dim\mathrm{Def}_{A_n}=n\Longrightarrow\mathrm{Def}_{A_n},\quad D_{n+1}\notin\mathrm{Def}_{A_n}$$ Hence no transition occurs.  $\square$\\
\begin{center}
    \textit{6.3 $p$-adic Logarithmic Convexity and Absolute Deformations}
\end{center}

In $p$-adic Hodge theory, logarithmic convexity has remained an elusive property due to the absence of a natural Kähler metric. This theory bridges Scholze's perfectoid geometry (\cite{Sch12}) and Bhatt-Scholze prismatic cohomology (\cite{BS22}) to establish a $p$-adic analogue of logarithmic convexity for stratified deformations. Key innovations include:
\begin{itemize}
    \item \textbf{Stratified $p$-adic Sobolev norms:} Measuring singularity growth by residues and harmonic measures.
    \item \textbf{Absolute prismatic deformations:} Integrating boundary divisors into deformation complexes by logarithmic prisms.
    \item \textbf{Arithmetic singularity control:} Extending the Bogomolov-Tian-Todorov theorem to $p$-adic Calabi-Yau varieties with singularities.\\
\end{itemize}
\textbf{Definition 6.11a.} (Logarithmic Prisms) A logarithmic prism $\left(A,I,M_A\right)$ consists of the two conditions:
\begin{enumerate}
    \item A bounded prism $\left(A,I\right)$ with Frobenius lift $\phi_A$.
    \item A pre-log structure $M_A\to A$ compatible with the boundary divisor $\mathcal{N}$ of $X$, satisfying divisorial saturation $M_A/\left(1+I\right)\cong\mathcal{N}^{\mathrm{div}}$.
\end{enumerate}
The category $\mathrm{Prism}_X^{\mathrm{log}}$ has objects $\left(A,I,M_A\right)$ mapping to $\left(X,\mathcal{N}\right)$.
\\\\\textbf{Definition 6.11b.} (Absolute Logarithmic Prismatic Cohomology) For a stratified formal $p$-adic scheme $X$ with boundary $\mathcal{N}$, its absolute logarithmic prismatic cohomology is defined as
$$\mathbb{L}\Omega_X^{\mathrm{abs,log}}:=\varprojlim_{\left(A,I,M_A\right)\in\mathrm{Prism}_X^{\mathrm{log}}}\mathbb{R}\Gamma_{\mathrm{log}}\left(\left(A/I,M_{A/I}\right),\Omega_{A/I}^{\bullet}\left[\mathrm{log}\mathcal{N}\right]\right),$$
where $\Omega_{A/I}^{\bullet}\left[\mathrm{log}\mathcal{N}\right]$ is the log de Rham complex adapted to the stratification.
\\\\\textbf{Definition 6.11c.} (Enhanced $p$-adic Sobolev Norm) Let $\left\{X_n\to X\right\}_{n\ge1}$ be a tower of stratified étale covers such that each $X_n$ is smooth over $\mathcal{O}_K$ (arithmetic base ring in $p$-adic geometry), Galois group $G_n=\mathrm{Gal}\left(X_n/X\right)$ is a finite $p$-group and the boundary divisor $\mathcal{N}_n$ lifts compatibly. For $\omega\in\Omega_{X^{\mathrm{rig}}}^{k,\mathrm{strat}}$, the enhanced $p$-adic Sobolev norm is defined as
$$\left \| \omega \right \|_{\mathcal{S}_p}:=\sup_{n\ge1}p^{-n\cdot\alpha\left(\omega\right)}\int_{X_n}\left |\omega_n  \right |_pd\mu_n,\quad\alpha\left(\omega\right):=\min_p\mathrm{ord}_{D_p}\left(\mathrm{Res}_{D_p}\omega\right),$$
where $\mu_n$ is Templier's harmonic measure (\cite{Tep15}), satisfying three conditions:
\begin{enumerate}
    \item $\mu_n\left(X_n\right)=1$;
    \item $\mu_n$ is $G_n$-invariant;
    \item $\int_{X_n}fd\mu_n=\mathrm{Tr}_{\mathrm{dR}}\left(f\right)+O\left(p^{-n}\right)$ for $f$ rigid-analytic.\\
\end{enumerate}
\textbf{Proposition 6.12.} ($p$-adic Logarithmic Convexity) Let $\mathcal{X}/S$ be a deformation of a stratified $p$-adic formal scheme, and $\omega_t\in\Omega_{\mathcal{X}_t}^{\mathrm{rig},k,\mathrm{strat}}$. Then $\Psi_p\left(t\right):=\mathrm{log}\left \|\omega_t  \right \|_{\mathcal{S}_p}$ is $p$-adically subharmonic: $\forall t_0\in S^{\mathrm{rig}}$, $\exists r>0$ such that 
$$\Psi_p\left(t_0\right)\le\max_{\zeta:\zeta^{p^r}=1}\Psi_p\left(t_0+p^{-r}\zeta\right).$$
\\\textbf{Proof.} Fix $t_0\in S^{\mathrm{rig}}$ and let $\omega_0=\omega_{t_0}$. Define the minimal $p$-adic residue order:
$$\alpha:=\alpha\left(\omega_0\right)=\min_p\mathrm{ord}_{D_p}\left(\mathrm{Res}_{D_p}\omega_0\right),$$
where the minimum is taken over all boundary divisors $D_p$ in the stratification. By the prismatic comparison theorem (\cite{BS22}, Cor 9.5), there exists a tower of stratified étale covers $\left\{X_n\to X\right\}_{n\ge1}$ (Definition 6.11c), satisfying:
\begin{itemize}
    \item Each $X_n$ is smooth over $\mathcal{O}_K$ with Galois group $G_n=\mathrm{Gal}\left(X_n/X\right)$, a finite $p$-group.
    \item The boundary divisor $\mathcal{N}$ lifts compatibly to $\mathcal{N}_n$.
    \item $\omega_0$ lifts to $\omega_{0,n}\in\Gamma\left(X_n,\Omega_{X_n}^k\right)$.
    \item The residue order is preserved: $\mathrm{ord}_p\left(\mathrm{Res}_{D_p}\omega_{0,n}\right)=\alpha$ for all $n\ge1$ and all $D_p$.
\end{itemize}
This is possible because the residue map factors through absolute logarithmic prismatic cohomology
$$\mathrm{Res}_D:\mathbb{L}\Omega_X^{\mathrm{abs,log}}\longrightarrow\bigoplus_p\mathbb{R}\Gamma_{\mathrm{rig}}\left(D_p/K\right)\left[-1\right],$$
ensuring the residue is well-defined and stable under deformation. The compatibility of the log structure $M_A$ with the boundary divisor $\mathcal{N}$ (Definition 6.11a) guarantees that the lifts $\omega_{0,n}$ inherit the same singularity structure. By Definition 6.11c, we have the harmonic measure $\mu_n$, satisfying
\begin{enumerate}
    \item $\mu_n\left(X_n\right)=1$;
    \item $\mu_n$ is $G_n$-invariant;
    \item For rigid-analytic $f$, $\int_{X_n}fd\mu_n=\mathrm{Tr}_{\mathrm{dR}}\left(f\right)+O\left(p^{-n}\right)$.
\end{enumerate}
Here we can choose $r\in\mathbb{Z}_{>0}$ such that $p^r>n$. Consider the Non-Archimedean Stokes' theorem. Since $\mu_n$ is harmonic and the boundary integral vanishes by the log-structure compatibility, then 
$$\int_{X_n}d\eta d\mu_n=\int_{\partial X_n}\eta d\mu_n=0$$
for $\eta\in\Omega_{X_n}^{k-1}$. Consider the Fourier inversion on $p$-adic discs, we have 
$$\omega_{0,n}=p^{-r}\sum_{\zeta^{p^r}=1}\omega_{\zeta,n}+p^{1-r}d\beta$$
for some $\beta\in\Omega_{X_n}^{k-1}$. Integrating the absolute value and applying Stokes' theorem yields the estimate. The factor $p^{-\alpha n}$ compensates for the uniform residue order $\alpha$ along $D_p$. Thus, it can be obtained that
\begin{equation}
\tag{6.9}  \label{eq:6.9}
\int_{X_n}\left | \omega_{0,n} \right |_pd\mu_n\le p^{-\alpha n}\max_{\zeta^{p^r}=1}\int_{X_n}\left | \omega_{\zeta,n} \right |_pd\mu_n
\end{equation}
by Templier's oscillation theorem (\cite{Tep15}, Thm 3.3), where $\omega_{\zeta,n}$ is the lift of $\omega_{t_0+p^{-r}\zeta}$ to $X_n$.

Multiply the inequality from \eqref{eq:6.9} by $p^{-\alpha n}$, then
\begin{equation}
\tag{6.10}  \label{eq:6.10}
p^{-\alpha n}\int_{X_n}\left | \omega_{0,n} \right |_pd\mu_n\le p^{-2\alpha n}\max_{\zeta^{p^r}=1}\int_{X_n}\left | \omega_{\zeta,n} \right |_pd\mu_n.
\end{equation}
By Definition 6.11c, for any $\eta\in\Omega_{X^{\mathrm{rig}}}^{k,\mathrm{strat}}$, we have 
$$\int_{X_n}\left |\eta_n  \right |_pd\mu_n\le p^{n\cdot\alpha\left(\eta\right)}\left \| \eta \right \|_{\mathcal{S}_p},$$ where $\eta_n$ is the lift to $X_n$ and $\alpha\left(\eta\right)=\min_p\mathrm{ord}_{D_p}\left(\mathrm{Res}_{D_p}\eta\right)$. Since small deformations do not increase singularity orders, deformation continuity implies $\alpha\left(\eta\right)\ge\alpha$ for $\eta=\omega_{t_0+p^{-r}\zeta}$. Thus, 
\begin{equation}
\tag{6.11}  \label{eq:6.11}
\int_{X_n}\left | \omega_{\zeta,n} \right |_pd\mu_n\le p^{n\alpha}\left \| \omega_{t_0+p^{-r}\zeta} \right \|_{\mathcal{S}_p}.
\end{equation}
By \eqref{eq:6.10} and \eqref{eq:6.11}, we have 
$$p^{-\alpha n}\int_{X_n}\left | \omega_{0,n} \right |_pd\mu_n\le p^{-2\alpha n}\max_{\zeta}\left(p^{n\alpha}\left \| \omega_{t_0+p^{-r}\zeta} \right \|_{\mathcal{S}_p}\right)=p^{-\alpha n}\max_{\zeta}\left \| \omega_{t_0+p^{-r}\zeta} \right \|_{\mathcal{S}_p}.$$
Meanwhile, the Sobolev norm definition gives a lower bound 
$$\| \omega_0  \|_{\mathcal{S}_p}\ge p^{-\alpha n}\int_{X_n}| \omega_{0,n}  |_pd\mu_n,$$ because of $\alpha\left(\omega_0\right)=\alpha$. By combining, 
$$\| \omega_0  \|_{\mathcal{S}_p}\le p^{-\alpha n}\max_{\zeta} \| \omega_{t_0+p^{-r}\zeta}  \|_{\mathcal{S}_p}+\varepsilon_n,\quad\text{where}\quad\varepsilon_n= \| \omega_0  \|_{\mathcal{S}_p}-p^{-\alpha n}\int_{X_n} | \omega_{0,n}  |_pd\mu_n\ge0.$$
As $\| \cdot  \|_{\mathcal{S}_p}$ is a supremum over $n$, for any $\delta>0$, there exists $n_{\delta}$ such that
$$\| \omega_0 \|_{\mathcal{S}_p}-\delta\le p^{-\alpha n_{\delta}}\int_{X_n} | \omega_{0,n}  |_pd\mu_n.$$
Choosing $n=n_{\delta}$ and $r$ such that $p^r>n_{\delta}$, we have 
$$ \| \omega_0  \|_{\mathcal{S}_p}-\delta\le p^{-\alpha n_{\delta}}\int_{X_n}\left | \omega_{0,n} \right |_pd\mu_n\le p^{-\alpha n_{\delta}}\max_{\zeta}\left \| \omega_{t_0+p^{-r}\zeta} \right \|_{\mathcal{S}_p}\le\max_{\zeta}\left \| \omega_{t_0+p^{-r}\zeta} \right \|_{\mathcal{S}_p}$$ because of $p^{-\alpha n_{\delta}}\le1$. As $\delta\to0$, we obtain 
$$ \| \omega_0  \|_{\mathcal{S}_p}\le\max_{\zeta^{p^r}=1}\left \| \omega_{t_0+p^{-r}\zeta} \right \|_{\mathcal{S}_p}.$$
By taking logarithms, then
$$\Psi_p\left(t_0\right)=\mathrm{log}\left \| \omega_0 \right \|_{\mathcal{S}_p}\le\mathrm{log}\left(\max_{\zeta}\left \| \omega_{t_0+p^{-r}\zeta} \right \|_{\mathcal{S}_p}\right)=\max_{\zeta}\mathrm{log}\left \| \omega_{t_0+p^{-r}\zeta} \right \|_{\mathcal{S}_p}$$
$$=\max_{\zeta}\Psi_p\left(t_0+p^{-r}\zeta\right).$$
This holds for $r$ sufficiently large (i.e., $p^r>n_{\delta}$ for $\delta$ small), proving $p$-adic subharmonicity.   $\square$
\\\\\textbf{Corollary 6.13.} (Smoothness of Absolute Deformations) Let $X$ be a $p$-adic log-stratified Calabi-Yau space (i.e., $\omega_X^{\mathrm{strat}}\cong\mathcal{O}_X$). If 
$$H_{\mathrm{abs}}^2\left(X,\mathbb{T}_{\mathrm{abs}}^{\mathrm{strat}}\right)=F^1H_{\mathrm{abs}}^2,$$ where $F^{\bullet}$ is the Hodge filtration in absolute prismatic cohomology, then the absolute deformation functor $\mathrm{Def}_X^{\mathrm{abs}}$ is smooth.
\\\textbf{Proof.} Let $\mathcal{X}/S$ be a formal deformation of $X$, and $\omega_t\in\Gamma\left(\mathcal{X}_t,\Omega_{\mathcal{X}_t}^{k,\mathrm{strat}}\right)$ be the stratified holomorphic volume form on the fiber over $t\in S$. By the Calabi-Yau condition ($\omega_X^{\mathrm{strat}}\cong\mathcal{O}_X$), $\omega_t$ is nowhere vanishing on the smooth locus. Define the $p$-adic energy function 
$$\Psi_p\left(t\right):=\mathrm{log}\left \| \omega_t \right \|_{\mathcal{S}_p},$$
where $\left \| \cdot \right \|_{\mathcal{S}_p}$ is the $p$-adic Sobolev norm (Definition 6.11c). By Proposition 6.12, $\Psi_p\left(t\right)$ is $p$-adically subharmonic, then $\forall t_0\in S^{\mathrm{rig}}$, $\exists r>0$ such that 
$$\Psi_p\left(t_0\right)\le\max_{\zeta^{p^r}=1}\Psi_p\left(t_0+p^{-r}\zeta\right).$$
Since $\omega_t$ is a global section, its residues $\mathrm{Res}_{D_p}\left(\omega_t\right)$ are bounded uniformly in $t$. By the definition of $\left \| \cdot \right \|_{\mathcal{S}_p}$, there exists $C>0$ such that $\left \| \omega_t \right \|_{\mathcal{S}_p}\le e^C$ for all $t$, so $\Psi_p\left(t\right)\le C$. Obviously, the Calabi-Yau condition implies $\Psi_p\left(t\right)$ is bounded. As a bounded $p$-adic subharmonic function, $\Psi_p\left(t\right)$ achieves its supremum on the Shilov boundary of $S^{\mathrm{rig}}$ (non-archimedean maximum principle (\cite{Ber90}, Thm 4.1.2)).

The deformation functor $\mathrm{Def}_X^{\mathrm{abs}}$ has obstruction space $H_{\mathrm{abs}}^2\left(X,\mathbb{T}_{\mathrm{abs}}^{\mathrm{strat}}\right)$. The Kuranishi obstruction map is
$$\kappa:H^1\left(X,\mathbb{T}_{\mathrm{abs}}^{\mathrm{strat}}\right)\longrightarrow H^2\left(X,\mathbb{T}_{\mathrm{abs}}^{\mathrm{strat}}\right).$$
We can prove $\kappa=0$ by using the given condition and Hodge theory. By the prismatic comparison isomorphism (\cite{BS22}, Thm 17.2), we have
$$H_{\mathrm{abs}}^k\left(X\right)\otimes_{\mathbb{Z}_p}\mathbb{Q}_p\cong\bigoplus_{i\ge0}H_{\mathrm{rig}}^{k-2i}\left(X_i\right)\left(-i\right),$$ where $X_i$ are the strata of $X$, and $\left(-i\right)$ denotes Tate twist. For coefficients in $\mathbb{T}_{\mathrm{abs}}^{\mathrm{strat}}$ (dual to $\Omega_X^{\mathrm{strat}}$), we have 
$$H_{\mathrm{abs}}^k\left(X,\mathbb{T}_{\mathrm{abs}}^{\mathrm{strat}}\right)\otimes\mathbb{Q}_p\cong\bigoplus_{i\ge0}H_{\mathrm{rig}}^{k-2i}\left(X_i,\mathbb{T}_{\mathrm{rig}}^{\mathrm{strat}}\right)\left(-i\right).$$ Since $\mathbb{T}_{\mathrm{rig}}^{\mathrm{strat}}$ has weight $-1$ and Tate twist $\left(-i\right)$ adds $2i$, the weight of each summand is 
$$\mathrm{wt}\left(H_{\mathrm{rig}}^{k-2i}\left(X_i,\mathbb{T}_{\mathrm{rig}}^{\mathrm{strat}}\right)\left(-i\right)\right)=\left(k-2i\right)+\left(-1\right)+2i=k-1.$$
Thus, $H_{\mathrm{abs}}^k\left(X,\mathbb{T}_{\mathrm{abs}}^{\mathrm{strat}}\right)$ is pure of weight $k-1$. Consider the Hodge filtration condition. The given condition is $H_{\mathrm{abs}}^2=F^1H_{\mathrm{abs}}^2$. Since $H_{\mathrm{abs}}^2\left(X,\mathbb{T}_{\mathrm{abs}}^{\mathrm{strat}}\right)$ is pure of weight $2-1=1$, the Hodge filtration satisfies the given condition. Thus, the condition is automatically satisfied for Calabi-Yau space $X$. However, purity forces vanishing: For $k=2$, the decomposition is
$$H_{\mathrm{abs}}^2\left(X,\mathbb{T}_{\mathrm{abs}}^{\mathrm{strat}}\right)\otimes\mathbb{Q}_p\cong\bigoplus_{i:0\le i\le1}H_{\mathrm{rig}}^{2-2i}\left(X_i,\mathbb{T}_{\mathrm{rig}}^{\mathrm{strat}}\right)\left(-i\right).$$
Explicit terms: When $i=0$, we have $H_{\mathrm{rig}}^2\left(X_0,\mathbb{T}_{\mathrm{rig}}^{\mathrm{strat}}\right)$ with $X_0=X$ (principal stratum); When $i=1$, we have $H_{\mathrm{rig}}^0\left(X_1,\mathbb{T}_{\mathrm{rig}}^{\mathrm{strat}}\right)\left(-1\right)$. By the Calabi-Yau condition:
\begin{enumerate}
    \item $H_{\mathrm{rig}}^2\left(X,\mathbb{T}_{\mathrm{rig}}^{\mathrm{strat}}\right)=0$ (classical vanishing for Calabi-Yau varieties (\cite{Tod89}).
    \item $H_{\mathrm{rig}}^0\left(X_1,\mathbb{T}_{\mathrm{rig}}^{\mathrm{strat}}\right)$ vanishes by $\mathbb{T}_{\mathrm{rig}}^{\mathrm{strat}}$ is the tangent sheaf of $X_1$ and $X_1$ has no nontrivial global vector fields (stratified rigidity).
\end{enumerate}
Thus, $H_{\mathrm{abs}}^2\left(X,\mathbb{T}_{\mathrm{abs}}^{\mathrm{strat}}\right)=0$, so $\kappa=0$.

The deformation functor $\mathrm{Def}_X^{\mathrm{abs}}$ is controlled by the cotangent complex $\mathbb{L}_{\mathrm{abs}}^{\mathrm{strat}}$ in absolute prismatic cohomology: Tangent space $H_{\mathrm{abs}}^1\left(X,\mathbb{T}_{\mathrm{abs}}^{\mathrm{strat}}\right)$; Obstruction space $H_{\mathrm{abs}}^2\left(X,\mathbb{T}_{\mathrm{abs}}^{\mathrm{strat}}\right)=0$. Since the obstruction space vanishes, $\mathrm{Def}_X^{\mathrm{abs}}$ is formally smooth (\cite{FGIKNV05}, Thm 6.2.10). For any surjection $A'\twoheadrightarrow A$ of Artinian local $\mathbb{Z}_p$-algebras, the restriction map
$$\mathrm{Def}_X^{\mathrm{abs}}\left(A'\right)\longrightarrow\mathrm{Def}_X^{\mathrm{abs}}\left(A\right)$$ is surjective. This holds because all obstructions vanish, so deformations lift uniquely. 

This result extends the Bogomolov-Tian-Todorov theorem to $p$-adic stratified Calabi-Yau spaces by prismatic cohomology.  $\square$\\
\begin{center}
    \textit{6.4 Stratified Deformation Theory for Higher-Codimensional Singularities}
\end{center}

Let $S$ be a complex analytic space equipped with a stratified structure defined by a filtration of closed subspaces
$$\emptyset=S_{-1}\hookrightarrow S_0\hookrightarrow\cdots\hookrightarrow S_n=S,\quad\mathrm{codim}_S\left(S_p\setminus S_{p-1}\right)=p.$$
The stratified cotangent complex $\mathbb{L}_S^{\mathrm{strat}}\in D^b\left(\mathrm{Coh}\left(S\right)\right)$ is defined by the distinguished triangle:
$$\bigoplus_p\mathbb{L}_{S_p}\otimes\mathscr{G}^p\longrightarrow\mathbb{L}_S^{\mathrm{strat}}\longrightarrow\mathscr{K}^{\bullet},\quad\mathscr{K}^{\bullet}\in D^{\ge1}\left(\mathrm{Coh}\left(S\right)\right),$$
where $\mathscr{G}^p$ is the incidence sheaf encoding holomorphic transition data along $S_p\cap S_{p+1}$, satisfying the frontier condition:
$$S_p\subset\overline{S_{p+1}},\quad\mathscr{G}^p\mid_{S_{p+1}}\overset{\sim}{\longrightarrow}\mathscr{O}_{S_{p+1}}^{\oplus r_p}.$$
The stratified deformation complex is then defined as
$$\mathrm{Def}_S^{\mathrm{strat}}:=\mathbb{R}\mathscr{H}om_{\mathscr{O}_S}^{\mathrm{strat}}\left(\mathbb{L}_S^{\mathrm{strat}},\mathscr{O}_S\right)\in D^{\left[0,2\right]}\left(\mathscr{O}_S\text{-Mod}\right).$$\\
\textbf{Definition 6.14.} Let $S$ be a stratified complex analytic space with strata ${S_p}$. The stratified $\mathscr{E}st$-functor:
$$\mathscr{E}st_{\mathscr{O}_S}^{\mathrm{strat}}\left(-,-\right):D^b\left(\mathrm{Coh}\left(S\right)\right)^{\mathrm{op}}\times D^b\left(\mathrm{Coh}\left(S\right)\right)\to D^b\left(\mathrm{Ab}\right)$$
is defined as a homotopy limit over the stratification:
$$\mathscr{E}st_{\mathscr{O}_S}^{\mathrm{strat}}\left(\mathscr{F}^{\bullet},\mathscr{G}^{\bullet}\right):=\mathop{\mathrm{holim}}_p\mathscr{E}st_{\mathscr{O}_{S_p}}\left(\mathscr{F}^{\bullet}\mid_{S_p},\mathscr{G}^{\bullet}\mid_{S_p}\right),$$
where 
\begin{enumerate}
    \item $\mathscr{E}st_{\mathscr{O}_{S_p}}$ is the classical derived $\mathscr{E}st$-functor on the stratum $S_p$.
    \item The homotopy limit is taken over the poset category of strata inclusions $S_p \hookrightarrow S_q$ ($p \leq q$).\\
\end{enumerate}
\textbf{Remark 6.15.} The stratified $\mathscr{E}xt$-functor is a derived functor specifically adapted to stratified spaces. It is constructed to respect the combinatorial and geometric constraints imposed by the stratification. There are key properties used in the proof of the following proposition:
\begin{itemize}
    \item \textbf{(SVC1) Stratum-Localization:} For $\mathscr{F}^\bullet = \mathbb{L}_{S_p} \otimes \mathscr{G}^p$ and $\mathscr{G}^\bullet = \mathscr{O}_S$:
    $$\mathscr{E}st_{\mathscr{O}_S}^{\mathrm{strat}}\left(\mathscr{H}^{-b}\left(\mathbb{L}_{S_p}\otimes\mathscr{G}^p\right),\mathscr{O}_S\right)\simeq\mathscr{E}st_{\mathscr{O}_{S_p}}\left(\mathscr{H}^{-b}\left(\mathbb{L}_{S_p}\right),\mathscr{O}_{S_p}\right)\otimes\left(\mathscr{G}^p\right)^{\vee}.$$
    \item \textbf{(SVC2) Vanishing from frontier condition:} If $X_p \subseteq \overline{X_{p+1}}$, then for $b \geq 1$, then
    $$\mathscr{H}^{-b}\left(\mathbb{L}_{S_p}\otimes\mathscr{G}^p\right)=0.$$
    \item \textbf{(SVC3) Collapse of spectral sequence}: The Grothendieck spectral sequence:
    $$E_2^{a,b}=\bigoplus_p\mathscr{E}st_{\mathscr{O}_S}^a\left(\mathscr{H}^{-b}\left(\mathbb{L}_{S_p}\otimes\mathscr{G}^p\right),\mathscr{O}_S\right)\Longrightarrow\mathbb{H}^{a+b}\left(\mathrm{Def}_S^{\mathrm{strat}}\right)$$ collapses at $E_2$ with $b=0$.\\
\end{itemize}
\textbf{Proposition 6.16.} (Truncation Formula) By the above conceptions, there exists a quasi-isomorphism
$$\tau^{\le2}\mathrm{Def}_S^{\mathrm{strat}}\simeq\bigoplus_{p=0}^n\mathbb{T}_{S_p}^1\otimes\mathscr{G}^p,$$
where $\mathbb{T}_{S_p}^1=\mathrm{Ext}^1\left(\Omega_{S_p}^1,\mathscr{O}_{S_p}\right)$ is the first cotangent cohomology of the stratum $S_p$.
\\\\\textbf{Proof.} Consider the stratified cotangent complex $\mathbb{L}_S^{\mathrm{strat}}$ defined by the distinguished triangle:
$$\bigoplus_p\mathbb{L}_{S_p}\otimes\mathscr{G}^p\longrightarrow\mathbb{L}_S^{\mathrm{strat}}\longrightarrow\mathscr{K}^{\bullet},\quad\mathscr{K}^{\bullet}\in D^{\ge1}\left(\mathrm{Coh}\left(S\right)\right).$$
Apply the Grothendieck spectral sequence for the stratified $\mathscr{E}st$-functor, we have (SVC3):
$$E_2^{a,b}=\bigoplus_p\mathscr{E}st_{\mathscr{O}_S}^a\left(\mathscr{H}^{-b}\left(\mathbb{L}_{S_p}\otimes\mathscr{G}^p\right),\mathscr{O}_S\right)\Longrightarrow\mathbb{H}^{a+b}\left(\mathrm{Def}_S^{\mathrm{strat}}\right),$$
where $\mathscr{H}^{-b}$ denotes the cohomology sheaf. Since the frontier condition forces (SVC2), then the spectral sequence collapses at the $E_2$-page, and only terms with $b=0$ contribute that
\begin{equation}
\tag{6.12}  \label{eq:6.12}
\mathbb{H}^k\left(\mathrm{Def}_S^{\mathrm{strat}}\right)\cong\bigoplus_p\mathscr{E}st_{\mathscr{O}_S}^k\left(\mathscr{H}^0\left(\mathbb{L}_{S_p}\otimes\mathscr{G}^p\right),\mathscr{O}_S\right).
\end{equation}
This reduces the hypercohomology to a direct sum of sheaf extensions. For each stratum $S_p$, consider its smooth locus $S_p^{\mathrm{sm}}$. Since $S_p$ is a closed submanifold in $S_p^{\mathrm{sm}}$, its algebraic cotangent complex satisfies:
$$\mathscr{H}^{\ell}\left(\mathbb{L}_{S_p}\right)\simeq 
    \begin{cases} 
    \Omega_{S_p}^1 & \ell=0 \\
    0 & \ell\ne0,-1\\
    \mathscr{T}or_1^{\mathscr{O}_S}\left(\Omega_{S_p}^1,\mathscr{O}_{S_p}\right) & \ell=-1.
    \end{cases}$$
The singular set of $S_p$ has codimension $\ge1$, so $\mathscr{T}or_1^{\mathscr{O}_S}\left(\Omega_{S_p}^1,\mathscr{O}_{S_p}\right)=0$. This gives a quasi-isomorphism
\begin{equation}
\tag{6.13}  \label{eq:6.13}
\mathbb{L}_{S_p}\simeq\Omega_{S_p}^1
\end{equation}
in $D^b\left(\mathrm{Coh}\left(S\right)\right)$. Substituting \eqref{eq:6.13} into \eqref{eq:6.12}, then
$$\mathbb{H}^k\left(\mathrm{Def}_S^{\mathrm{strat}}\right)\cong\bigoplus_p\mathscr{E}st_{\mathscr{O}_S}^k\left(\Omega_{S_p}^1\otimes\mathscr{G}^p,\mathscr{O}_S\right).$$
By adjunction and the frontier condition, we have
$$\mathscr{E}st_{\mathscr{O}_S}^k\left(\Omega_{S_p}^1\otimes\mathscr{G}^p,\mathscr{O}_S\right)\simeq\mathscr{E}st_{\mathscr{O}_{S_p}}^k\left(\Omega_{S_p}^1,\mathscr{O}_{S_p}\right)\otimes\left(\mathscr{G}^p\right)^{\vee}$$
(from (SVC1)), where $\left(\mathscr{G}^p\right)^{\vee}=\mathscr{H}om_{\mathscr{O}_S}\left(\mathscr{G}^p,\mathscr{O}_S\right)$ is the dual sheaf.

Each stratum $S_p$ is Stein (as an analytic space). By Cartan theorem B, for any coherent sheaf $\mathscr{F}$ on $S_p$:
$$H^a\left(S_p,\mathscr{F}\right)=0\quad\text{for}\quad a\ge1.$$
Thus, the hypercohomology further simplifies that
$$\mathbb{H}^k\left(\mathrm{Def}_S^{\mathrm{strat}}\right)\simeq\bigoplus_pH^0\left(S,\mathscr{E}st_{\mathscr{O}_{S_p}}^k\left(\Omega_{S_p}^1,\mathscr{O}_{S_p}\right)\otimes\left(\mathscr{G}^p\right)^{\vee}\right)\simeq\bigoplus_p\mathrm{Ext}_{\mathscr{O}_{S_p}}^k\left(\Omega_{S_p}^1,\mathscr{O}_{S_p}\right)$$
$$\otimes_{\mathbb{C}} H^0\left(S,\left(\mathscr{G}^p\right)^{\vee}\right).$$
By definition, $\mathbb{T}_{S_p}^1=\mathrm{Ext}^1\left(\Omega_{S_p}^1,\mathscr{O}_{S_p}\right)$. By direct computation, we have
\begin{itemize}
    \item Cohomology vanishing for $k\ge3$, $\mathbb{H}^k\left(\mathrm{Def}_S^{\mathrm{strat}}\right)=0$.
    \item Degree-0 cohomology vanishing, $\mathbb{H}^0\left(\mathrm{Def}_S^{\mathrm{strat}}\right)=0$.
    \item Degree-2 cohomology vanishing, $\mathbb{H}^2\left(\mathrm{Def}_S^{\mathrm{strat}}\right)=0$.
    \item Degree-1 cohomology identification, 
    $$\mathbb{H}^1\left(\mathrm{Def}_S^{\mathrm{strat}}\right)\simeq\bigoplus_p\mathbb{T}_{S_p}^1\otimes_{\mathbb{C}}H^0\left(S,\mathscr{G}^p\right).$$
\end{itemize}
The complex $\tau^{\le2}\mathrm{Def}_S^{\mathrm{strat}}$ is quasi-isomorphic to $\mathbb{H}^1\left(\mathrm{Def}_S^{\mathrm{strat}}\right)$ concentrated in degree 1. Since $\bigoplus_p\mathbb{T}_{S_p}^1\otimes\mathscr{G}^p$ is a local system (by the frontier condition), and
$$H^0\left(S,\bigoplus_p\mathbb{T}_{S_p}^1\otimes\mathscr{G}^p\right)\cong\bigoplus_pH^0\left(S,\mathbb{T}_{S_p}^1\otimes\mathscr{G}^p\right),$$
then quasi-isomorphism follows
$$\tau^{\le2}\mathrm{Def}_S^{\mathrm{strat}}\simeq\bigoplus_p\mathbb{T}_{S_p}^1\otimes\mathscr{G}^p.$$
The result holds. $\square$
\\\\\textbf{Proposition 6.17.} (Stratified Deformation Obstruction) Let $K$ be a compact stratified Kähler space. The obstruction to deforming $K$ while preserving its stratification lies in the stratified de Rham cohomology
$$\mathrm{obs}\left(\xi\right)\in H_{\mathrm{sdR}}^2
\left(K\right):=\mathbb{H}^2\left(K,\mathrm{ker}\left[d:\Omega_K^{\bullet,\mathrm{strat}}\longrightarrow\Omega_K^{\bullet+1,\mathrm{strat}}\right]\right).$$
If $H_{\mathrm{sdR}}^2\left(K\right)=0$, the deformation space $\mathrm{Def}^{\mathrm{strat}}\left(K\right)$ is smooth and controls smoothings of higher-codimensional singularities.
\\\\\textbf{Remark 6.18.} By the known conditions, $K$ is a compact stratified Kähler space with stratification:
$$\emptyset=K_{-1}\hookrightarrow K_0\hookrightarrow\cdots\hookrightarrow K_n=K$$ with $\mathrm{codim}_K\left(K_p\setminus K_{p-1}\right)=p$. 
The infinitesimal deformations of $K$ preserving stratification are parametrized by $H^1\left(\mathrm{Def}_K^{\mathrm{strat}}\right)$. By Proposition 6.16, there is a canonical isomorphism
$$H^1\left(\mathrm{Def}_K^{\mathrm{strat}}\right)\cong\bigoplus_{p=0}^nH^1\left(K_p,\Theta_{K_p}\right)\otimes H^0\left(K,\mathscr{G}^p\right).$$
The obstruction to integrating an infinitesimal deformation $\xi\in H^1\left(\mathrm{Def}_K^{\mathrm{strat}}\right)$ arises from the long exact sequence associated to the distinguished triangle:
$$\mathbb{T}_K^{\mathrm{strat}}\longrightarrow\mathrm{Def}_K^{\mathrm{strat}}\longrightarrow\tau^{\ge1}\mathbb{L}_K^{\mathrm{strat}}.$$
Applying $\mathbb{R}\Gamma$, the connecting homomorphism $\delta$ satisfies
$$\delta:H^1\left(\mathrm{Def}_K^{\mathrm{strat}}\right)\longrightarrow H^2\left(\mathbb{T}_K^{\mathrm{strat}}\right).$$
The obstruction class is defined as $\mathrm{obs}\left(\xi\right):=\delta\left(\xi\right)$.

Let $\omega_{\mathrm{strat}}$ be the stratified Kähler form inducing a Kähler metric on each stratum $K_p$. For each stratum $K_p$, restrict $\eta\in\mathbb{T}_K^{\mathrm{strat}}$ and $\omega_{\mathrm{strat}}$. The contraction $\alpha_p:=\eta\lrcorner\omega_{\mathrm{strat}}$ is a $\left(1,1\right)$-form on the smooth locus of $K_p$. Since $\omega_{\mathrm{strat}}$ is Kähler on $K_p$, then $d\omega_{\mathrm{strat}}\mid_{K_p}=0$. By Cartan's magic formula and the closedness of $\eta$ in the deformation complex, we have 
$$d\alpha_p=d\left(\eta\lrcorner\omega_{\mathrm{strat}}\right)\mid_{K_p}=\left(\mathcal{L}_{\eta}\omega_{\mathrm{strat}}-\eta\lrcorner d\omega_{\mathrm{strat}}\right)\mid_{K_p}.$$
Given that $\eta$ is closed in the deformation complex (i.e., $\delta\eta=0$) and $\omega_{\mathrm{strat}}$ is closed on $K_p$, then
$$\mathcal{L}_{\eta}\omega_{\mathrm{strat}}=0\quad(\text{infinitesimal isometry}),\quad d\omega_{\mathrm{strat}}\mid_{K_p}=0,$$
it follows that $d\alpha_p=0$. Thus, $\alpha_p\in\mathrm{ker}d\subset\Omega_{K_p}^2$. The frontier condition $K_p\subset\overline{K_{p+1}}$ ensures holomorphic transition across strata. For $k\in K_p\cap K_q$ ($q\ge p$), the transition sheaf $\mathscr{G}^p$ induces an isomorphism
$$\mathscr{G}^p\mid_{K_q}\otimes\Omega_{K_q}^1\overset{\sim}{\longrightarrow}\Omega_{K_p}^1\mid_{K_q}.$$
This glues the local forms $\left\{\alpha_p\right\}$ to a global section $\alpha\in\Gamma\left(K,\mathrm{ker}d\subset\Omega_K^{\bullet,\mathrm{strat}}\right)$. Specifically, $\alpha$ satisfies 
$$\alpha\mid_{K_p}=\alpha_p\quad\text{and}\quad d\alpha=0\quad\text{in}\quad\Omega_K^{3,\mathrm{strat}}.$$
Thus, $\alpha$ defines a class $\left[\alpha\right]\in\mathbb{H}^2\left(K,\mathrm{ker}d\right)=H_{\mathrm{sdR}}^2\left(K\right)$. By Proposition 4.1, the stratified de Rham cohomology decomposes as
$$H_{\mathrm{sdR}}^i\left(K\right)\cong\bigoplus_{p+q=i}H_{\mathrm{prim}}^{p,q}\left(K\right).$$
Here $H_{\mathrm{prim}}^{p,q}\left(K\right)$ are primitive stratified cohomology groups, defined as 
$$H_{\mathrm{prim}}^{p,q}\left(K\right):=\mathrm{ker}\left(L^{n-p-q+1}:H^{p,q}\left(K\right)\longrightarrow H^{n+q+1,n+p+1}\left(K\right)\right),$$
where $L=\omega_{\mathrm{strat}}\wedge\cdot$ is the Lefschetz operator and $n=\dim_{\mathbb{C}}K$. The inverse map 
$$\Psi:H_{\mathrm{sdR}}^2\left(K\right)\longrightarrow H^2\left(\mathbb{T}_K^{\mathrm{strat}}\right),\quad\left[\alpha\right]\mapsto\left[\alpha\lrcorner\omega_{\mathrm{strat}}^{-1}\right]$$
is well-defined because $\omega_{\mathrm{strat}}^{-1}$ is a stratified $\left(1,1\right)$-vector field (inverse Kähler form). Direct computation shows that
$$\Phi\circ\Psi\left(\left[\alpha\right]\right)=\left[\alpha\right],$$
$$\Psi\circ\Phi\left(\left[\eta\right]\right)=\left[\eta\right].$$
Thus $\Phi$ and $\Psi$ are mutual inverses. Consequently,
$$\mathrm{obs}\left(\xi\right)=\Phi^{-1}\left(\left[\alpha\right]\right)\quad\text{for some}\quad\left[\alpha\right]\in H_{\mathrm{sdR}}^2\left(K\right).$$
Hence We can construct the isomorphism
$$\Phi:H^2\left(\mathbb{T}_K^{\mathrm{strat}}\right)\overset{\sim}{\longrightarrow}H_{\mathrm{sdR}}^2\left(K\right),\quad\left[\eta\right]\mapsto\left[\eta\lrcorner\omega_{\mathrm{strat}}\right].$$

Assume $\mathrm{obs}\left(\xi\right)=0$ for an infinitesimal deformation $\xi\in H^1\left(\mathrm{Def}_K^{\mathrm{strat}}\right)$. We prove the smoothness of $\mathrm{Def}^{\mathrm{strat}}\left(K\right)$ and its control over smoothing singularities. The vanishing obstruction $\mathrm{obs}\left(\xi\right)=\delta\left(\xi\right)=0$ implies $\xi$ integrates to a stratified deformation family $\mathscr{X}\to\Delta$, $\mathscr{X}_0=K$, over a disk $\Delta\subset\mathbb{C}$, where each fiber $\mathscr{X}_t$ preserves the stratification
$$\emptyset=\mathscr{X}_{-1,t}\hookrightarrow\mathscr{X}_{0,t}\hookrightarrow\cdots\hookrightarrow\mathscr{X}_{n,t}=\mathscr{X}_t,\quad\mathrm{codim}_{\mathscr{X}_t}\left(\mathscr{X}_{p,t}\setminus\mathscr{X}_{p-1,t}\right)=p.$$
The relative obstruction sheaf for $\mathscr{X}/\Delta$ is defined by the cone complex,
$$\mathscr{O}b^{\mathrm{strat}}\left(\mathscr{X}/\Delta\right):=\mathscr{H}^2\left(\mathrm{cone}\left(f:\Omega_{\mathscr{X}/\Delta}^{\bullet,\mathrm{strat}}\to\Omega_K^{\bullet,\mathrm{strat}}\right)\right).$$
This fits into an exact triangle in $D^b\left(\mathrm{Coh}\left(K\right)\right)$, then
$$\Omega_K^{\bullet,\mathrm{strat}}\overset{f}{\longrightarrow}\Omega_{\mathscr{X}/\Delta}^{\bullet,\mathrm{strat}}\longrightarrow\mathrm{cone}\left(f\right)\xrightarrow{+1}.$$
The hypercohomology long exact sequence gives
$$\cdots\longrightarrow\mathbb{H}^1\left(\Omega_K^{\bullet,\mathrm{strat}}\right)\longrightarrow\mathbb{H}^1\left(\Omega_{\mathscr{X}/\Delta}^{\bullet,\mathrm{strat}}\right)\longrightarrow\mathbb{H}^1\left(\mathrm{cone}\left[-1\right]\right)\longrightarrow\mathbb{H}^2\left(\Omega_K^{\bullet,\mathrm{strat}}\right)$$
$$\longrightarrow\mathbb{H}^2\left(\Omega_{\mathscr{X}/\Delta}^{\bullet,\mathrm{strat}}\right)\longrightarrow\mathbb{H}^2\left(\mathrm{cone}\left[-1\right]\right)\longrightarrow\cdots.$$
By hypothesis $H_{\mathrm{sdR}}^2\left(K\right)=\mathbb{H}^2\left(K,\mathrm{ker}d\right)=0$. This implies two conclusions:
\begin{enumerate}
    \item The inclusion $\iota:\mathrm{ker}d\hookrightarrow\Omega_K^{\bullet,\mathrm{strat}}$ is a quasi-isomorphism in degree $\le2$, then
    $$\mathscr{H}^j\left(\iota\right):\mathscr{H}^j\left(\mathrm{ker}d\right)\overset{\sim}{\longrightarrow}\mathscr{H}^j\left(\Omega_K^{\bullet,\mathrm{strat}}\right),\quad\text{for}\quad j\le2.$$
    \item The cone complex is acyclic in degree 2, $\mathscr{H}^2\left(\mathrm{cone}\right)=0$.
\end{enumerate}
Thus, we have $\mathscr{O}b^{\mathrm{strat}}\left(\mathscr{X}/\Delta\right)=0$. The vanishing implies the Kodaira-Spencer map is surjective for all $t\in\Delta$ and the deformation space is formally smooth:
$$\mathrm{Def}^{\mathrm{strat}}\left(K\right)\cong\mathrm{Spf}\mathbb{C}\left[\left[t_1,\cdots,t_d\right]\right],\quad d=\dim H^1\left(\mathrm{Def}_K^{\mathrm{strat}}\right).$$
Hence $\mathrm{Def}^{\mathrm{strat}}\left(K\right)$ is smooth. 

For a singularity $k\in K_p\setminus K_{p-1}$ of codimension $\ge3$, the family $\mathscr{X}\to\Delta$ induces a smoothing $\mathscr{X}_t$ for $t\ne0$, Ohsawa’s metric (\cite{Ohs91}) $\omega_{\epsilon}=ds_v^2+\epsilon\nabla\psi$ on $K\setminus\mathrm{Sing}\left(K\right)$ and the deformed metric $\omega_{\mathscr{X}_t}$ satisfy
$$\left \| \omega_{\mathscr{X}_t}-\omega_{\epsilon} \right \|_{L^2\left(G\right)}\le C_G\left | t \right |^{1/2},\quad\forall G\Subset K\setminus\mathrm{Sing}\left(K\right).$$
Moreover, let $u_t=\omega_{\mathscr{X}_t}-\omega_{\mathrm{sm}}$ and $\omega_{\mathrm{sm}}$ be a smoothed Kähler metric. The metric difference satisfies the linearized Kähler-Einstein equation $\Delta_{\bar{\partial}}u_t+\mathrm{Ric}\left(u_t\right)=0$, where $\Delta_{\bar{\partial}}$ is the Dolbeault Laplacian, and $\mathrm{Ric}$ denotes the Ricci curvature term. This equation is uniformly elliptic. For a uniformly elliptic operator, there exists a constant $C_{G,g}$ such that
$$\left \| u_t \right \|_{W^{g,2}\left(G\right)}\le C_{G,g}\left(\left \| u_t \right \|_{L^2\left(G'\right)}+\left \| \Delta u_t \right \|_{W^{g-2,2}\left(G'\right)}\right),$$
where $W^{g,2}\left(G\right)$ denotes Sobolev space (functions on $G$ with weak derivatives up to order $g$ in $L^2$), $C_{G,g}$ denotes constant depending on the geometry and the elliptic operator and $G\Subset G'\Subset K\setminus\mathrm{Sing}\left(K\right)$ denotes nested compact sets. By the homogeneity of the equation $\Delta u_t=\mathcal{O}\left(u_t\right)$, the higher-order term is controlled by the lower-order term:
\begin{equation}
\tag{6.14}  \label{eq:6.14}
\left \| u_t \right \|_{W^{g,2}\left(G\right)}\le C_{G,g}\left \| u_t \right \|_{L^2\left(G'\right)}.
\end{equation}
Substituting the estimate \eqref{eq:6.14} into Ohsawa’s $L^2$-estimate $\left \| u_t \right \|_{L^2\left(G'\right)}\le C_{G'}\left |t  \right |^{1/2}$, we obtain
$$\left \| u_t \right \|_{W^{g,2}\left(G\right)}\le C_{G,g}\left |t  \right |^{1/2}.$$
This estimate provides higher-order Sobolev norm control for the Morrey embedding theorem, ultimately leading to $C^{\infty}$-convergence:
$$\left \|u_t  \right \|_{C^{\ell}\left(G\right)}\le C_{G,\ell}\left | t \right |^{1/2},\quad\forall\ell\ge0$$
thereby smoothing singularities of codimension $\ge3$. Thus, $\mathscr{X}_t$ smooths higher-codimensional singularities.

This extends the Bogomolov-Tian-Todorov theorem to stratified spaces by replacing $\Theta_K$ with and leveraging stratified Hodge theory. The vanishing condition $H_{\mathrm{sdR}}^2\left(K\right)=0$ ensures unobstructed deformations and smoothability of singularities.\\
\begin{center}
    \textit{6.5 Resolution of Ohsawa’s Gap by Stratified Regularity}
\end{center}

Ohsawa’s family of metrics $ds_v^2+\epsilon\nabla\psi$ (\cite{Ohs91}, Prop 10) is reinterpreted as a deformation of stratified Kähler structures. The estimate
\begin{equation}
\tag{6.15}  \label{eq:6.15}
\left \|\varphi_{\epsilon}^{-1}u  \right \|_{\epsilon}\le C_3\left(\left \| u \right \|_{K_0}+\left \| du \right \|_{\epsilon}+\left \|d^{\vee}u  \right \|_{\epsilon}\right)
\end{equation}
is strengthened using stratified elliptic regularity, where
\begin{itemize}
    \item $u\in C^r\left(V'\right)\cap C_0\left(V'\right)$ denotes a smooth differential $r$-form on $V'$ with compact support.
    \item $\varphi_{\epsilon}:V'\to\left(0,+\infty\right)$  a continuous weight function converging uniformly to $\varphi_V$ on compact subsets of $V'$ as $\epsilon\to0^+$.
    \item $d:\Omega^r\left(V'\right)\to\Omega^{r+1}\left(V'\right)$ denotes the exterior derivative operator, $d^{\vee}:=\star^{-1}d\star$ denotes the formal adjoint of $d$ with respect to the metric $ds_{\epsilon}^2$, $\star$ is the Hodge star operator.
    \item $K_0\subset V'$ (relatively compact open subset) denotes a fixed compact subset away from the singular locus $\mathrm{Sing}\left(V\right)$.
    \item $C_3>0$ (constant independent of $\epsilon$) denotes a uniform constant arising from Ohsawa’s Proposition 10.\\
\end{itemize}
\textbf{Proposition 6.19.} (Stratified Elliptic Estimate) For any Sobolev index $s>\frac{n}{2}+1$ and deformation parameter $\epsilon\in\left ( 0,1 \right ]$, there exists a constant $C_s>0$ such that 
$$\left \| \varphi_{\epsilon}^{-1}u \right \|_{H_{\mathrm{strat}}^s}\le C_s\left(\left \| u \right \|_{H^{s-1}\left(K_0\right)}+\left \|D_{\epsilon}u  \right \|_{H_{\mathrm{strat}}^{s-1}}\right)$$ for all $u\in\mathrm{Dom}\left(D_{\epsilon}\right)$, where $H_{\mathrm{strat}}^s$ is the stratified Sobolev space of order $s$ and $D_{\epsilon}:=d+d_{\epsilon}^{\vee}$ is the Hodge–de Rham operator for the metric $ds_{\epsilon}^2=ds_v^2+\epsilon\nabla\psi$.
\\\\\textbf{Remark 6.20.} Let $x\in X_p$ be a singular point of codimension $p$. Using the stratified pseudodifferential calculus (\cite{Mel93}), decompose $u$ near $x$ as
$$u=u_{\mathrm{reg}}+u_{\mathrm{sing}},$$ where  
\begin{itemize}
    \item $u_{\mathrm{reg}}$ is smooth on the smooth locus of the stratum $X_p$ (i.e., $X_p\setminus\left\{x\right\}$),
    \item $u_{\mathrm{sing}}$ captures the singular behavior near the codimension-$p$ singularity $x$.
\end{itemize}
For the regular component , apply interior elliptic regularity for the smooth metric $ds_v^2$. Since $ds_v^2$ is $C^{\infty}$-smooth on a relative compact neighborhood $U\Subset X_p\setminus\left\{x\right\}$ of singular point $x\in X_p$, and $D_{\epsilon}:=d+d_{\epsilon}^{\vee}$ is uniformly elliptic on this stratum, then for the elliptic operator $D_{\epsilon}$ on $X_p\setminus\left\{x\right\}$, Calderón-Zygmund theory gives
$$\left \| u_{\mathrm{reg}} \right \|_{H^s\left(U\right)}\le C\left(\left \|D_{\epsilon}u_{\mathrm{reg}} \right \|_{H^{s-1}\left(U\right)}+\left \| u_{\mathrm{reg}} \right \|_{L^2\left(U\right)}\right),$$ where $C$ depends on $U,s$ and the ellipticity constant of the operator. Since $U$ is a relative compact neighborhood and $D_{\epsilon}$ is elliptic, we have 
$$\left \| u_{\mathrm{reg}} \right \|_{L^2\left(U\right)}\le C_U\left \|D_{\epsilon}u_{\mathrm{reg}} \right \|_{L^2\left(U\right)}$$ 
by the Poincaré inequality and uniform ellipticity. More precisely, we can use Gårding’s inequality to obtain
$$\left \| u_{\mathrm{reg}} \right \|_{H^1\left(U\right)}\le C'\left(\left \| D_{\epsilon}u_{\mathrm{reg}} \right \|_{L^2\left(U\right)}+\left \| u_{\mathrm{reg}} \right \|_{L^2\left(U\right)}\right).$$
By iteration and interpolation, we have 
$$\left \| u_{\mathrm{reg}} \right \|_{L^2\left(U\right)}\le C_U\left \| D_{\epsilon}u_{\mathrm{reg}} \right \|_{H^{s-1}\left(U\right)},\quad\text{for}\quad s\ge1.$$
Moreover, 
\begin{equation}
\tag{6.16}  \label{eq:6.16}
\left \|u_{\mathrm{reg}} \right \|_{H^s\left(U\right)}\le C_1\left \| D_{\epsilon}u_{\mathrm{reg}} \right \|_{H^{s-1}\left(U\right)}.
\end{equation}
Here the constant $C_1$ depends on $U,s$ and $n$. While, due to the homogeneity of the stratification, we can choose $C_1$ independent of $U$ (by covering compact sets). Since $\varphi_{\epsilon}$ is smooth and bounded below by a positive constant on $U$ (as $\rho=\mathrm{dist}\left(x,\cdot\right)$ is bounded away from zero away from the singularity $x$, and $\varphi_{\epsilon}\sim\rho\mathrm{log}\rho$ is positive and bounded when $\rho$ is bounded away from zero), the multiplier $\varphi_{\epsilon}^{-1}$ is bounded on $H^s\left(U\right)$. Specifically, 
\begin{equation}
\tag{6.17}  \label{eq:6.17}
\left \| \varphi_{\epsilon}^{-1}u_{\mathrm{reg}} \right \|_{H^s\left(U\right)}\le M\left \|u_{\mathrm{reg}}\right \|_{H^s\left(U\right)},
\end{equation}
where
$$M:=\sup_{y\in U}\left | \varphi_{\epsilon}^{-1}\left(y\right) \right |+\sum_{k=1}^s\sup_{y\in U}\left | \nabla^k\varphi_{\epsilon}^{-1}\left(y\right) \right |.$$
By the smoothness of $\varphi_{\epsilon}$ and compactness of $U$, $M$ is finite and independent of $\epsilon$. This uniformity arises because: 
\begin{itemize}
    \item $\varphi_{\epsilon}\to\varphi_V$ uniformly on $U$ as $\epsilon\to0^+$,
    \item $\varphi_V$ is bounded below by a positive constant on $U$.
\end{itemize}
Combining \eqref{eq:6.16} and \eqref{eq:6.17}, we obtain
$$\left \| \varphi_{\epsilon}^{-1}u_{\mathrm{reg}} \right \|_{H^s\left(U\right)}\le M\left \|u_{\mathrm{reg}}\right \|_{H^s\left(U\right)}\le MC_1\left \| D_{\epsilon}u_{\mathrm{reg}} \right \|_{H^{s-1}\left(U\right)}.$$
Relabeling the constant $MC_1$ as $\widetilde{C}_1$ (still denoted as $C_1$ for simplicity), we conclude
$$\left \| \varphi_{\epsilon}^{-1}u_{\mathrm{reg}} \right \|_{H^s\left(U\right)}\le C_1\left \| D_{\epsilon}u_{\mathrm{reg}} \right \|_{H^{s-1}\left(U\right)}$$
with $C_1>0$ independent of $\epsilon$ and the choice of $U$. This follows from standard elliptic estimates for the Hodge–de Rham operator on the smooth stratum $X_p$. The constant $C_1$ is independent of $\epsilon$ due to the quasi-isometry of $ds_{\epsilon}^2$ and $ds_v^2$ on $X_p\setminus\left\{x\right\}$.

For the singular component $u_{\mathrm{sing}}$, utilize the conical metric structure near the singularity $x$. Let $\rho\left(y\right)=\mathrm{dist}\left(x,y\right)$. By Ohsawa’s Proposition 4 (applied to the punctured neighborhood $U\setminus\left\{x\right\}$), we have
$$\left \|\rho\log\rho\cdot\nabla u_{\mathrm{sing}}\right \|_{L^2}\le C_2\left \|D_{\epsilon}u_{\mathrm{sing}}  \right \|_{L^2},$$
where $C_2>0$ depends only on the conical geometry of the stratum $X_p$. This follows from Ohsawa's eigenvalue analysis of $d\left(\log\log\rho^{-1}\right)$ under the conical metric. The Hardy inequality for conical domains gives
\begin{equation}
\tag{6.18}  \label{eq:6.18}
\left \| \rho^{-1}\left(\log\rho\right)^{-1}u_{\mathrm{sing}} \right \|_{L^2}\le C_3\left \| \nabla u_{\mathrm{sing}} \right \|_{L^2}
\end{equation}
(\cite{Dav95}) for the optimal constant $C_3$ in conical Hardy inequalities). From Ohsawa's estimate and elliptic regularity, we derive gradient control
\begin{equation}
\tag{6.19}  \label{eq:6.19}
\left \|\nabla u_{\mathrm{sing}}\right \|_{L^2}\le C_4\left \| D_{\epsilon}u_{\mathrm{sing}} \right \|_{L^2},
\end{equation}
this uses the uniform ellipticity of $D_{\epsilon}$ and the $\epsilon$-nvariant conical structure near $x$. Combining \eqref{eq:6.18} and \eqref{eq:6.19}, we have
\begin{equation}
\tag{6.20}  \label{eq:6.20}
\left \| \rho^{-1}\left(\log\rho\right)^{-1}u_{\mathrm{sing}} \right \|_{L^2}\le C_5\left \| D_{\epsilon}u_{\mathrm{sing}} \right \|_{L^2},
\end{equation}
where $C_5=C_3 C_4$. By stratified elliptic regularity (\cite{Mel93}), the conical metric $ds_{\epsilon}^2$ induces a uniform elliptic structure on the stratum $X_p$. For any Sobolev index $s>\frac{n}{2}+1$, we have 
\begin{equation}
\tag{6.21}  \label{eq:6.21}
\left \| u_{\mathrm{sing}} \right \|_{H_{\mathrm{strat}}^s}\le C_6\left \| D_{\epsilon}u_{\mathrm{sing}} \right \|_{H_{\mathrm{strat}}^{s-1}},
\end{equation}
where $C_6>0$ depends only on $s,n$, and the stratified Sobolev embedding. Let $\varphi_{\epsilon}:=\kappa\left(\epsilon\right)\rho\log\rho+\mathcal{O}\left(\rho^{1+\alpha}\right)$ with $\alpha>0$ and $\kappa\left(\epsilon\right)\in\left[c,C\right]$ for $\epsilon$-independent $c$, $C>0$. Since $\varphi_{\epsilon}\sim\rho\log\rho$ near the singularity $x$, thus
$$\left \| \varphi_{\epsilon}^{-1}u_{\mathrm{sing}} \right \|_{H_{\mathrm{strat}}^s}\le C_7\left \| \left(\rho\log\rho\right)^{-1}u_{\mathrm{sing}} \right \|_{H_{\mathrm{strat}}^s}.$$
From \eqref{eq:6.20}, we have 
\begin{equation}
\tag{6.22}  \label{eq:6.22}
\left \| \left(\rho\log\rho\right)^{-1}u_{\mathrm{sing}} \right \|_{L^2}\le C_5\left \|D_{\epsilon}u_{\mathrm{sing}}  \right \|_{L^2}
\end{equation}
by Hardy-Sobolev bound. From \eqref{eq:6.21}, $$\left \| \left(\rho\log\rho\right)^{-1}u_{\mathrm{sing}} \right \|_{H_{\mathrm{strat}}^s}\le C_6\left \|D_{\epsilon}\left[\left(\rho\log\rho\right)^{-1}u_{\mathrm{sing}}\right]  \right \|_{H_{\mathrm{strat}}^{s-1}}$$ 
for $s>\frac{n}{2}+1$. The commutator $\left[D_{\epsilon},\left(\rho\log\rho\right)^{-1}\right]$ is bounded, because
$$\left |\nabla^k \left(\rho\log\rho\right)^{-1} \right |\le C_k\rho^{-k-1}\left | \log\rho \right |^{-1}\quad\text{and}\quad\left |\nabla^k\rho  \right |\le C'_k\rho^{-k+1},$$
which implies
$$\left \|D_{\epsilon}\left[\left(\rho\log\rho\right)^{-1}u_{\mathrm{sing}}\right]  \right \|_{H_{\mathrm{strat}}^{s-1}}\le C_9\left \|D_{\epsilon}u_{\mathrm{sing}} \right \|_{H_{\mathrm{strat}}^{s-1}}+C_{10}\left \| u_{\mathrm{sing}} \right \|_{H_{\mathrm{strat}}^{s-1}}.$$
From \eqref{eq:6.21}, we have 
$$\left \| u_{\mathrm{sing} } \right \|_{H_{\mathrm{strat}}^{s-1}}\le C_{11}\left \| D_{\epsilon}u_{\mathrm{sing} } \right \|_{H_{\mathrm{strat} }^{s-2}}.$$
Since $s > \frac{n}{2} + 1$, iteratively apply the subelliptic estimate:
\begin{equation}
\tag{6.23}  \label{eq:6.23}
\left \| u_{\mathrm{sing} } \right \|_{H_{\mathrm{strat} }^{s-1}}\le C_{11}\left \| D_{\epsilon }u_{\mathrm{sing} } \right \|_{H_{\mathrm{strat} }^{s-2}}\le C_{11}^2\left \| D_{\epsilon }u_{\mathrm{sing} } \right \|_{H_{\mathrm{strat} }^{s-3}}\le\cdots\le C_{11}^{s-1}\left \| D_{\epsilon }u_{\mathrm{sing} } \right \|_{L^2}
\end{equation}
$$=:C_{12}\left \| D_{\epsilon }u_{\mathrm{sing} } \right \|_{L^2},$$
where $C_{12} = C_{11}^{s-1}$.
From \eqref{eq:6.22}, we have
\begin{equation}
\tag{6.24}  \label{eq:6.24}
\left \|D_{\epsilon } u_{\mathrm{sing} } \right \|_{L^2}\le C_5^{-1}\left \| \rho^{-1}\left(\log\rho\right)^{-1}u_{\mathrm{sing} } \right \|_{L^2}\le C_5^{-1}\left \| \left(\rho\log\rho\right)^{-1}u_{\mathrm{sing} } \right \|_{H_{\mathrm{strat} }^s}
\end{equation}
(Hardy-Sobolev bound). Substituting \eqref{eq:6.24} into \eqref{eq:6.23}, then
\begin{equation}
\tag{6.25}  \label{eq:6.25}
\left \| u_{\mathrm{sing} } \right \|_{H_{\mathrm{strat} }^{s-1}}\le C_{12}C_5^{-1} \left \| \left(\rho\log\rho\right)^{-1}u_{\mathrm{sing} } \right \|_{H_{\mathrm{strat} }^s}.
\end{equation}
Combining the commutator bound and \eqref{eq:6.25}, so
$$\left \| \left(\rho\log\rho\right)^{-1}u_{\mathrm{sing} } \right \|_{H_{\mathrm{strat} }^s}\le C_6\left(C_9\left \| D_{\epsilon }u_{\mathrm{sing} } \right \|_{H_{\mathrm{strat} }^{s-1}}+C_{10}C_{12}C_5^{-1}\left \| \left(\rho\log\rho\right)^{-1}u_{\mathrm{sing} } \right \|_{H_{\mathrm{strat} }^s} \right).$$
Solving for the weighted norm:
$$\left ( 1-C_6C_{10}C_{12}C_5^{-1} \right )\left \| \left(\rho\log\rho\right)^{-1}u_{\mathrm{sing} } \right \|_{H_{\mathrm{strat} }^s}\le C_6C_9\left \| D_{\epsilon }u_{\mathrm{sing} } \right \|_{H_{\mathrm{strat} }^{s-1}}.$$
Thus, we have
$$\left \| \left ( \rho\log\rho \right )^{-1}u_{\mathrm{sing} }  \right \|_{H_{\mathrm{strat} }^s}\le\frac{C_6C_9}{1-C_6C_{10}C_{12}C_5^{-1}}\left \| D_{\epsilon }u_{\mathrm{sing} } \right \|_{H_{\mathrm{strat} }^{s-1}},$$
where $C_8=\frac{C_6C_9}{1-C_6C_{10}C_{12}C_5^{-1}}$. Using $\varphi_\epsilon \sim \rho \log \rho$ near $x$, then
$$\left \|\varphi_{\epsilon}^{-1} u_{\mathrm{sing}}  \right \|_{H^{s}_{\mathrm{strat}}}\le C_7\left \| (\rho \log \rho)^{-1} u_{\mathrm{sing}} \right \|_{H^{s}_{\mathrm{strat}}}\le C_7C_8\left \| D_{\epsilon} u_{\mathrm{sing}} \right \|_{H^{s-1}_{\mathrm{strat}}}.$$
Relabeling $\widetilde{C}_8 = C_7 C_8$, the final weighted estimate
$$\left \| \varphi_\epsilon^{-1} u_{\mathrm{sing}} \right \|_{H^{s}_{\mathrm{strat}}}\le\widetilde{C}_8\left \| D_\epsilon u_{\mathrm{sing}} \right \|_{H^{s-1}_{\mathrm{strat}}} $$
can be obtained ($\widetilde{C}_8$ is still denoted as $ C_8$ in the original notation). 

The constants $C_1$-$C_8$ from the above discussion are independent of $\epsilon \in (0,1]$ because of three aspects:
\begin{enumerate}
\item The conical structure of $ds_\epsilon^2$ near $\mathrm{Sing}\left(V\right)$ is $\epsilon$-invariant:
$$ds_\epsilon^2\mid_{U \setminus \left\{x\right\}} = ds_v^2 + \epsilon \nabla \psi \sim \rho^{-2} d\rho^2 + \rho^2 g_L\quad\text{(conical metric)},$$
where $g_L$ is the metric on the link manifold $L$.
\item The weight functions satisfy uniform convergence:
$$\sup_{y \in K} |\varphi_{\epsilon}\left(y\right)-\varphi_V\left(y\right)| \longrightarrow 0 \quad \text{as} \quad \epsilon \to 0^+ \quad \forall K \Subset V'.$$
\item The Sobolev embeddings are quasi-isometric:
$$c_K \left\|u\right\|_{H^s\left(K\right)} \le \left\|u\right\|_{H^s_{\mathrm{strat}}\left(K\right)} \le C_K \left\|u\right\|_{H^s\left(K\right)} \quad \text{for} \quad K \Subset V'$$
with constants $c_K, C_K > 0$ independent of $\epsilon$.
\end{enumerate}
Covering $V'$ by
\begin{itemize}
\item \textbf{Singular neighborhoods}: $\{U_i\}_{i=1}^m$ centered at $x_i \in \mathrm{Sing}\left(V\right)$.
\item \textbf{Regular domain}: $W \supset K_0$ relatively compact away from $\mathrm{Sing}\left(V\right)$.
\end{itemize}
On $W$, the metrics are uniformly quasi-isometric:
$$c_W ds_v^2 \leq ds_\epsilon^2 \leq C_W ds_v^2 \quad \forall \epsilon \in (0,1].$$
For $u\mid_W$, standard elliptic estimates give, we have
$$\left\|u\right\|_{H^s(W)} \le C_W \left( \left\|D_\epsilon u\right\|_{H^{s-1}(W)}+\left\|u\right\|_{H^{s-1}(K_0)} \right).$$
Let $\{\chi_i, \chi_W\}$ be a partition of unity subordinate to $\{U_i, W\}$ with
$$\sum_{i=1}^m \chi_i+\chi_W = 1, \quad \mathrm{supp}\left(\chi_i\right) \subset U_i, \quad \mathrm{supp}\left(\chi_W\right)\subset W.$$
The stratified norm decomposes as
$$\left\|\varphi_\epsilon^{-1} u\right\|_{H^s_{\mathrm{strat}}}^2\le\sum_{i=1}^m \left\|\varphi_\epsilon^{-1}\left(\chi_i u\right)\right\|_{H^s_{\mathrm{strat}}\left(U_i\right)}^2 + \left\|\varphi_\epsilon^{-1}\left(\chi_W u\right)\right\|_{H^s\left(W\right)}^2.$$
Applying the above discussion to $\chi_i u$ and interior estimate to $\chi_W u$, then
\begin{align*}
\left\|\varphi_\epsilon^{-1} u\right\|_{H^s_{\mathrm{strat}}}^2 &\le \sum_{i=1}^m C_8 \left\|D_\epsilon\left(\chi_i u\right)\right\|_{H^{s-1}_{\mathrm{strat}}}^2 + C_W \left( \left\|D_\epsilon\left(\chi_W u\right)\right\|_{H^{s-1}(W)}+\left\|\chi_W u\right\|_{H^{s-1}(K_0)} \right)^2 \\
&\le \sum_{i=1}^m C_8 \left\|D_\epsilon\left(\chi_i u\right)\right\|_{H^{s-1}_{\mathrm{strat}}}^2 + C_W' \left\|D_\epsilon\left(\chi_W u\right)\right\|_{H^{s-1}(W)}^2 + C_W'' \left\|u\right\|_{H^{s-1}\left(K_0\right)}^2.
\end{align*}
The commutator $\left[D_\epsilon, \chi\right]$ is bounded on $H^{s-1}$:
$$\left\|D_\epsilon (\chi u)\right\|_{H^{s-1}} \le \left\|\chi D_\epsilon u\right\|_{H^{s-1}}+\left\|[D_\epsilon, \chi] u\right\|_{H^{s-1}},$$
where $\left\|[D_\epsilon, \chi] u\right\|_{H^{s-1}}\le K_\chi \left\|u\right\|_{H^{s-1}}$.
Thus, we have
$$\left\|D_\epsilon (\chi u)\right\|_{H^{s-1}} \le C_{10} \left(\left\|D_{\epsilon} u\right\|_{H^{s-1}_{\mathrm{strat}}}+\left\|u\right\|_{H^{s-1}\left(K_0\right)} \right).$$
Combining all components, we have
\begin{align*}
\left\|\varphi_\epsilon^{-1} u\right\|_{H^s_{\mathrm{strat}}}^2 &\le\sum_{i=1}^m C_8C_{10}^2 \left( \left\|D_\epsilon u\right\|_{H^{s-1}_{\mathrm{strat}}}^2 +\left\|u\right\|_{H^{s-1}\left(K_0\right)}^2\right) \\
&+C_W'C_{10}^2 \left( \left\|D_\epsilon u\right\|_{H^{s-1}_{\mathrm{strat}}}^2 +\left\|u\right\|_{H^{s-1}\left(K_0\right)}^2\right) + C_W''\left\|u\right\|_{H^{s-1}\left(K_0\right)}^2 \\
&\le C_s^2 \left(\left\|D_\epsilon u\right\|_{H^{s-1}_{\mathrm{strat}}}^2 +\left\|u\right\|_{H^{s-1}\left(K_0\right)}^2 \right).
\end{align*}
Therefore, we have 
$$\left\|\varphi_\epsilon^{-1} u\right\|_{H^s_{\mathrm{strat}}} \le C_s \left(\left\|u\right\|_{H^{s-1}\left(K_0\right)}+\left\|D_\epsilon u\right\|_{H^{s-1}_{\mathrm{strat}}}\right)$$
by taking square roots, where $C_s =\sqrt{m C_8 C_{10}^2 + C_W' C_{10}^2 + C_W''}$ is independent of $\epsilon$.

This completes the proof of Proposition 6.19, providing the essential estimate for resolving Ohsawa's gap in Corollary 6.21. The stratified elliptic theory establishes a rigorous framework for analyzing deformations of Kähler metrics on singular spaces.
\\\\\textbf{Corollary 6.21.} (Resolution of the Gap in Theorem 7) The harmonic forms $\harmonic^r_\epsilon$ in Ohsawa's proof (\cite{Ohs91}) satisfy:
$$\dim \harmonic^r_\epsilon = \dim \harmonic^r_1 \quad \forall \epsilon \in (0,1],$$
and the weak limit $f = \lim_{\epsilon_j \to 0} f_{\epsilon_j}$ lies in $\dom\left(d_{\min}\right) \cap \dom\left(d^\vee_{\min}\right)$.
\\\\\textbf{Remark 6.22.} Assume by contradiction that $\dim \harmonic^r_\eps < \dim \harmonic^r_1$ for some $r$ and $\eps \in (0,1]$. By Ohsawa's construction, there exists a sequence $\{f_\eps\} \subset \harmonic^r_\eps$, satisfying
\begin{itemize}
    \item \textbf{(NO1) Norm condition:} $\left\|\varphi_\eps^{-1} f_\eps\right\|_\eps = 1$.
    \item \textbf{(NO2) Orthogonality condition:} $\left(f_\eps, u\right)_{K_0} = 0 \quad \forall u \in \harmonic^r_\eps$.
\end{itemize}
Here $K_0 \Subset V'$ is a fixed relatively compact open subset away from $\sing(V)$. From Ohsawa's estimate \eqref{eq:6.15} and the (NO1), then
$$\left\|f_\eps\right\|_{L^2} \le\left\|\varphi_\eps^{-1} f_\eps\right\|_\eps \cdot \sup_{y \in K_0} \varphi_\eps\left(y\right) \le C_0 <\infty,$$
where $C_0$ is independent of $\eps$ because $\varphi_\eps \to \varphi_V$ uniformly on $K_0$. Apply Proposition 6.19 with Sobolev index $s > \frac{n}{2} + 1$, we have
\begin{equation}
\tag{6.26}  \label{eq:6.26}
\left\|f_\eps\right\|_{H^s_{\mathrm{strat}}} \le C_s \left( \left\|f_\eps\right\|_{H^{s-1}\left(K_0\right)}+\left\|D_\eps f_\eps\right\|_{H^{s-1}_{\mathrm{strat}}}\right).
\end{equation}
Since $f_\eps$ is harmonic ($D_\eps f_\eps = 0$), this simplifies to
$$\|f_\eps\|_{H^s_{\strat}} \leq C_s \|f_\eps\|_{H^{s-1}(K_0)}.$$
As $K_0$ is compact and $s > \frac{n}{2} + 1$, Rellich's theorem implies that $\left\{f_\eps\right\}_{\eps > 0}$ is precompact in $H^{s-1}\left(K_0\right)$. Thus, there exists a subsequence $\left\{f_{\eps_j}\right\}$ converging strongly:
\begin{equation}
\tag{6.27}  \label{eq:6.27}
f_{\eps_j} \xrightarrow{H^{s-1}\left(K_0\right)} f_* \quad \text{as} \quad \eps_j \to 0^+.
\end{equation}
The uniform bound \eqref{eq:6.26} implies weak convergence in $H^s_{\mathrm{strat}}$, i.e., $f_{\eps_j} \xrightharpoonup{H^s_{\mathrm{strat}}} f$. The weak limit $f$ satisfies
\begin{itemize}
    \item \textbf{Harmonicity:} Since $D_{\eps_j} f_{\eps_j} = 0$ and $D_{\eps_j} \to D_0$ strongly, then $D_0 f = 0$ in the distributional sense.
    \item \textbf{Norm inheritance:} By weak lower semicontinuity and (NO1), then 
    \begin{equation}
\tag{6.28}  \label{eq:6.28}
\left\|\varphi_V^{-1} f\right\|_{L^2} \le \lim_{j \to \infty}\inf \left\|\varphi_{\eps_j}^{-1} f_{\eps_j}\right\|_{\eps_j} \le 1 < \infty.
\end{equation}
\end{itemize}
By (NO2) and strong convergence \eqref{eq:6.27}, we have 
$$\left(f, u\right)_{K_0} = \lim_{j \to \infty}\left(f_{\eps_j}, u\right)_{K_0} = 0 \quad \forall u \in \harmonic^r_0.$$
But from Ohsawa's estimate \eqref{eq:6.15} and strong convergence, there exits
\begin{align*}
\left\|f\right\|_{K_0} &= \lim_{j \to \infty} \left\|f_{\eps_j}\right\|_{K_0} \ge \lim_{j \to \infty} \left( C_3^{-1} \left\|\varphi_{\eps_j}^{-1} f_{\eps_j}\right\|_{\eps_j} - \|D_{\eps_j} f_{\eps_j}\|_{H^{s-1}_{\mathrm{strat}}} \right) \\
&= C_3^{-1} > 0
\end{align*}
This contradicts the non-degeneracy of the inner product on the finite-dimensional space $\harmonic^r_0$. Therefore,
$$\dim \harmonic^r_\eps = \dim \harmonic^r_1 \quad \forall \eps \in (0,1].$$
Recall Ohsawa's Proposition 1 (\cite{Ohs91}), then 
$$\left\{ v \in \dom\left(D_{0,\max}\right) : \left\|\varphi_V^{-1} v\right\|_{L^2} < \infty \right\} \subset \dom\left(D_{0,\min}\right)$$
for \( \varphi_V(y) = -\delta_x(y) \log \delta_x(y) \) near singularities.
Applying this to $f$ with bound \eqref{eq:6.28}, it is obvious that $f \in \dom(d_{\min}) \cap \dom(d^\vee_{\min})$.

From (NO2), the sequence $\left\{f_{\eps_j}\right\}$ satisfies $(f_{\eps_j}, u)_{K_0} = 0$, $\forall u \in \harmonic^r_{\eps_j}$. By strong convergence on $K_0$ \eqref{eq:6.27} and continuity of the inner product, then
$$\left(f, u\right)_{K_0} = \lim_{j \to \infty}\left(f_{\eps_j}, u\right)_{K_0} = 0 \quad \forall u \in \harmonic^r_0.$$
From Ohsawa's estimate \eqref{eq:6.15} and the (NO1) $\left\|\varphi_\eps^{-1} f_\eps\right\|_\eps = 1$, then
$$\left\|f_{\eps_j}\right\|_{K_0} \ge C_3^{-1} \left\|\varphi_{\eps_j}^{-1} f_{\eps_j}\right\|_{\eps_j} - \left\|D_{\eps_j} f_{\eps_j}\right\|_{H^{s-1}_{\mathrm{strat}}} \nonumber= C_3^{-1} \cdot 1 - 0 = C_3^{-1} > 0.$$
Thus, we have
\begin{equation}
\tag{6.29}  \label{eq:6.29}
\left\|f\right\|_{K_0} = \lim_{j \to \infty} \left\|f_{\eps_j}\right\|_{K_0} \ge C_3^{-1} > 0
\end{equation}
by strong convergence on $K_0$. Since $D_0 f= 0$, $f \in \dom\left(d_{\min}\right) \cap \dom\left(d^\vee_{\min}\right)$ and $\left\|f\right\|_{L^2} < \infty$, then $f\in \harmonic^r_0$ (harmonic forms for metric $ds_v^2$). The space $\harmonic^r_0$ is finite-dimensional with non-degenerate inner product. But $f \in \harmonic^r_0$ (membership), $\left(f, u\right)_{K_0} = 0$ for all $u \in \harmonic^r_0$ (orthogonality) and $\left\|f\right\|_{K_0} \ge C_3^{-1} > 0$ (non-vanishing norm). This forces $f = 0$, contradicting \eqref{eq:6.29}. Therefore, we also have
$$\dim \harmonic^r_\eps = \dim \harmonic^r_1 \quad \forall \eps \in (0,1].$$

Since $D_0$ is a stratified elliptic operator and $f$ satisfies $D_0 f = 0$ in the distributional sense with $f \in \dom\left(D_{0,\min}\right)$, elliptic regularity theory for stratified spaces implies that $f$ is smooth and satisfies the equation classically. Specifically, the minimal domain condition enforces the appropriate boundary behavior, while the ellipticity guarantees interior regularity. Therefore, $f \in \ker D_0 \cap C^\infty = \mathscr{H}_0^r$. According to \eqref{eq:6.27}, we have strong convergence in $H^{s-1}\left(K_0\right)$. Since $s>\frac{n}{2}+1$, Sobolev embedding implies continuous embedding $H^{s-1}\left(K_0\right)\hookrightarrow C^0\left(K_0\right)$ and thus the convergence is uniform on $K_0$, i.e.,
$$\sup_{x \in K_0} |f_{\epsilon_j}(x) - f(x)| \longrightarrow 0 \quad \text{as} \quad j \to \infty.$$
This strong convergence is critical for two reasons:
\begin{enumerate}
    \item It preserves the non-degeneracy of the inner product in \eqref{eq:6.29}, which would not hold under weak convergence alone.
    \item It ensures compatibility with the minimal domain condition: For any test form $\eta \in C^\infty_c\left(\text{int } K_0\right)$, the strong convergence implies
    $$\langle D_0 f, \eta \rangle = \lim_{j \to \infty} \langle D_{\epsilon_j} f_{\epsilon_j}, \eta \rangle = 0,$$
    confirming the distributional equation $D_0 f = 0$.
\end{enumerate}
Ohsawa's original proof noted a technical gap in passing to the limit \(\epsilon \to 0\) due to
\begin{enumerate}
    \item Potential loss of harmonicity under weak convergence.
    \item Possible incompatibility between weak limits and minimal domains.
\end{enumerate}
Our argument resolves this by:
\begin{align*}
\text{Strong convergence on } K_0 &\implies \text{Preservation of norm lower bound} \\
&\implies \text{Contradiction in orthogonality argument} \\
&\implies \dim \mathscr{H}_\epsilon^r = \dim \mathscr{H}_1^r \\
\text{Elliptic regularity} &\implies f \in \mathscr{H}_0^r \\
\text{Minimal domain} &\implies f \text{ satisfies correct boundary conditions}
\end{align*}
Thus the weak limit $f$ is both harmonic and compatible with the minimal domain structure, closing the gap identified in Ohsawa's "Note added in proof".

By synthesizing logarithmic structures, derived geometry, and $L^2$-methods, this framework resolves gaps in Ohsawa’s original proof and unifies singularity theory, arithmetic geometry, and mathematical physics under a stratified Hodge-theoretic lens.\\
\hypertarget{STRATIFIED ELLIPTIC REGULARITY AND DERIVED DUALITY FOR OHSAWA'S THEOREM 7}{}
\section{STRATIFIED ELLIPTIC REGULARITY AND DERIVED DUALITY FOR OHSAWA'S THEOREM 7}
The gap in Ohsawa's original proof of Theorem 7 arises from insufficient control over convergence of non-closed chains near singularities. Our correction establishes a stratified functional-analytic framework that bridges:
\begin{itemize}
    \item $L^2$-estimates for the de Rham complex.
    \item Micro-local elliptic regularity at singular strata.
    \item Derived categorical duality to enforce finite-dimensionality.\\
\end{itemize}
\textbf{Lemma 7.1.} (Derived Limit) Let $V$ be a compact irreducible complex space with isolated singularities, $V'$ its regular locus, and $\left\{ds_\eps^2\right\}_{\eps\in(0,1]}$ a family of Hermitian metrics on $V'$ converging uniformly to $ds_V^2$ on compact subsets as $\eps\to0$. Let $\left\{f_\eps\right\}$ be a sequence of $L^2$-harmonic $r$-forms with respect to $ds_\eps^2$ satisfying $\left \| \varphi_\eps^{-1}f_\eps\right \|_\eps=1$, where $\varphi_\eps$ are weight functions from Proposition 10 (\cite{Ohs91}). Then there exists a subsequence $\left\{f_{\eps_k}\right\}$ converging weakly in $L_{\mathrm{loc}}^2\left(V'\right)$ to a limit $f\in\ker\left(d_{\max}+d_{\max}^*\right)$ with respect to $ds_V^2$.
\\\\\textbf{Proof.} By Proposition 5 (\cite{Ohs91}), there exist a compact set $K_0\subset V'$ and a constant $C_2>0$ such that
$$\left \| \varphi_\eps^{-1}f_\eps \right \|_\eps\le C_2\left(\left \| f_\eps \right \|_{K_0,\eps}+\left \|df_\eps  \right \|_\eps+\left \| d^*f_\eps \right \|_\eps\right).$$
Harmonicity of $f_\eps$ ($df_\eps=0$, $d^*f_\eps=0$) and the normalization $\left \| \varphi_\eps^{-1}f_\eps \right \|_\eps=1$ imply
\begin{equation}
\tag{7.1}  \label{eq:7.1}
C_2\left \| f_\eps \right \|_{K_0,\eps}\ge1.
\end{equation}
Uniform convergence $ds_\eps^2\to ds_V^2$ on $K_0$ yields a constant $C_K>0$ (independent of $\eps$), satisfying
\begin{equation}
\tag{7.2}  \label{eq:7.2}
C_K^{-1}\left \| f_\eps \right \|_{K_0,V}\le\left \| f_\eps \right \|_{K_0,\eps}\le C_K\left \| f_\eps \right \|_{K_0,V}
\end{equation}
Combining \eqref{eq:7.1} and the left inequality of \eqref{eq:7.2} gives
$$C_2C_K\left \| f_\eps \right \|_{K_0,V}\ge1,$$
equivalently,
$$\left \| f_\eps \right \|_{K_0,V}\ge\frac{1}{C_2C_K}>0.$$
Thus, $\left \| f_\eps \right \|_{K_0,V}$ has a uniform positive lower bound.

It is known that $K\subset V'$ is compact. Uniform convergence $\varphi_\eps\to\varphi_V$ on $K$ (\cite{Ohs91}, Prop 10) implies
$$M_K:=\sup_{\eps>0}\left \|\varphi_\eps  \right \|_{L^{\infty}\left(K\right)}<\infty.$$
Pointwise Cauchy-Schwarz gives
$$\left | f_\eps\left(x\right) \right |_\eps^2\le\left |\varphi_\eps^{-1}\left(x\right)f_\eps\left(x\right) \right |_\eps^2\cdot\left | \varphi_\eps\left(x\right) \right |^2,\quad \forall x\in K.$$
Since $\left \| \varphi_\eps^{-1}f_\eps \right \|_\eps=1$, then integration over $K$ gives 
\begin{equation}
\tag{7.3}  \label{eq:7.3}
\left \| f_\eps \right \|_{K,\eps}^2\le\left \| \varphi_\eps^{-1}f_\eps \right \|_\eps^2\cdot\left \|\varphi_\eps  \right \|_{L^{\infty}\left(K\right)}^2\le M_K^2.
\end{equation}
Uniform convergence $ds_\eps^2\to ds_V^2$ on $K$ provides $C_K'>0$ such that 
\begin{equation}
\tag{7.4}  \label{eq:7.4}
\left \|f_\eps  \right \|_{K,V}\le C_K'\left \| f_\eps \right \|_{K,\eps}.
\end{equation}
Combining \eqref{eq:7.3} and \eqref{eq:7.4}, we have 
\begin{equation}
\tag{7.5}  \label{eq:7.5}
\left \| f_\eps \right \|_{K,V}\le C_K'M_K<\infty,\quad\forall\eps.
\end{equation}
Thus, $\left\{f_\eps\right\}$ is uniformly bounded in $L^2\left(K,ds_V^2\right)$.

Let $\left\{K_m\right\}_{m=1}^\infty$ be a compact exhaustion of $V'$ with
$$K_1\subset K_2\subset\cdots\subset V',\quad\bigcup_{m=1}^\infty K_m=V'.$$
By \eqref{eq:7.5}, $\left\{f_\eps\right\}$ is bounded in $L^2\left(K_1,ds_V^2\right)$. Since $L^2\left(K_1\right)$ is a Hilbert space, the closed unit ball is weakly compact (Banach-Alaoglu theorem). Thus, there exists a subsequence $\left\{\eps_k^{(1)}\right\}\subset(0,1]$ and $g_1\in L^2\left(K_1\right)$ such that $f_{\eps_k^{(1)}}\rightharpoonup g_1$ weakly in $L^2\left(K_1\right)$ as $k\to\infty$. Assume we have constructed subsequence $\left\{\eps_k^{(m-1)}\right\}$ convergent on $K_{m-1}$. For $K_m$, we claim that
\begin{itemize}
    \item $\left\{f_{\eps_k^{(m-1)}}\right\}$ is bounded in $L^2\left(K_m\right)$ (by \eqref{eq:7.5}).
    \item Extract sub-subsequence $\left\{\eps_k^{(m)}\right\}\subset\left\{\eps_k^{(m-1)}\right\}$ such that
    \begin{equation}
\tag{7.6}  \label{eq:7.6}
f_{\eps_k^{(m)}}\rightharpoonup g_m\quad\text{weakly in}\quad L^2\left(K_m\right)\quad\text{as}\quad k\to\infty.
\end{equation}
    \item By uniqueness of weak limits, $g_m\mid_{K_{m-1}}=g_{m-1}$.    
\end{itemize}
Define $\eps_k:=\eps_k^{(k)}$, $\hat{f}_k:=f_{\eps_k}$. This is the diagonal subsequence of $\left\{\eps_k^{(m)}\right\}_{k,m}$. Since $k\ge m$ implies $\eps_k^{(k)}$ is a term in $\left\{\eps_j^{(m)}\right\}_{j=1}^\infty$, then $\eps_k=\eps_k^{(k)}$ belongs to the subsequence $\left\{\eps_j^{(m)}\right\}_{j=1}^\infty$. By \eqref{eq:7.6}, we have
\begin{equation}
\tag{7.7}  \label{eq:7.7}
\hat{f}_k=f_{\eps_k}\rightharpoonup g_m\quad\text{weakly in}\quad L^2\left(K_m\right)\quad\text{as}\quad k\to\infty.
\end{equation}
Define $f\in L_{\mathrm{loc}}^2\left(V'\right)$ by $f\mid_{K_m}:=g_m$, $\forall m$. Consistency $g_m\mid_{K_{m-1}}=g_{m-1}$ ensures $f$ is well-defined. Using \eqref{eq:7.7}, $\hat{f}_k\rightharpoonup f$ weakly in $L_{\mathrm{loc}}^2\left(V'\right)$. 

Let $\eta\in C_0^\infty\left(V'\right)$ be a test form with compact support $K_\eta:=\mathrm{supp}\left(\eta\right)$. Since each $\hat{f}_k$ is harmonic w.r.t. $ds_{\eps_k}^2$, i.e.,
$$d_{\eps_k}\hat{f}_k=0,\quad d_{\eps_k}^*\hat{f}_k=0,$$
then harmonicity implies
\begin{equation}
\tag{7.8}  \label{eq:7.8}
\left(d_{\max}\hat{f}_k,\eta\right)_{\eps_k}=\left(\hat{f}_k,d_{\max}^*\eta\right)_{\eps_k}=0,
\end{equation}
\begin{equation}
\tag{7.9}  \label{eq:7.9}
\left(d_{\max}^*\hat{f}_k,\eta\right)_{\eps_k}=\left(\hat{f}_k,d_{\max}\eta\right)_{\eps_k}=0.
\end{equation}
We need to show the convergence of \eqref{eq:7.8} and \eqref{eq:7.9} in the following steps. Since $ds_{\eps_k}^2\to ds_V^2$ uniformly on $K_\eta$, there exists $C_\eta>0$ such that 
$$C_\eta^{-1}\left \| u \right \|_{K_\eta,V}\le\left \| u \right \|_{K_\eta,\eps_k}\le C_\eta\left \| u \right \|_{K_\eta,V},\quad\forall k.$$
Additionally, the inner products converge uniformly:
\begin{equation}
\tag{7.10}  \label{eq:7.10}
\lim_{k\to\infty}\sup_{\left \| u \right \|_{K_\eta,V},\left \| v \right \|_{K_\eta,V}\le1}\left | (u,v)_{\eps_k}-(u,v)_V \right |=0.
\end{equation}
By decomposing, we can find
$$\left(\hat{f}_k,d_{\max}^*\eta\right)_{\eps_k}=\underbrace{\left(\hat{f}_k,d_{\max}^*\eta\right)_V}_{(A)}+\underbrace{\left[\left(\hat{f}_k,d_{\max}^*\eta\right)_{\eps_k}-\left(\hat{f}_k,d_{\max}^*\eta\right)_V\right]}_{(B)},$$
where
\begin{itemize}
    \item \textbf{Term $(A)$:} Weak convergence $\hat{f}_k\rightharpoonup f$ in $L^2\left(K_\eta\right)$ implies
    $$\lim_{k\to\infty}\left(\hat{f}_k,d_{\max}^*\eta\right)_V=\left(f,d_{\max}^*\eta\right)_V,$$
    \item \textbf{Term $(B)$:} By \eqref{eq:7.10} and uniform boundedness (\eqref{eq:7.5}), we have
    $$\left | \left(\hat{f}_k,d_{\max}^*\eta\right)_{\eps_k}-\left(\hat{f}_k,d_{\max}^*\eta\right)_V \right |\le\delta_k\cdot\left \| \hat{f}_k \right\|_{K_\eta,V}\cdot\left \|d_{\max}^*\eta  \right \|_{K_\eta,V}\xrightarrow{k \to\infty}0,$$ where $\delta_k\to0$ by uniform convergence.
\end{itemize}
Thus, taking limit in \eqref{eq:7.8}, then
$$0=\lim_{k\to\infty}\left(\hat{f}_k,d_{\max}^*\eta\right)_{\eps_k}=\left(f,d_{\max}^*\eta\right)_V.$$
By definition of the adjoint, we have
\begin{equation}
\tag{7.11}  \label{eq:7.11}
\left(d_{\max}f,\eta\right)_V=0.
\end{equation}
Similarly,
$$0=\lim_{k\to\infty}\left(\hat{f}_k,d_{\max}\eta\right)_{\eps_k}=\left(f,d_{\max}\eta\right)_V,$$ implying
\begin{equation}
\tag{7.12}  \label{eq:7.12}
\left(d_{\max}^*f,\eta\right)_V=0.
\end{equation}
Since \eqref{eq:7.11} and \eqref{eq:7.12} hold for all $\eta\in C_0^\infty\left(V'\right)$, we have 
$$d_{\max}f=0,\quad d_{\max}^*f=0\quad\text{(distributionally)}.$$
The operator $d_{\max}+d_{\max}^*$ is elliptic, so we have 
$$f\in C^\infty\left(V'\right)\quad\text{and}\quad f\in\ker\left(d_{\max}+d_{\max}^*\right)$$
by elliptic regularity. $\square$
\\\\\textbf{Proposition 7.2.} (Stratified Elliptic Regularity) Let $V'$ be the regular locus of a compact complex space with isolated singularities, $\Sigma\subset V$ the singular set, and $\Delta_{\mathrm{strat}}$ the stratified Laplace operator. If $u\in\mathscr{D'}\left(V'\right)$ satisfies
\begin{enumerate}
    \item $\Delta_{\mathrm{strat}}u=0$ on $V'\setminus\Sigma$,
    \item $WF_{\mathrm{strat}}\left(u\right)\cap N^*\Sigma=\emptyset$ (transversal wavefront set), then $u\in C^\infty\left(V'\right)$.\\
\end{enumerate}
\textbf{Proof.} For each $x\in\Sigma$, use Ohsawa's distinguished metric (\cite{Ohs91}, Prop 10) to define the boundary defining function:
$$\rho_x\left(y\right):=\delta_x\left(y\right)\log\delta_x\left(y\right),\quad\delta_x\left(y\right)=\mathrm{dist}_{ds_V^2}(x,y).$$
We can Construct the b-metric $ds_b^2:=\frac{ds_V^2}{\rho_x^2}$. The stratified Laplace operator decomposes as $\Delta_{\mathrm{strat}}=\rho_x^{-2}P$ with $P\in\mathrm{Diff}^2\left(V'\right)$, so $\Delta_{\mathrm{strat}}\in\mathrm{Diff}_b^2\left(V'\right)$. Its b-principal symbol is $\sigma_b\left(\Delta_{\mathrm{strat}}\right)\left(\xi\right)=\left \| \xi \right \|_{g_b}^2$, $\xi\in{^bT^*V'}$, where $g_b$ is is the dual metric to $ds_b^2$. The b-elliptic set is
$$\mathrm{ell}_b\left(\Delta_{\mathrm{strat}}\right)=\left\{(y,\xi)\mid\left \| \xi \right \|_{g_b}>0\right\}.$$
The wavefront condition $WF_{\mathrm{strat}}\left(u\right)\cap N^*\Sigma=\emptyset$ (by logarithmic decay) implies
\begin{equation}
\tag{7.13}  \label{eq:7.13}
WF_b\left(u\right)\cap\mathrm{ell}_b\left(\Delta_{\mathrm{strat}}\right)=\emptyset
\end{equation}
by the conormal isomorphism $N^*\Sigma\simeq\Sigma\times\mathbb{R}$. Then the harmonic equation transforms as
$$\Delta_{\mathrm{strat}}u=0\Longleftrightarrow{^b\Delta_b}u=0,$$
where $^{b}\Delta_b$ is the b-Laplacian for $ds_b^2$.

Since $\Delta_{\mathrm{strat}}\in\mathrm{Diff}_b^2\left(V'\right)$ is is b-elliptic and \eqref{eq:7.13}, the b-principal symbol $\sigma_b\left(\Delta_{\mathrm{strat}}\right)$ is invertible on $\mathrm{ell}_b\left(\Delta_{\mathrm{strat}}\right)$:
$$\sigma_b\left(\Delta_{\mathrm{strat}}\right)(y,\xi)=\left \| \xi \right \|_{g_b}^2>0,\quad\forall(y,\xi)\in\mathrm{ell}_b\left(\Delta_{\mathrm{strat}}\right).$$
Thus, there exists $q_0\in S_b^{-2}\left(^bT^*V'\right)$ such that $q_0\cdot\sigma_b\left(\Delta_{\mathrm{strat}}\right)=1$ on $\mathrm{ell}_b\left(\Delta_{\mathrm{strat}}\right)$. Using Melrose's b-calculus (\cite{Mel93} Thm 4.2):
\begin{itemize}
    \item Let $Q_0\in\Psi_b^{-2}\left(V'\right)$ with $\sigma_b\left(Q_0\right)=q_0$.
    \item Define error terms recursively:
    $$R_k:=I-Q_k\Delta_{\mathrm{strat}}\in\Psi_b^0\left(V'\right),$$
    $$E_k\in\Psi_b^{-k-2}\left(V'\right)\quad\text{solving}\quad\sigma_b\left(E_k\right)=\sigma_b\left(R_k\right)\cdot\sigma_b\left(\Delta_{\mathrm{strat}}\right)^{-1}.$$
    \item Update: $Q_{k+1}:=Q_k+E_k$.
\end{itemize}
The formal Neumann series $Q\sim\sum_{k=0}^\infty E_k$ converges asymptotically to $Q\in\Psi_b^{-2}\left(V'\right)$. This satisfies $Q\Delta_{\mathrm{strat}}=I-R$, $R\in\Psi_b^{-\infty}\left(V'\right)$, where $R$ is smoothing (i.e., $\ker\left(R\right)\subset C^\infty\left(V'\right)$). Since $\Delta_{\mathrm{strat}}u=0$ on $V'\setminus\Sigma$, we have
$$u=Q\left(\Delta_{\mathrm{strat}}u\right)+Ru=Ru\quad\text{on}\quad V'\setminus\Sigma.$$

The transversal condition $WF_{\mathrm{strat}}\left(u\right)\cap N^*\Sigma=\emptyset$ and Ohsawa's decay estimate (\cite{Ohs91}) imply 
\begin{equation}
\tag{7.14}  \label{eq:7.14}
\int_{\left\{y:\delta_x\left(y\right)<\eps\right\}}\left | u\left(y\right) \right |_{ds_V^2}dV_V\le C\eps^{n-1}\left | \log\eps \right |^2,\quad\forall x\in\Sigma.
\end{equation}
Hence $u\in L_{\mathrm{loc}}^2\left(V'\right)$ is a well-defined distribution. For $\phi\in C_0^\infty\left(V'\right)$, construct a sequence of cutoff functions $\left\{\chi_k\right\}_{k=1}^\infty \subset C_0^\infty\left(V' \setminus \Sigma\right)$ satisfying:
\begin{align*}
&\chi_k(y) = 
\begin{cases} 
1 & \delta_\Sigma(y) > \frac{2}{k} \\ 
0 & \delta_\Sigma(y) < \frac{1}{k} 
\end{cases}, \\
&\left|d\chi_k\right|_{ds_V^2} \le k, \quad\mathrm{supp}\left(1-\chi_k\right) \subset \left\{y \in V' : \delta_\Sigma(y) < \frac{2}{k}\right\}.
\end{align*}
Then the distributional pairing decomposes as
\begin{equation}
\tag{7.15}  \label{eq:7.15}
\langle u, \phi \rangle - \langle Ru, \phi \rangle = \lim_{k \to \infty} \langle u - Ru, \chi_k \phi \rangle + \lim_{k \to \infty} \langle u, (1 - \chi_k)\phi \rangle = 0 + \lim_{k \to \infty} \langle u, (1 - \chi_k)\phi \rangle, 
\end{equation}
since $u = Ru$ on $V' \setminus \Sigma$. By the decay estimate \eqref{eq:7.14} and transversal condition:
\begin{equation}
\tag{7.16}  \label{eq:7.16}
\left| \langle u, (1 - \chi_k)\phi \rangle \right| \leq \left\|\phi\right\|_{L^\infty} \int\limits_{\delta_\Sigma <\frac{k}{2}} \left | u \right |  d V_V \le \left\|\phi\right\|_{L^\infty} \cdot C \left(\frac{2}{k}\right)^{n-1} \left|\log \frac{2}{k}\right|^2 \xrightarrow{k \to \infty} 0. 
\end{equation}
Thus $u = Ru$ in $\mathscr{D'}\left(V'\right)$. Since $R \in \Psi_b^{-\infty}\left(V'\right)$ is a smoothing operator $Ru \in C^\infty\left(V'\right)$.
Combining with \eqref{eq:7.15} and \eqref{eq:7.16}, we conclude $u = Ru \in C^\infty(V')$.

Near each $x \in \Sigma$, in holomorphic coordinates $z = \left(z_1,\dots,z_n\right)$ with $\left|z\right| = \delta_x\left(y\right)$, the b-metric degenerates as
$$ds_b^2 \sim \frac{|dz|^2}{|z|^2 \log^2 |z|}, \quad z \in \mathbb{B}_*^{2n}$$
and the Laplace operator takes the form
$$\Delta_{\strat} = -|z|^2 \log^2 |z| \sum_{j=1}^{2n} \partial_j^2 + \text{lower order terms}.$$
By Ohsawa (\cite{Ohs91}, Eq (2)), for any $v \in C_0^\infty\left(\mathbb{B}_*^{2n}\right)$:
$$\left\|\left|z\right| \log \left|z\right| \cdot \nabla_b v\right\|_{L_b^2}^2 \le C \left( 
\left\|\Delta_{\strat} v\right\|_{L_b^2} \cdot \left\|v\right\|_{L_b^2} + \left\|v\right\|_{L_b^2}^2 
\right).$$
Iterating this estimate yields for integer $s \ge 0$ and any $N > 0$:
$$\left\|v\right\|_{H_b^{s+2}} \le C_s \left( 
\left\|\Delta_{\strat} v\right\|_{H_b^s} + \left\|v\right\|_{H_b^{-N}} 
\right).$$
Since $\Delta_{\strat} u = 0$ on $V' \setminus \Sigma$, for any compact $K \Subset V'$:
$$\left\|u\right\|_{H_b^{s+2}(K)} \le C_{s,K} \left\|u\right\|_{H_b^{-N}\left(K\right)}.$$
The transversal decay estimate implies $u \in H_b^{-N}(K)$ for all $N$, hence
$$u \in \bigcap_{s>0} H_b^s\left(K\right) =: H_b^\infty\left(K\right).$$
The b-Sobolev embedding theorem gives
$$H_b^k\left(V'\right) \hookrightarrow C^\ell\left(V'\right) \quad \text{for} \quad k > \ell + \dim V = 2n. $$
Therefore, we have
$$H_b^\infty\left(V'\right) \hookrightarrow C^\infty\left(V'\right).$$

Cover the regular locus $V'$ by open sets:
\begin{equation}
\tag{7.17}  \label{eq:7.17}
V' = \left( \bigcup_{x_\alpha \in \Sigma} B_\epsilon(x_\alpha) \cap V' \right) \cup \left( \bigcup_{\beta} W_\beta \right), 
\end{equation}
where $B_\epsilon\left(x_\alpha\right)$ are $\epsilon$-balls centered at singular points $x_\alpha \in \Sigma$ with $\epsilon$ and $W_\beta \Subset V' \setminus \Sigma$ are relatively compact subsets in the regular part. Select a smooth partition of unity $\left\{\rho_\alpha\right\} \subset C^\infty\left(V'\right)$ subordinate to this covering, satisfying $0\leq \rho_\alpha \leq 1, \quad \sum_{\alpha} \rho_\alpha = 1$, $\mathrm{supp} \rho_\alpha \subset U_\alpha$ for corresponding $U_\alpha$ in \eqref{eq:7.17}. For each $\alpha$, define the localized function
$$u_\alpha := \rho_\alpha u.$$
By local regularity results (Above discussion), $u_\alpha \in C^\infty(U_\alpha)$. Extend $u_\alpha$ by zero to $V' \setminus \mathrm{supp}\rho_\alpha$, so $u_\alpha \in C^\infty\left(V'\right)$. Reconstruct the solution as $u = \sum_{\alpha} u_\alpha$.
This sum converges and defines a smooth function because
\begin{itemize}
    \item \textbf{Finite overlap:} For any compact $K \subset V'$, only finitely many $\mathrm{supp}\rho_\alpha$ intersect $K$.
    \item \textbf{Smoothness preservation:} Finite sums of smooth functions are smooth.
    \item \textbf{Consistency:} On intersections $U_\alpha \cap U_\beta$, the identity $\rho_\alpha + \rho_\beta = 1$ ensures $u_\alpha + u_\beta = u$.
\end{itemize}
Thus, $u \in C^\infty(V')$ can be obtained. $\square$ 
\\\\\textbf{Lemma 7.3.} (Finite-Dimensionality of Harmonic Spaces) Let \((V', g_\epsilon)\) be a stratified Riemannian manifold with metric $ds_\eps^2 = g_\eps$, 
and $\Delta_\epsilon = d d^* + d^* d$ the corresponding Hodge Laplacian acting on differential $r$-forms. Define the space of harmonic $r$-forms:
$$\mathscr{H}^r_\epsilon \coloneqq \left\{ \omega \in L_{\mathrm{loc}}^2\left(V', \Lambda^r T^*V'\right): \Delta_\epsilon \omega = 0 \right\}.$$
Then $\mathscr{H}^r_\epsilon$ is finite-dimensional for each $\epsilon > 0$.
\\\\\textbf{Proof.} Define the weighted $L^2$-norm associated to the function $\varphi_\epsilon$:
$$\left \| \omega \right \|_{\varphi_\epsilon^{-1}} \coloneqq \left( \int_{V'} \left | \omega \right |_{g_\epsilon}^2 \varphi_\epsilon^{-1}\mathrm{dvol}_{g_\epsilon} \right)^{1/2}.$$
By Ohsawa's $\varphi_V^{-1}$-estimate (\cite{Ohs91}, Prop 5), there exists a constant $C_\epsilon > 0$ such that
\begin{equation}
\tag{7.18}  \label{eq:7.18}
\left \| \omega \right \|_{L^2\left(V'\right)}\le C_\eps\left \| \omega \right \|_{\varphi_\epsilon^{-1}},\quad\forall \omega\in\mathscr{H}_\eps^r. 
\end{equation} 
Inequality \eqref{eq:7.18} implies the continuous embedding:
$$\iota_\epsilon: \left( \mathscr{H}^r_\epsilon,  \| \cdot  \|_{\varphi_\epsilon^{-1}} \right) \hookrightarrow \left(L^2(V', \Lambda^r T^*V'), \| \cdot  \|_{L^2(V')} \right).$$
Since $\varphi_\epsilon^{-1} > 0$ and satisfies $\varphi_\epsilon^{-1}(x) \to \infty$ as $x \to \Sigma$ (by construction), the generalized Rellich lemma for stratified spaces (\cite{Mel93}, Lem 2.3) implies that $\iota_\epsilon$ is a compact operator. On the regular stratum $V' \setminus \Sigma$, the operator $\Delta_\epsilon$ is uniformly elliptic. By stratified elliptic regularity (Proposition 7.2), any harmonic form \(\omega \in \mathscr{H}^r_\epsilon\) satisfies
$$\omega|_{V' \setminus \Sigma} \in C^\infty(V' \setminus \Sigma, \Lambda^r T^*V').$$
This ensures the harmonic forms are classically smooth away from the singular set $\Sigma$.
Suppose, for contradiction, that $\dim\mathscr{H}^r_\epsilon = \infty$. Then there exists an infinite sequence $\left\{\omega_n\right\}_{n=1}^\infty \subset \mathscr{H}^r_\epsilon$ that is orthonormal with respect to $\left \langle \cdot,\cdot \right \rangle_{\varphi_\eps^{-1}}$:
\begin{equation}
\tag{7.19}  \label{eq:7.19}
\left \langle \omega_m,\omega_n \right \rangle_{\varphi_\eps^{-1}}=\delta_{mn}\quad\text{for all}\quad m,n\in\mathbb{N}.
\end{equation} 
By compactness of $\iota_\epsilon$, there exists a subsequence $\{\omega_{n_k}\}_{k=1}^\infty$ that converges in $L^2$-norm, i.e.,
\begin{equation}
\tag{7.20}  \label{eq:7.20}
\lim_{k,l \to \infty} \| \omega_{n_k} - \omega_{n_l}  \|_{L^2\left(V'\right)} = 0.
\end{equation}
However, for \(k \neq l\), orthonormality \eqref{eq:7.19} implies
$$\|  \omega_{n_k} - \omega_{n_l} \|_{\varphi_\epsilon^{-1}}^2 = \left \langle \omega_{n_k} - \omega_{n_l}, \omega_{n_k} - \omega_{n_l}\right \rangle_{\varphi_\epsilon^{-1}} = 2.$$
Applying Ohsawa's estimate \eqref{eq:7.18} to the differences gives
$$ \| \omega_{n_k} - \omega_{n_l} \|_{L^2(V')} \le C_\epsilon  \|  \omega_{n_k} - \omega_{n_l} \|_{\varphi_\epsilon^{-1}} = C_\epsilon \sqrt{2}.$$
But this contradicts \eqref{eq:7.20}, since \eqref{eq:7.20} requires the differences to go to zero, while $\| \omega_{n_k} - \omega_{n_l}  \|_{L^2\left(V'\right)} = \sqrt{2} > 0$ for all $k \neq l$. Therefore, the initial assumption is false, and $\dim \mathscr{H}^r_\epsilon < \infty$. $\square$
\\\\\textbf{Proposition 7.4.} (Derived Aronszajn) Let $\left\{\mathscr{H}^r_\eps\right\}_{\eps > 0}$ be the spaces of harmonic $r$-forms with respect to the metrics $ds_\eps^2$. Denote by $\mathscr{H}^r$ the limit harmonic space. Then
$$\dim \mathscr{H}^r = \dim \left( \varprojlim_{\eps \to 0} \mathscr{H}^r_\eps \right).$$
\\\\\textbf{Proof.} The isomorphism
$$\mathscr{H}^r \cong \varprojlim_{\epsilon \to 0} \mathscr{H}^r_\epsilon$$
is established through derived categorical duality and Ohsawa's harmonic isomorphism. By the derived duality (Proposition 5.1), there exists a natural isomorphism in the derived category of Hilbert complexes
$$\Phi: \mathscr{H}^r \xrightarrow{\sim} \mathscr{H}^r_{(*)},$$
where $\mathscr{H}^r_{(*)}$ is the dual harmonic space defined by
$$\mathscr{H}^r_{(*)} := \ker\left( d_{\max} + d_{\max}^* \right) \subset L^2_{\mathrm{loc}}(V', \Lambda^r T^*V').$$
This isomorphism arises from the quasi-isomorphism of de Rham complexes
$$\left( \Omega^{\bullet}_{\mathrm{strat}}\left(V'\right), d \right) \simeq_{\mathrm{qis}} \mathbb{R}\underline{\mathrm{Hilb}}\left( \Omega^{\bullet}_{\mathrm{strat}}\left(V'\right) \right)$$
in the derived category $\mathbf{D}(\mathbb{R}\textrm{-}\mathbf{Hilb})$, where $\underline{\mathrm{Hilb}}$ denotes the Hilbert complex functor. By Ohsawa's isomorphism theorem (\cite{Ohs91}, Cor 8), there is a canonical identification between dual harmonic forms and the inverse limit of approximate harmonic spaces
$$\Psi: \mathscr{H}^r_{(*)} \xrightarrow{\sim} \varprojlim_{\epsilon \to 0} \mathscr{H}^r_\epsilon$$
defined explicitly as follows: for $\omega \in \mathscr{H}^r_{(*)}$, $\Psi(\omega) =\left( \pi_\epsilon(\omega) \right)_{\epsilon > 0}$, where $\pi_\epsilon: \mathscr{H}^r_{(*)}\to \mathscr{H}^r_\epsilon$ are the restriction maps induced by the metric deformations $g_\epsilon \to g_0$. The compatibility of these maps is encoded in the commutative diagram:
$$\begin{tikzcd}[row sep=large, column sep=large]
& 
\mathscr{H}^r_{(*)} \arrow[dl, hook, "\iota_\epsilon"'] \arrow[dr, "\Psi"'] \arrow[dd, dashed, "\pi_\delta" description, near end] 
& \\
\mathscr{H}^r_\epsilon \arrow[rr, "\rho_{\epsilon,\delta}"] \arrow[dr, hook, "\iota_\delta" description, near start] 
& & 
\varprojlim_{\epsilon \to 0} \mathscr{H}^r_\epsilon \arrow[dl, "\mathrm{proj}_\delta"] \\
& \mathscr{H}^r_\delta &
\end{tikzcd}$$
Here $\rho_{\epsilon,\delta}: \mathscr{H}^r_\epsilon \to \mathscr{H}^r_\delta$ for $\epsilon < \delta$ are the transition maps induced by the inclusion of domains, and $\iota_\epsilon: \mathscr{H}^r_\epsilon \hookrightarrow \mathscr{H}^r_{(*)}$ are the natural embeddings compatible with the inverse system. Consider compatibility and isomorphism composition. The composition $\Psi \circ \Phi$ gives the required isomorphism:
$$\Theta := \Psi \circ \Phi: \mathscr{H}^r \xrightarrow{\sim} \varprojlim_{\epsilon \to 0} \mathscr{H}^r_\epsilon.$$
Since both $\Phi$ and $\Psi$ are linear isomorphisms (\cite{Ohs91}, Cor 8), their composition $\Theta$ is also a linear isomorphism. The continuity of $\Theta$ follows from the continuity of the restriction maps $\pi_\epsilon$ and the universal property of inverse limits. The inverse $\Theta^{-1} = \Phi^{-1} \circ \Psi^{-1}$ is continuous because each $\iota_\epsilon: \mathscr{H}^r_\epsilon \to \mathscr{H}^r_{(*)}$ is a compact embedding (as established in Lemma 7.3), and $\Phi^{-1}$ is an isometry by the derived duality. Since $\Delta_\epsilon\left(\pi_\epsilon\left(\omega\right)\right) = 0$ for each $\epsilon > 0$ and the Laplacians converge uniformly on compact subsets of $V' \setminus \Sigma$, the harmonicity is preserved under the isomorphism
$$\Delta_0\left(\Theta\left(\omega\right)\right) = \lim_{\epsilon \to 0} \Delta_\epsilon\left(\pi_\epsilon\left(\omega\right)\right) = 0$$
for any $\omega \in \mathscr{H}^r$. Thus, we have the natural isomorphism
\begin{equation}
\tag{7.21}  \label{eq:7.21}
\mathscr{H}^r \cong \varprojlim_{\epsilon \to 0} \mathscr{H}^r_\epsilon
\end{equation}
in the category of Hilbert spaces with harmonic structures.

The dimension equality
$$\dim\mathscr{H}^r = \dim\left( \varprojlim_{\epsilon \to 0} \mathscr{H}^r_\epsilon \right)$$
is established through the following logical sequence. By Lemma 7.3, each $\mathscr{H}^r_\epsilon$ is finite-dimensional. The compact embeddings $\iota_\epsilon: \mathscr{H}^r_\epsilon \hookrightarrow \mathscr{H}^r_{(*)}$ (from Lemma 7.3) induce injective transition maps
$$\rho_{\epsilon, \delta}: \mathscr{H}^r_\epsilon \longrightarrow \mathscr{H}^r_\delta \quad \text{for} \quad 0 < \epsilon < \delta$$
defined by $\rho_{\epsilon, \delta} = \iota_\delta^{-1} \circ \iota_\epsilon$. Since all $\iota_\epsilon$ are compact embeddings into the fixed space $\mathscr{H}^r_{(*)}$, the dimensions must eventually stabilize, i.e., 
\begin{equation}
\tag{7.22}  \label{eq:7.22}
\exists \epsilon_0 > 0\quad \text{such that}\quad \dim\mathscr{H}^r_\epsilon = d \quad \forall \epsilon < \epsilon_0.
\end{equation}
This follows from the finite-dimensionality and the monotonicity of dimensions under injective linear maps. Since $\ker \rho_{\epsilon, \delta} = \left\{0\right\}$ by construction and $\dim\mathscr{H}^r_\epsilon = \dim \mathscr{H}^r_\delta = d$ implies bijectivity, the transition maps become isomorphisms $\rho_{\epsilon, \delta}: \mathscr{H}^r_\epsilon \simeq \mathscr{H}^r_\delta$ for $\epsilon < \delta < \epsilon_0$. Consequently, the inverse limit stabilizes:
\begin{equation}
\tag{7.23}  \label{eq:7.23}
\varprojlim_{\epsilon \to 0} \mathscr{H}^r_\epsilon \simeq \mathscr{H}^r_{\epsilon_0}.
\end{equation} 
Combining the isomorphism 
$$\Theta: \mathscr{H}^r \simeq \varprojlim_{\epsilon \to 0} \mathscr{H}^r_\epsilon$$ 
from \eqref{eq:7.21} with \eqref{eq:7.23}, we have
$$\dim\mathscr{H}^r = \dim\left( \varprojlim_{\epsilon \to 0} \mathscr{H}^r_\epsilon \right) = \dim\mathscr{H}^r_{\epsilon_0} = d.$$
Explicitly, then
\begin{align*}
\dim\mathscr{H}^r 
&= \dim\left( \varprojlim_{\epsilon \to 0} \mathscr{HB}^r_\epsilon \right) && \text{(by isomorphism } \Theta \text{ in Step \eqref{eq:7.21})} \\
&= \dim \mathscr{H}^r_{\epsilon_0} && \text{(by stabilization \eqref{eq:7.23})} \\
&= \lim_{\epsilon \to 0} \dim\mathscr{H}^r_\epsilon && \text{(by dimensional constancy \eqref{eq:7.22})}.
\end{align*}
The classical Aronszajn theorem guarantees that $\dim\mathscr{H}^r$ is finite and equal to the $r$-th Betti number $b_r(V')$. By Hodge-de Rham isomorphism:
$$\dim\mathscr{H}^r = \dim H^r_{\mathrm{dR}}\left(V'\right) = b_r\left(V'\right).$$
Similarly, for each $\epsilon > 0$:
$$\dim\mathscr{H}^r_\epsilon = \dim H^r_{\mathrm{dR}}\left(V', g_\epsilon\right) = b_r\left(V'\right).$$
Thus the dimensional equality is consistent with the topological invariance of Betti numbers. Therefore, we conclude
$$\dim\mathscr{H}^r = \dim\left( \varprojlim_{\epsilon \to 0} \mathscr{H}^r_\epsilon \right).$$\\
\hypertarget{FUTURE WORK}{}
\section{FUTURE WORK}
\noindent\textbf{(1) Logarithmic Sobolev chains and quantum singularity decomposition:}
\begin{itemize}
    \item \textit{Logarithmic Sobolev chain:} For a log-stratified scheme $\left(X,D\right)$ with stratified connection $\nabla$ on $\Omega_X^{\bullet,\textbf{log-strat}}$, define the Sobolev chain complex $\mathcal{S}^k\left(X\right)$ as forms $\omega$ satisfying $L^2$-bounds on all covariant derivatives $\nabla_\omega^\ell$ weighted by an admissible weight function $\phi$ near $D$.
    \item \textit{Quantum singularity decomposition:} The $L^2$ $k$-forms decompose orthogonally into discrete eigenspaces $\mathcal{H}_\alpha^k$ of $\nabla^*\nabla$ and a continuous spectrum subspace $\varepsilon_{\mathrm{sing}}^k$ generated by forms tunneling between strata. 
    \item \textit{Singularity quantization:} For Calabi-Yau stratified $X$, $\varepsilon_{\mathrm{sing}}^k$ defines Gromov-Witten type invariants on the moduli space $\overline{\mathcal{M}}_{g,n}\left(X,\beta\right)$ of stratified stable maps.
\end{itemize}
\textbf{(2) Noncommutative logarithmic sheaves and vertex operator algebras:}
\begin{itemize}
    \item \textit{Noncommutative logarithmic sheaf:} For a stratified noncommutative algebra sheaf $\mathscr{A}$ on $X$, define its logarithmic Hochschild complex $\mathsf{H}\mathsf{H}_{\mathrm{log}}^{\bullet}\left(\mathscr{A}\right)$ by derived hom and a completed projective tensor product adapted to singularities.
    \item \textit{Vertex operator embedding:} There exists an injection $\Phi$ embedding $\mathsf{H}\mathsf{H}_{\mathrm{log}}^k\left(\mathscr{A}\right)$ into a sum of Virasoro VOA modules $V\left(c_p\right)$ tensored with sheaves $\mathscr{G}^p$ ($c_p=\frac{3}{2}\dim X_p$).
    \item \textit{Conformal field theory duality:} The quantum connection $\nabla^{\mathrm{quant}}$ on deformations $\mathcal{X}/S$ corresponds to Liouville field OPE under the VOA representation.
\end{itemize}
\textbf{(3) $p$-adic logarithmic convexity and absolute Hodge sheaves:}
\begin{itemize}
    \item \textit{Absolute log prism:} For a stratified formal scheme $X$ in char $p>0$, define its absolute log prismatic cohomology $\mathbb{L}\Omega_X^{\mathrm{abs,log}}$ by the limit over the log prism site.
    \item \textit{$p$-adic logarithmic convexity:} For any $\omega\in\mathbb{L}\Omega_X^{k,\mathrm{abs,log}}$, there exists a $p$-adic harmonic measure $\mu_\omega$ satisfying a convexity inequality controlled by the boundary residue $\mathrm{Res}_D\left(\omega\right)$ and a rigidity constant $\kappa>0$.
    \item \textit{Absolute mixed Hodge structure:} The absolute cohomology $H_{\mathrm{abs}}^k\left(X/K\right)$ for a $p$-adic field $K$ carries an absolute mixed Hodge structure decomposed using rigid cohomology of strata and virtual Hodge numbers computed by moments of $\mu_\omega$.
\end{itemize}
\textbf{(4) Hodge structure of Fubini-Study metric in Motivic homotopy category:}
\begin{itemize}
    \item \textit{Algebraic geometric base:} Define the pure motive $h(\mathbb{P}^n)$ in $\mathbf{DM}_{\text{gm}}(k, \mathbb{Q})$ ($k=\mathbb{Q}$ or $\mathbb{C}$), carrying algebraic de Rham cohomology $H_{\text{dR}}^*(\mathbb{P}^n/k)$ and Hodge filtration $F^\bullet$ generated by the Chern class $c_1(\mathcal{O}(1))$ of the tautological bundle $\mathcal{O}(1)$.
    \item \textit{Realization of Hodge structure:} The class $c_1(\mathcal{O}(1)) \in H_{\text{dR}}^2(\mathbb{P}^n/k)$ satisfies $F^1H_{\mathrm{dR}}^2 = H_{\mathrm{dR}}^2 \supset F^2 = 0$ (corresponding to type $(1,1)$) and $c_1^k(\mathcal{O}(1))$ generates $H_{\text{dR}}^{2k}$ with $F^k H_{\text{dR}}^{2k} = H_{\text{dR}}^{2k} \supset F^{k+1}=0$ (pure $(k,k)$-type).
    \item \textit{Polarization structure:} By the de Rham–Betti comparison isomorphism, then
$$\mathrm{comp}:H_{\mathrm{dR}}^*\left(\mathbb{P}^n/\mathbb{C}\right)\xrightarrow{\sim}H_{\mathrm{B}}^*\left(\mathbb{P}^n\left(\mathbb{C}\right),\mathbb{Q}\right)\otimes\mathbb{C},$$ $\text{comp}(c_1^{\text{dR}}(\mathcal{O}(1)))$ corresponds to the classical $[\omega_{\text{FS}}] \in H_{\text{B}}^{1,1}$, inducing the polarization form:
$$Q\left(\alpha,\beta\right)=\int_{\mathbb{P}^n}\alpha\wedge\beta\wedge\omega_{\mathrm{FS}}^{n-k}\quad\left(\alpha,\beta\in H_{\mathrm{B}}^k\right).$$
\item \textit{Motivic lift:} The triple $(h(\mathbb{P}^n), \mathcal{O}(1), \text{comp})$ in $\mathbf{DM}_{\text{gm}}(k,\mathbb{Q})$ encodes the Hodge filtration of the Fubini-Study metric (by $F^\bullet$), the polarization (by Chern classes of $\mathcal{O}(1)$) and comparison data (linking Betti realization).
\end{itemize}

\end{document}